\documentclass{amsart}
\usepackage{amsfonts,amsmath,amssymb,mathrsfs,array}
\usepackage[english]{babel}
\usepackage[colorlinks=true]{hyperref}
\newtheorem{thm}{Theorem}[section]
\newtheorem{lem}[thm]{Lemma}
\newtheorem{prop}[thm]{Proposition}
\newtheorem{proposition}[thm]{Proposition}
\newtheorem{cor}[thm]{Corollary}
\newtheorem{corollary}[thm]{Corollary}
\theoremstyle{definition}
\newtheorem{definition}[thm]{Definition}
\newtheorem{cond}[thm]{Condition}
\newtheorem{rem}[thm]{Remark}
\newtheorem{remark}[thm]{Remark}
\numberwithin{equation}{section}
\newcommand{\C}{\mathcal{C}}

\begin{document}
\title[Convolution on quantum tori II]
{Convolution on quantum tori II: Noncommutative version of  O'Neil inequality and its applications to fractional Sobolev embeddings}
\author[F. Sukochev]{F. Sukochev}
\address{School of Mathematics and Statistics, University of New South Wales, Kensington, NSW 2052, Australia}
\email{f.sukochev@unsw.edu.au}

\author[R. Tastankul]{R. Tastankul}
\address{Institute of Mathematics and Mathematical Modeling,
050010 Almaty, Kazakhstan}
\email{ramazan.tastankul@gmail.com}

\author[K. Tulenov]{K. Tulenov}
\address{Institute of Mathematics and Mathematical Modeling,
050010 Almaty, Kazakhstan; and School of Mathematics and Statistics,
University of New South Wales, Kensington, NSW 2052, Australia}
\email{tulenov@math.kz}

\author[D. Zanin]{D. Zanin}
\address{School of Mathematics and Statistics, Central South University,
Changsha 410075, China}
\email{d.zanin@csu.edu.cn}

\date{}

\begin{abstract}
We prove a sharp O'Neil-type rearrangement inequality for operator convolution on quantum tori, with optimal constant one. We apply this inequality to obtain convolution estimates in Lorentz and Besov--Lorentz spaces, including logarithmic limiting estimates. We also study Sobolev and Besov spaces built on symmetric operator spaces and establish fractional Sobolev embeddings and Besov multiplication estimates. Finally, we apply these results to fractional diffusion and prove local well-posedness for semilinear heat equations.
\end{abstract}

\maketitle
\begingroup\small
\tableofcontents
\endgroup

\section{Introduction}

Convolution inequalities are useful for proving Sobolev
embeddings and studying evolution equations. In the classical
setting, O'Neil's rearrangement inequality
\cite[Lemma~2.5]{ONeil} gives convolution estimates in Lorentz
spaces, while the limiting estimates of Br\'ezis and Wainger
\cite[Theorem~2]{BW} describe logarithmic integrability at
critical exponents. We establish corresponding inequalities
on quantum tori and develop their applications to symmetric
operator spaces, Sobolev embeddings, and heat equations.

Quantum tori are basic examples in noncommutative geometry.
Their canonical torus action provides both a Fourier analysis
and a differential structure; see Connes
\cite{Connes1980,Connes1994} and
\cite{ChenXuYin,HaLeePonge,McDSX}.
Sobolev, Besov, and Triebel--Lizorkin spaces based on
noncommutative $L_p$ norms were studied systematically by
Xiong, Xu, and Yin \cite{XiongXuYin}.
Here we study Sobolev and Besov spaces with more general
symmetric operator norms.

This paper continues the study of operator convolution begun
in Part~I \cite{ConvolutionI}. For real skew-symmetric
$d\times d$ matrices $\theta_1,\theta_2$, that construction gives
\begin{equation*}
 \star_{\theta_1,\theta_2}:
 L_1(\mathbb T_{\theta_1}^d)\times
 L_1(\mathbb T_{\theta_2}^d)
 \longrightarrow
 L_1(\mathbb T_{\theta_1+\theta_2}^d),
\end{equation*}
with Fourier coefficients
\begin{equation*}
 \widehat{x\star_{\theta_1,\theta_2}y}(n)
 =\widehat x(n)\widehat y(n),
 \ \ n\in\mathbb Z^d.
\end{equation*}
It recovers classical convolution when
$\theta_1=\theta_2=0$ and scalar semi-convolution when
$\theta_1=0$.
The canonical trace on each quantum torus is normalized
by $\tau_\theta(\mathbf 1)=1$, and Haar measure on
$\mathbb T^d$ has total mass one.
We write $\mu(\cdot,x)$ for the singular value function
of an operator $x$, which in the commutative case is the
decreasing rearrangement of $|x|$; see
\cite{DoddsDePagterSukochev,LSZ,LMSZ}.

Our main result is the sharp inequality
\begin{equation}\label{eq:intro-sharp-ONeil}
\begin{split}
 \int\limits_0^t
 \mu(s,x\star_{\theta_1,\theta_2}y)\,ds
 &\leq
 \int\limits_0^t\mu(s,x)\,ds
 \cdot
 \int\limits_0^t\mu(s,y)\,ds\\
 &+
 t\int\limits_t^1\mu(s,x)\mu(s,y)\,ds,
 \ \ 0<t<1,
\end{split}
\end{equation}
for $x\in L_1(\mathbb T_{\theta_1}^d)$ and
$y\in L_1(\mathbb T_{\theta_2}^d)$; see
Theorem~\ref{thm:operator-oneil}.
Equality holds when $x$ and $y$ are the identity operators
in their respective algebras, so the constant one is optimal.
We prove an abstract bilinear operator-valued version
of O'Neil's scalar argument and then apply the endpoint
convolution bounds from Part~I \cite{ConvolutionI}.

The main advantage of \eqref{eq:intro-sharp-ONeil} is that it controls the singular value distribution of the convolution, rather than only a single $L_p$ norm. This estimate therefore yields convolution inequalities in several function-space scales, including spaces that describe limiting integrability.

These consequences include Young--O'Neil inequalities
in Lorentz spaces
(Corollary~\ref{cor:Young-ONeil}),
logarithmic limiting estimates
(Theorem~\ref{thm:working-Brezis-Wainger}),
and Besov--Lorentz convolution estimates
(Theorem~\ref{thm:Besov-Lorentz}).
Theorem~\ref{thm:working-Brezis-Wainger}(ii)
is a quantum-torus counterpart of the limiting convolution
estimate in \cite[Theorem~2]{BW}.
For other related classical convolution results, see
\cite{Burenkov,KuhnSchilling,Krepela}.
Our logarithmic targets belong to the Lorentz--Zygmund
scale \cite{OpicPick}. They differ from the oscillation
spaces used in \cite{Nursultanov}, whose defining functionals
involve $f^{**}-f^*$, where
$f^*=\mu(f)$ and
$f^{**}(t)=t^{-1}\int_0^t f^*(s)\,ds$.

For the Sobolev and Besov applications, let $E(0,1)$
be a fully symmetric Banach function space.
The associated symmetric operator space is
\begin{equation*}
 \mathcal E(\mathbb T_\theta^d)
 =
 \{x:\mu(x)\in E(0,1)\},
 \ \ \|x\|_{\mathcal E}=\|\mu(x)\|_E,
\end{equation*}
where $x$ ranges over the $\tau_\theta$-measurable operators
affiliated with $\mathbb T_\theta^d$.
Thus $E=L_p(0,1)$ gives
$\mathcal E=L_p(\mathbb T_\theta^d)$.
We define
\begin{equation*}
 H_{\mathcal E}^s(\mathbb T_\theta^d)
 =
 \left\{
 x\in\mathcal S'(\mathbb T_\theta^d):
 J_\theta^s x\in\mathcal E(\mathbb T_\theta^d)
 \right\},
 \ \ J_\theta^s=(1-\Delta_\theta)^{\frac s2},
\end{equation*}
with norm
$\|x\|_{H_{\mathcal E}^s}=\|J_\theta^s x\|_{\mathcal E}$.
Proposition~\ref{prop:prelim-completeness}
establishes completeness and the Bessel lifting property.
When $E$ is separable,
Theorem~\ref{thm:HEs-properties} also gives density
of Fourier polynomials.

Potential-kernel estimates and
\eqref{eq:intro-sharp-ONeil} yield subcritical Lorentz
and critical logarithmic Sobolev embeddings; see
Theorems~\ref{thm:Lorentz-Sobolev-embedding}
and~\ref{thm:critical-logarithmic-Sobolev}.
More generally, at the critical order $\frac d2$,
we use the scalar Hardy operator
\begin{equation}\label{eq:intro-Hardy}
 (\mathcal Th)(t)
 =
 t^{-\frac12}\int\limits_0^t h(v)\,dv
 +
 \int\limits_t^1v^{-\frac12}h(v)\,dv,
 \ \ 0<t<1,\ \ h\in L_1(0,1).
\end{equation}
Theorem~\ref{thm:critical-fractional-distribution}
controls the singular values of
$J_\theta^{-\frac{d}{2}}x$ through $\mathcal T\mu(x)$.
Consequently, if $\mathcal T:E\to F$ is bounded and
$F$ is fully symmetric, then
\begin{equation*}
 H_{\mathcal E}^{\frac{d}{2}}(\mathbb T_\theta^d)
 \hookrightarrow
 \mathcal F(\mathbb T_\theta^d),
\end{equation*}
where $\mathcal F$ is the operator space associated
with $F$; see
Theorem~\ref{thm:fractional-symmetric-embedding}.

This gives a scalar sufficient condition for a
noncommutative Sobolev embedding.
Different target spaces can therefore be obtained
by studying the same one-dimensional operator,
without repeating the noncommutative argument.

Classical reductions of Sobolev embeddings to
one-dimensional inequalities can be found in
\cite{Klimov1969,Klimov1970,EKP,Cianchi,CPS}.

In particular, Corollary~\ref{cor:critical-L2-Lambda}(ii),
using the scalar estimate in \cite[Lemma~9]{SZ}, gives
\begin{equation*}
 H_2^{\frac{d}{2}}(\mathbb T_\theta^d)
 \hookrightarrow
 \Lambda_\psi^{(2)}(\mathbb T_\theta^d),
 \ \ \psi(t)\asymp\frac1{\log(e/t)}.
\end{equation*}
Here $\psi$ is the concave gauge specified in the
preliminaries, and $\Lambda_\psi^{(2)}$ is the
$2$-convexification of the Lorentz space $\Lambda_\psi$.
Its norm satisfies
\begin{equation*}
 \|x\|_{\Lambda_\psi^{(2)}}
 \asymp
 \left(
 \int\limits_0^1
 \left[\frac{\mu(t,x)}{\log(e/t)}\right]^2
 \frac{dt}{t}
 \right)^{\frac12}.
\end{equation*}
This is a noncommutative counterpart of the classical
critical embedding of Hansson and Br\'ezis--Wainger;
see \cite[Theorem~3.3 and equation~(3.13)]{Hansson}
and \cite[Theorem~2]{BW}.
The explicit logarithmic target describes integrability
at the critical order, where an embedding into $L_\infty$
is unavailable in general.

For the Besov scale, we define
\begin{equation*}
 \|x\|_{B_{\mathcal E,r}^s(\mathbb T_\theta^d)}
 =
 \left\|
 \left(2^{js}\|\Delta_j^\theta x\|_{\mathcal E}\right)_{j\geq0}
 \right\|_{\ell_r},
\end{equation*}
where $(\Delta_j^\theta)_{j\geq0}$ is a smooth dyadic
decomposition of the Fourier frequencies.
Proposition~\ref{lem:besov-app-basic} establishes
completeness, independence of the chosen decomposition,
and density of Fourier polynomials when $r<\infty$,
without assuming separability of $E$.

For $s>0$, Lemma~\ref{lem:symmetric-besov-embedding}
gives a sufficient condition on $E$, $s$, and $r$ for
$B_{\mathcal E,r}^s(\mathbb T_\theta^d)$ to embed
into $L_\infty(\mathbb T_\theta^d)$.
Under this condition,
Proposition~\ref{prop:besov-app-algebra} shows that
$B_{\mathcal E,r}^s(\mathbb T_\theta^d)$ is also a
Banach algebra under operator multiplication.
For $E=L_p$, this includes $s>\frac{d}{p}$ and the endpoint
$s=\frac{d}{p}$, $r=1$, for $1\leq p<\infty$.
The product estimates follow by adapting Bony's
decomposition to preserve the order of multiplication;
see \cite{Bony,SchmeisserTriebel,BahouriCheminDanchin}
for the classical theory.

Finally, we consider linear fractional diffusion
and semilinear heat equations
\begin{equation*}
 \partial_tu-\Delta_\theta u=\mathcal P(u,u^*),
 \ \ u(0)=u_0.
\end{equation*}
Here $u=u(t)$ takes values in
$B_{\mathcal E,r}^s(\mathbb T_\theta^d)$,
$u^*$ is its adjoint, $\partial_t$ denotes differentiation
in time, and $\mathcal P$ is a polynomial with fixed
coefficients and a prescribed order of multiplication.
Theorem~\ref{thm:besov-app-algebra-pde}
gives local well-posedness and a blow-up alternative
in the Besov algebra range, including the Sobolev
setting of Yang--Yang \cite[Theorem~1.2]{YY}.

Corollary~\ref{cor:besov-app-power} also allows less
regular initial data. For example, cubic equations
with smooth coefficients in dimension three admit
initial data in $H_\theta^1$.
Such data need not be bounded operators.
The smoothing of the heat semigroup nevertheless
allows the nonlinear equation to be solved, extending
the local theory beyond the Banach algebra setting.

Below the algebra threshold, uniqueness holds in the
solution class of Theorem~\ref{thm:besov-app-weighted}.
For smooth coefficients,
Corollary~\ref{cor:besov-app-positive-regularity}
gives spatial smoothness at every positive time.
These nonlinear arguments use the endpoint convolution
bounds from Part~I and the Besov product estimates;
the full rearrangement inequality is used for the
refined convolution and Sobolev estimates.

Section~\ref{sec:preliminaries} gives the notation
and background.
Section~\ref{sec:convolution} proves the main
rearrangement inequality.
Section~\ref{sec:potential-structure} develops the
Sobolev and Besov spaces.
Sections~\ref{sec:applications} and
\ref{sec:sobolev-applications} contain the convolution
and Sobolev applications, respectively.
Section~\ref{sec:evolution} treats the evolution equations.

\section{Preliminaries}\label{sec:preliminaries}

We collect the notation and background on symmetric spaces, quantum tori, and Fourier multipliers used throughout the paper.

All finite traces in this paper are normalized by $\tau(\mathbf{1})=1$. The classical torus is $\mathbb{T}^{d}=\mathbb{R}^{d}/(2\pi\mathbb{Z})^{d}$, with normalized Haar measure.

\subsection{Symmetric (quasi-)Banach function spaces}

Let $(\Omega,m)$ be a nonatomic probability space, and let $S(\Omega)$ denote the measurable complex-valued functions modulo equality almost everywhere. For $f\in S(\Omega)$ its distribution and decreasing rearrangement are
\begin{equation*}
d_f(\lambda)=m\left\{\omega:|f(\omega)|>\lambda\right\},\ \ 
\mu(t,f)=\inf\left\{\lambda\geq0:d_f(\lambda)\leq t\right\}, \ \  0<t<1.
\end{equation*}
A symmetric (quasi-)Banach function space $E(\Omega)$ is a complete (quasi-)normed linear subspace of $S(\Omega)$ such that $\mu(g)\leq\mu(f)$ and $f\in E$ imply $g\in E$ and $\|g\|_E\leq\|f\|_E$. We use the corresponding representation on $(0,1)$ and the same notation $E$ for it; see \cite{BennettSharpley,KPS,LT}.

For $0<p<\infty$ and $0<q\leq\infty$, the Lorentz quasi-norm is
\begin{equation}\label{eq:prelim-Lorentz}
\|f\|_{L_{p,q}}=
\begin{cases}
\displaystyle\left(\int\limits_{0}^{1}[t^{\frac{1}{p}}\mu(t,f)]^q\frac{dt}{t}\right)^{\frac{1}{q}},
&q<\infty,\\
\displaystyle\sup_{0<t<1}t^{\frac{1}{p}}\mu(t,f),&q=\infty.
\end{cases}
\end{equation}
In particular, $L_{p,p}=L_p$ with equality of these norms, and
\begin{equation}\label{LorentzEqvNorm}
\|f\|_{L_{p,\infty}}=\sup_{\lambda>0}\lambda d_f(\lambda)^{\frac{1}{p}}.
\end{equation}
The supremum here is over all positive levels $\lambda$. The Hardy--Ces\`aro operator is denoted by
\begin{equation}\label{eq:Cesaro-prelim}
(\mathcal{C}h)(t)=\frac{1}{t}\int\limits_0^t h(s)\,ds.
\end{equation}
For $1<p<\infty$ and $1\leq q\leq\infty$,
\begin{equation}\label{LpqNormsIneq}
\|f\|_{L_{p,q}}\leq
\|t^{\frac{1}{p}}\mathcal{C}\mu(f)(t)\|_{L_q((0,1),\frac{dt}{t})}
\leq\frac{p}{p-1}\|f\|_{L_{p,q}}.
\end{equation}
The middle expression is a fully symmetric Banach norm equivalent to \eqref{eq:prelim-Lorentz}. We use the usual norm for $L_{1,1}=L_{1}$. For $q>1$, $L_{1,q}$ need not be contained in $L_{1}$.

We say that $g$ is Hardy--Littlewood--P\'olya submajorized by $f$, and write $g\prec\!\prec f$, if
\begin{equation}\label{eq:prelim-HLP}
\int\limits_0^t\mu(s,g)\,ds\leq\int\limits_0^t\mu(s,f)\,ds,
\ \  0<t\leq1.
\end{equation}
A symmetric space is \emph{fully symmetric} if this relation, with $f\in E$, implies $g\in E$ and $\|g\|_E\leq\|f\|_E$. Every symmetric Banach function space on $(0,1)$ satisfies $L_\infty\hookrightarrow E\hookrightarrow L_{1}$.

Let $\eta:[0,1]\to[0,\infty)$ be increasing and concave, with $\eta(0)=0$ and $\eta(t)>0$ for $t>0$. The associated Lorentz and
Marcinkiewicz spaces have norms
\begin{equation*}
\|f\|_{\Lambda_\eta}=\int\limits_{0}^{1}\mu(t,f)\,d\eta(t),\ \ 
\|f\|_{\mathtt M_\eta}=\sup_{0<t\leq1}
\frac1{\eta(t)}\int\limits_0^t\mu(s,f)\,ds.
\end{equation*}
They are fully symmetric Banach spaces and $(\Lambda_\eta)'=\mathtt M_\eta$ for the K\"othe associate norm. For $r\geq1$, the $r$-convexification is
\begin{equation*}
E^{(r)}=\left\{f:|f|^r\in E\right\},\ \ 
\|f\|_{E^{(r)}}=\||f|^r\|_E^{\frac{1}{r}}.
\end{equation*}
In particular, $\|f\|_{E^{(2)}}^2=\||f|^2\|_E$. We refer to \cite[Chapter~6]{DoddsDePagterSukochev} for these constructions.

\subsection{The Hardy-type operator}

Following \cite{SZ}, for $f\in S(0,1)$ we define
\begin{equation}\label{operatorT}
(\mathcal{T}f)(t)=\frac1{\sqrt t}\int\limits_0^t f(s)\,ds
       +\int\limits_t^1\frac{f(s)}{\sqrt s}\,ds,\ \  0<t<1.
\end{equation}
If $f\in L_{1}(0,1)$ both terms are finite for every $t>0$. Its kernel is $(\max\left\{s,t\right\})^{-\frac{1}{2}}$, whose square has integral $2$ over $(0,1)^2$; hence $T$ is Hilbert--Schmidt on $L_2(0,1)$.

Put $\ell(t)=\log\left( e/t \right)$, and use the increasing concave gauge
\begin{equation}\label{PhiPsiFunc}
\psi(t)=
\begin{cases}
\ell(t)^{-1},&0<t\leq e^{-1},\\
\frac{1}{2}+\frac e4(t-e^{-1}),&e^{-1}<t\leq1,
\end{cases}
\ \  \psi(0)=0.
\end{equation}
This is equivalent to $\ell(t)^{-1}$ on $(0,1)$. The latter function has second derivative $(2-\ell(t))/(t^2\ell(t)^3)$, so the linear extension in \eqref{PhiPsiFunc} supplies a concave gauge on the whole interval. It defines the same spaces, up to equivalent norms, as the logarithmic gauge used in \cite{SZ}. In particular,
\begin{equation}\label{eq:fractional-Marcinkiewicz}
\|f\|_{\mathtt M_\psi} = \sup_{0<t\leq1}\frac1{\psi(t)}\int\limits_0^t\mu(s,f)\,ds.
\end{equation}

\begin{prop}\label{ThmTBdd}
The operator $\mathcal{T}$ satisfies:
\begin{itemize}
\item[(i)] $\|\mathcal{T}\|_{L_p(0,1)\to L_p(0,1)}\leq2$ for $1\leq p\leq\infty$;
\item[(ii)] $\mathcal{T}:L_2(0,1)\to\Lambda_\psi^{(2)}(0,1)$ is bounded.
\end{itemize}
\end{prop}
\begin{proof}
The positive kernel $K(t,s)=(\max\left\{t,s\right\})^{-\frac{1}{2}}$ satisfies
\begin{equation*}
 \int\limits_{0}^{1}K(t,s)\,dt=2-\sqrt s\leq2,\ \ 
 \int\limits_{0}^{1}K(t,s)\,ds=2-\sqrt t\leq2.
\end{equation*}
Tonelli's theorem gives the $L_1$ bound, and the second identity gives the $L_\infty$ bound. Interpolation proves (i) \cite[Theorem~1.3.4]{G2014}. Part (ii) is \cite[Lemma~9]{SZ}, since the concave gauge in \eqref{PhiPsiFunc} is equivalent to $\log^{-1}(e/t)$.
\end{proof}

\subsection{Noncommutative symmetric (quasi-)Banach spaces}
\label{SubsNCOpSp}

Let $\mathcal M\subset B(\mathcal H)$ be a finite von Neumann algebra, that is, a unital $*$-subalgebra closed in the weak operator topology, and let $\tau$ be a normal faithful tracial state. Its commutant is
\begin{equation*}
\mathcal M'=\left\{y\in B(\mathcal H):yx=xy\text{ for every }x\in\mathcal M\right\}.
\end{equation*}
A closed densely defined operator $x$ is \emph{affiliated} with $\mathcal M$ if $u^\ast xu=x$ for every unitary $u\in\mathcal M'$, with equality of domains. It is \emph{$\tau$-measurable} if, for every $\varepsilon>0$, there is a projection $p\in\mathcal M$ such that $p\mathcal H\subset\operatorname{dom}(x)$ and $\tau(\mathbf{1}-p)<\varepsilon$. Denote these operators by $S(\mathcal M,\tau)$. Products and sums mean the closures of the corresponding densely defined operations.

For $x\in S(\mathcal M,\tau)$, define
\begin{equation}\label{eq:prelim-mu}
 d_x(\lambda)=\tau(\mathbf{1}_{(\lambda,\infty)}(|x|)),\ \ 
 \mu(t,x)=\inf\left\{\lambda\geq0:d_x(\lambda)\leq t\right\}, \ \  0<t<1.
\end{equation}
This singular-value function is decreasing and right-continuous. The noncommutative submajorization relation is again \eqref{eq:prelim-HLP}. See \cite{FackKosaki,DoddsDePagterSukochev,LSZ}.

For a symmetric function space $E(0,1)$, set
\begin{equation}\label{eq:operator-space-norm}
\mathcal{E}(\mathcal M)=\left\{x\in S(\mathcal M,\tau):\mu(x)\in E\right\},
\ \  \|x\|_{\mathcal{E}(\mathcal M)}=\|\mu(x)\|_E.
\end{equation}
This is the associated symmetric operator space; the usual completeness and symmetry properties are discussed in \cite{KS,F,DoddsDePagterSukochev}. 
When $E$ is fully symmetric, so is $\mathcal{E}(\mathcal M)$. For $\mathcal M=L_\infty(0,1)$ this construction gives the original function space. In the general semifinite theory the rearrangement interval is $(0,\tau(\mathbf{1}))$ rather than necessarily $(0,1)$; see \cite{PXu,LSZ,DoddsDePagterSukochev}.

In particular, for $0<p<\infty$,
\begin{equation*}
L_p(\mathcal M)=\left\{x\in S(\mathcal M,\tau):\tau(|x|^p)<\infty\right\},
\ \  \|x\|_p=\tau(|x|^p)^{\frac{1}{p}}.
\end{equation*}
We set $L_\infty(\mathcal M)=\mathcal M$ with its operator norm. For $p\geq1$ these are Banach spaces; for $0<p<1$ they are quasi-Banach spaces. We take the inner product on $L_2$ to be $\langle x,y\rangle=\tau(x^*y)$, linear in its second variable. The spaces $L_{p,q}(\mathcal M)$, $\Lambda_\eta(\mathcal M)$ and $\mathtt M_\eta(\mathcal M)$ are defined by applying their function norms to $\mu(x)$. In particular, noncommutative K\"othe duality gives
\begin{equation*}
\|xy\|_1\leq\|x\|_{\mathtt M_\eta}\|y\|_{\Lambda_\eta}.
\end{equation*}
For the compatible couple $(L_{1}(\mathcal M),L_\infty(\mathcal M))$,
\begin{equation}\label{eq:K-prelim}
K(t,x;L_{1},L_\infty)=\int\limits_0^t\mu(s,x)\,ds,\ \  0<t\leq1.
\end{equation}
Consequently, a linear operator of norm at most $C$ at both endpoints sends $x$ to an operator submajorized by $C\mu(x)$.

\subsection{Noncommutative tori}

Let $d\geq2$ and let $\theta=(\theta_{jk})$ be a real skew-symmetric $d\times d$ matrix. The quantum torus $C^*$-algebra $\mathcal A_\theta$ is generated by unitaries $U_1,\ldots,U_d$ satisfying
\begin{equation*}
U_{j}U_k=e^{2\pi i\theta_{jk}}U_kU_{j}.
\end{equation*}
To agree with Part I, use the symmetrized Fourier unitaries
\begin{equation*}
e_n^{\theta} = \exp\left(-\pi i\sum_{j<k}\theta_{jk}n_{j}n_k\right)
U_1^{n_1}\cdots U_d^{n_d},\ \  n\in\mathbb{Z}^{d}.
\end{equation*}
They satisfy
\begin{equation}\label{eq:prelim-Weyl}
e_n^{\theta} e_m^{\theta}=e^{\pi i\langle n,\theta m\rangle}e_{n+m}^{\theta}, \ \  (e_n^{\theta})^*=e_{-n}^{\theta}, \ \ n\in\mathbb{Z}^{d}.
\end{equation}
The finite linear combinations form the dense polynomial algebra $\mathcal{P}_\theta\subset\mathcal A_\theta$. The functional
\begin{equation*}
\tau_\theta\left(\sum_n a_ne_n^{\theta}\right)=a_{\textbf{0}}, \ \ \textbf{0}=(0,...,0)\in\mathbb{Z}^{d}
\end{equation*}
extends to a faithful tracial state on $\mathcal A_\theta$. Let $\pi_{\tau_\theta}$ be its GNS representation. We write
\begin{equation*}
L_\infty(\mathbb{T}_{\theta}^{d}) = \pi_{\tau_\theta}(\mathcal A_\theta)''.
\end{equation*}
The trace extends to a normal faithful tracial state on this finite von Neumann algebra. The symbol $\mathbb{T}_{\theta}^{d}$ denotes the noncommutative space; $\mathcal A_\theta$ and $L_\infty(\mathbb{T}_{\theta}^{d})$ denote its $C^*$- and von Neumann algebras, respectively. See \cite{ChenXuYin,XiongXuYin,HaLeePonge,LSZ}. We abbreviate $S(L_\infty(\mathbb{T}_{\theta}^{d}),\tau_\theta)$ to $S(\mathbb{T}_{\theta}^{d})$ and use the corresponding notation for its symmetric operator spaces.

Each $e_n^{\theta}$ is unitary, and
\begin{equation*}
\tau_\theta((e_n^{\theta})^*e_m^{\theta})=\delta_{n,m}.
\end{equation*}
Thus these unitaries form an orthonormal basis of $L_2(\mathbb{T}_{\theta}^{d})$. For $x\in L_{1}(\mathbb{T}_{\theta}^{d})$,
\begin{equation*}
\widehat x(n)=\tau_\theta(x(e_n^{\theta})^*),\ \ 
x\sim\sum_{n\in\mathbb{Z}^{d}}\widehat x(n)e_n^{\theta}.
\end{equation*}
For $x\in L_2$, the series converges in $L_2$ and Parseval's identity is
\begin{equation}\label{Parceveal}
\|x\|_2^2=\sum_{n\in\mathbb{Z}^{d}}|\widehat x(n)|^2.
\end{equation}
Since $\tau_\theta(\mathbf{1})=1$, we have $L_q(\mathbb{T}_{\theta}^{d})\hookrightarrow L_p(\mathbb{T}_{\theta}^{d})$ for $0<p\leq q\leq\infty$. When $\theta=0$, $e_n(t)=e^{in\cdot t}$ and these constructions recover the classical torus.

The smooth algebra $C^\infty(\mathbb{T}_{\theta}^{d})$ consists of the Fourier series with rapidly decreasing coefficients:
\begin{equation*}
\sup_{n\in\mathbb{Z}^{d}}(1+|n|)^N|\widehat x(n)|<\infty
 \ \ \text{for every }N\in\mathbb N_0.
\end{equation*}
Its natural Fr\'echet topology is defined by these seminorms. Its continuous linear dual is denoted by $\mathcal S'(\mathbb{T}_{\theta}^{d})$. We use the bilinear pairing $(F,a)=F(a)$, with an integrable $x$ acting by $(x,a)=\tau_\theta(xa)$. A distribution is uniquely determined by
\begin{equation*}
\widehat F(n)=(F,(e_n^{\theta})^*).
\end{equation*}
Its coefficients have polynomial growth, and conversely every such sequence defines a distribution. The Fourier series converges in the distribution topology, for example through its finite symmetric partial sums; see \cite{XiongXuYin,McDSX,SchmeisserTriebel}.

\subsection{Fourier multipliers and semi-convolution}

For a symbol $\phi:\mathbb{Z}^{d}\to\mathbb{C}$, the Fourier multiplier on polynomials is defined by
\begin{equation*}
T_{\phi}^{\theta} e_n^{\theta}=\phi(n)e_n^{\theta}, \ \ n\in\mathbb{Z}^{d}.
\end{equation*}
If $\phi$ has polynomial growth, this defines continuous maps on $C^\infty(\mathbb{T}_{\theta}^{d})$ and $\mathcal S'(\mathbb{T}_{\theta}^{d})$. When it extends boundedly between given operator spaces, we call it a Fourier multiplier between those spaces. For bounded $\phi$ it acts boundedly on $L_2$, with norm $\sup_n|\varphi(n)|$. We use $T_{\phi}^{\theta}$ throughout; see \cite{ChenXuYin,Ricard, RST,STT, XiongXuYin}.

The canonical torus action is the trace-preserving $*$-automorphism
\begin{equation*}
\alpha_s^{\theta}(e_n^{\theta})=e^{in\cdot s}e_n^{\theta},
\ \  s\in\mathbb{T}^{d}.
\end{equation*}
It preserves singular values and hence every symmetric operator norm.
For $f\in L_{1}(\mathbb{T}^{d})$ and $x\in L_{1}(\mathbb{T}_{\theta}^{d})$,
define the \emph{semi-convolution}
\begin{equation}\label{eq:prelim-semiconvolution}
f\ast x=\int\limits_{\mathbb{T}^{d}}f(s)\alpha_{-s}^{\theta}(x)\,dm(s).
\end{equation}
The action is continuous in $L_{1}$, so this is a Bochner integral in
$L_{1}(\mathbb{T}_{\theta}^{d})$. For bounded $x$ it also agrees with the
weak-$*$ integral. The sign in the action gives
\begin{equation*}
\widehat{f\ast x}(n)=\widehat f(n)\widehat x(n),\ \ 
f\ast x=T_{\widehat f}^{\theta} x.
\end{equation*}
It is the case $\theta_1=0$ of the convolution recalled in
Section~\ref{sec:convolution}.

\begin{lem}[Endpoint bounds for semi-convolution]\label{FirstTechLemma}
Let $f\in L_{1}(\mathbb{T}^{d})$. Then
\begin{align}
\|f\ast x\|_{L_{1}(\mathbb{T}_{\theta}^{d})}
 &\leq\|f\|_{L_{1}(\mathbb{T}^{d})}\|x\|_{L_{1}(\mathbb{T}_{\theta}^{d})},
 \ \ x\in L_{1}(\mathbb{T}_{\theta}^{d}),\label{ConvIneqFirst}\\
\|f\ast x\|_{L_\infty(\mathbb{T}_{\theta}^{d})}
 &\leq\|f\|_{L_{1}(\mathbb{T}^{d})}\|x\|_{L_\infty(\mathbb{T}_{\theta}^{d})},
\ \   x\in L_\infty(\mathbb{T}_{\theta}^{d}).\label{ConvIneqSecond}
\end{align}
\end{lem}
\begin{proof}
The action preserves the trace and absolute values; hence
$\|\alpha_{-s}^{\theta}(x)\|_p=\|x\|_p$. For $x\in L_{1}$,
the Bochner integral gives
\begin{align*}
\|f\ast x\|_1
&=\left\|\int\limits_{\mathbb{T}^{d}}f(s)\alpha_{-s}^{\theta}(x)\,dm(s)\right\|_1\\
&\leq\int\limits_{\mathbb{T}^{d}}|f(s)|
       \tau_\theta(|\alpha_{-s}^{\theta}(x)|)\,dm(s)\\
&=\int\limits_{\mathbb{T}^{d}}|f(s)|\|x\|_1\,dm(s)
 =\|f\|_1\|x\|_1.
\end{align*}
For $x\in L_\infty$, use the weak-$*$ integral. For every
$a\in L_{1}(\mathbb{T}_{\theta}^{d})$,
\begin{align*}
|\tau_\theta((f\ast x)a)|
&\leq\int\limits_{\mathbb{T}^{d}}|f(s)|
       |\tau_\theta(\alpha_{-s}^{\theta}(x)a)|\,dm(s)\\
&\leq\int\limits_{\mathbb{T}^{d}}|f(s)|
       \|\alpha_{-s}^{\theta}(x)\|_\infty\|a\|_1\,dm(s)\\
&=\|f\|_1\|x\|_\infty\|a\|_1.
\end{align*}
Taking the supremum over $\|a\|_1\leq1$ proves the second estimate.
The two constructions agree on bounded $x$, since their pairings
with bounded test operators coincide.
\end{proof}

\begin{lem}[Submajorization for semi-convolution]\label{first oneal lemma}
For $f\in L_{1}(\mathbb{T}^{d})$ and $x\in L_{1}(\mathbb{T}_{\theta}^{d})$,
\begin{equation*}
\mu(f\ast x)\prec\!\prec\|f\|_1\mu(x).
\end{equation*}
\end{lem}
\begin{proof}
The case $f=0$ is immediate. Otherwise divide by $\|f\|_1$ and
assume $\|f\|_1=1$. Lemma~\ref{FirstTechLemma} shows that
$x\mapsto f\ast x$ is a contraction on both $L_{1}$ and $L_\infty$.
The interpolation argument in \cite[Theorem~7.12.6]{DoddsDePagterSukochev}
therefore applies. More explicitly, \eqref{eq:K-prelim} gives
\begin{equation*}
\int\limits_0^t\mu(s,f\ast x)\,ds
=K(t,f\ast x;L_{1},L_\infty)
\leq K(t,x;L_{1},L_\infty)
=\int\limits_0^t\mu(s,x)\,ds.
\end{equation*}
Rescaling $f$ proves the claim.
\end{proof}

\subsection{Logarithmic limiting classes}
\label{subsec:logarithmic-limiting-classes}

Put $\ell(t)=\log\left( e/t \right)$ for $0<t<1$. For $1<q<\infty$, define
\begin{equation}\label{eq:prelim-log-class}
 \|x\|_{L_{\infty,q}}
 =\left(
 \int\limits_{0}^{1}
 \left[\frac{\mu(t,x)}{\ell(t)}\right]^q
 \frac{dt}{t}
 \right)^{\frac{1}{q}}.
\end{equation}
In the commutative setting, this is the classical Lorentz--Zygmund space $L^{\infty,q}(\log L)^{-1}$, also denoted by $L_{\infty,q;-1}$; see \cite{OpicPick}. These spaces arise naturally as target spaces in limiting Sobolev embeddings; see \cite{Hansson,BW}.

The endpoint conventions are $L_{\infty,1}=L_\infty$ and $L_{\infty,\infty}=\operatorname{Exp}L$. For the latter, we use the rearrangement quasi-norm
\begin{equation*}
 \|x\|_{L_{\infty,\infty}}
 =\sup_{0<t<1}\frac{\mu(t,x)}{\ell(t)},
\end{equation*}
which is equivalent to the Luxemburg norm associated with $\Phi(u)=e^u-1$. The convention at $q=1$ is a separate definition: substituting $q=1$ into \eqref{eq:prelim-log-class} gives an infinite value for every nonzero element.

For $1<q\leq\infty$, these spaces admit equivalent fully symmetric Banach norms; see \cite[Theorem~7.4]{OpicPick} and, for the associated operator spaces, \cite{DoddsDePagterSukochev}. More explicitly, one may use
\begin{equation*}
 \|x\|_{L_{\infty,q}}^\# = \left\|t\mapsto
 \frac{(\mathcal{C}\mu(x))(t)}{\ell(t)}
 \right\|_{L_q((0,1),\frac{dt}{t})}.
\end{equation*}
Indeed, $\mu(x)\leq\mathcal{C}\mu(x)$ and Lemma~\ref{lem:logarithmic-Hardy-estimate} give
\begin{equation*}
 \|x\|_{L_{\infty,q}}
 \leq \|x\|_{L_{\infty,q}}^\#
 \leq 2\|x\|_{L_{\infty,q}}.
\end{equation*}
We retain the rearrangement quasi-norm \eqref{eq:prelim-log-class} in the estimates below.

Our notation should be distinguished from the oscillation spaces used in \cite{Nursultanov}. In the commutative setting, their defining functional is
\begin{equation*}
 \|f\|_{L_1(0,1)} + \left\|t\mapsto f^{**}(t)-f^*(t) \right\|_{L_q((0,1),\frac{dt}{t})},
 \ \   f^{**}(t)=\frac{1}{t}\int\limits_0^t f^*(s)\,ds,
\end{equation*}
where $f^*$ denotes the decreasing rearrangement of $|f|$, written $\mu(f)$ in our notation; thus $f^*(t)=\mu(t,f)$. For $1<q\leq\infty$, this gives a different limiting scale. Here we work with the logarithmically weighted rearrangement in \eqref{eq:prelim-log-class} and establish convolution estimates for this scale on quantum tori; see Theorem~\ref{thm:working-Brezis-Wainger}.

Finally, the notation $L_{\infty,\infty}$ used here does not change the definition of $B_{\infty,\infty}^{s}$ in \eqref{eq:Besov-infinity-norm}, whose dyadic blocks are measured in the ordinary operator norm.

\subsection{Differential calculus}
\label{NCSobolevSpace}

The canonical derivations, Laplacian and Bessel potentials are
\begin{equation}\label{eq:canonical-derivations}
\partial_{j}^{\theta} e_n^{\theta}=in_{j}e_n^{\theta},\ \ 
\Delta_\theta e_n^{\theta}=-|n|^2e_n^{\theta},\ \ 
J_\theta^{s} e_n^{\theta}=(1+|n|^2)^{\frac{s}{2}}e_n^{\theta}, \ \ n\in \mathbb{Z}^{d}.
\end{equation}
On $L_2$, each derivation has a closed skew-adjoint realization, and $-\Delta_\theta$ is positive and self-adjoint. The differential operators and Bessel potentials extend to distributions by duality. More precisely, for $F\in\mathcal S'(\mathbb{T}_{\theta}^{d})$ and $a\in\mathcal S(\mathbb{T}_{\theta}^{d})$, we define
\begin{equation*}
(\partial_{j}^{\theta} F,a)=-(F,\partial_{j}^{\theta} a) \ \textnormal{ and }\ 
(J_\theta^{s}F,a)=(F,J_\theta^{s}a),
\end{equation*}
where $1\leq j\leq d$, $s\in\mathbb{R}$, and $(\cdot,\cdot)$ denotes the bilinear pairing between distributions and smooth test elements. Both $J_\theta^{s}$ and $J_\theta^{-s}$ have symbols of polynomial growth.
They preserve rapidly decreasing Fourier coefficients and act continuously
on $\mathcal S$ and $\mathcal S'$, with
\begin{equation*}
 J_\theta^{s}J_\theta^{-s}=J_\theta^{-s}J_\theta^{s}=I;
\end{equation*}
see \cite{XiongXuYin,McDSX}.

We write $A\lesssim B$ for $A\leq CB$ and $A\asymp B$ for two-sided bounds. Constants are independent of the deformation parameters unless stated otherwise.

\section{Convolution and the sharp O'Neil inequality}
\label{sec:convolution}

The aim of this section is to obtain a rearrangement inequality for operator convolution from endpoint norm estimates. Our approach extends the classical inequality of O'Neil~\cite{ONeil} to an abstract bilinear map between noncommutative $L_1$ spaces. We then apply this result to the convolution introduced in Part~I~\cite{ConvolutionI}. The resulting sharp inequality is the main tool for the convolution and Sobolev estimates below.

We recall the operation of \cite{ConvolutionI} in order to fix both its
domain and its normalization. Let $\theta_3=\theta_1+\theta_2$. The
trace-preserving normal $*$-homomorphism
\begin{equation*}
\mathcal{T}_{\theta_1,\theta_2}:
L_\infty(\mathbb{T}_{\theta_3}^{d})\to 
L_\infty(\mathbb{T}_{\theta_1}^{d})\,\overline\otimes\,
L_\infty(\mathbb{T}_{\theta_2}^{d}), \ \ 
\mathcal{T}_{\theta_1,\theta_2}(e_n^{\theta_3})
=e_n^{\theta_1}\otimes e_n^{\theta_2}
\end{equation*}
is the cocycle-addition map. For $L_{1}$ inputs, its preadjoint defines
$x\star_{\theta_1,\theta_2}y$ by
\begin{equation}\label{eq:convolution-recall}
\tau_{\theta_3}((x\star_{\theta_1,\theta_2}y)a)
=(\tau_{\theta_1}\otimes\tau_{\theta_2})
((x\otimes y)\mathcal{T}_{\theta_1,\theta_2}(a)),
 \ \  a\in L_\infty(\mathbb{T}_{\theta_3}^{d}).
\end{equation}
In particular,
$\widehat{x\star_{\theta_1,\theta_2}y}(n)
=\widehat x(n)\widehat y(n)$ and
$f\star_{0,\theta}x=f\ast x$. Part~I \cite[Propositions~3.5 and~3.10]{ConvolutionI} proves the
following bounds, including their consistency on the intersections
of the domains:
\begin{align}
\|x\star_{\theta_1,\theta_2}y\|_1
&\leq\|x\|_1\|y\|_1,\label{eq:L1}\\
\|x\star_{\theta_1,\theta_2}y\|_\infty
&\leq\|x\|_1\|y\|_\infty,\label{eq:L1-Linfty-sharp}\\
\|x\star_{\theta_1,\theta_2}y\|_\infty
&\leq\|x\|_\infty\|y\|_1.\label{eq:Linfty-L1-sharp}
\end{align}
It also proves \cite[Theorem~4.1]{ConvolutionI}, with optimal constant one, that
\begin{equation}\label{eq:Young}
\|x\star_{\theta_1,\theta_2}y\|_r
\leq\|x\|_p\|y\|_q,
\ \  1+\frac1r=\frac1p+\frac1q,
 \ \  1\leq p,q,r\leq\infty.
\end{equation}
The endpoint maps are restrictions of the same $L_{1}\times L_{1}$
bilinear operation, since $L_\infty\subset L_{1}$ for the normalized
trace. The submajorization inequalities of
\cite[Corollary~5.1]{ConvolutionI} imply that
convolution by an $L_{1}$ input is a contraction on every fully
symmetric operator space (and on $L_{p,q}$ with equivalent Banach
norm when $p>1$); we refer to this property below as
\begin{align}
\mu(x\star_{\theta_1,\theta_2}y)
&\prec\!\prec\|x\|_1\mu(y),\label{eq:HLP-first}\\
\mu(x\star_{\theta_1,\theta_2}y)
&\prec\!\prec\|y\|_1\mu(x).\label{eq:HLP-second}
\end{align}

\begin{definition}[O'Neil relation]\label{def:ONeil-relation}
Let $(\mathcal M_{j},\tau_{j})$, $j=1,2,3$, be finite von Neumann
algebras with normalized traces. A triple
$(x,y,z)\in L_{1}(\mathcal M_1)\times L_{1}(\mathcal M_2)
\times L_{1}(\mathcal M_3)$ satisfies the \emph{O'Neil relation} if
\begin{equation}\label{eq:ONeil-relation}
\int\limits_0^t\mu(s,z)\,ds\leq
\int\limits_0^t\mu(s,x)\,ds\cdot
\int\limits_0^t\mu(s,y)\,ds
+t\int\limits_t^1\mu(s,x)\mu(s,y)\,ds,
 \ \  0<t<1.
\end{equation}
\end{definition}

\begin{proposition}[Fourier multiplier compatibility]
\label{prop:multiplier-compatibility}
Let $m_1,m_2:\mathbb{Z}^{d}\to\mathbb{C}$ have at most polynomial
growth. The Fourier coefficient rule extends convolution to
$\mathcal S'(\mathbb{T}_{\theta_1}^{d})\times
\mathcal S'(\mathbb{T}_{\theta_2}^{d})$ and gives
\begin{equation}\label{eq:multiplier-convolution-commutation}
(T_{m_1}^{\theta_1}x)\star_{\theta_1,\theta_2}
(T_{m_2}^{\theta_2}y)
=T_{m_1m_2}^{\theta_3}(x\star_{\theta_1,\theta_2}y).
\end{equation}
In particular, one may transfer a multiplier from the
output to either input. These identities also agree with the
$L_{1}$ convolution whenever the corresponding terms are integrable.
\end{proposition}
\begin{proof}
The product of two sequences of polynomial growth still has
polynomial growth; hence the Fourier coefficient rule defines
a distribution. Both sides of
\eqref{eq:multiplier-convolution-commutation} have $n$th
coefficient $m_1(n)m_2(n)\widehat x(n)\widehat y(n)$.
Uniqueness of Fourier coefficients proves the assertion.
\end{proof}

\begin{cond}\label{main condition general} Let $\mathcal{M}_{1},\mathcal{M}_{2},\mathcal{M}_{3}$ be finite von Neumann algebras equipped with normal faithful traces $\tau_{j}$ satisfying $\tau_{j}(\mathbf{1})=1$. Consider a bilinear mapping $\mathfrak{B}:L_{1}(\mathcal{M}_{1})\times L_{1}(\mathcal{M}_{2})\to L_{1}(\mathcal{M}_{3})$ such that
\begin{equation}\label{Gassump-1}
\|\mathfrak{B}(x,y)\|_{L_{1}(\mathcal{M}_{3})}\leq\|x\|_{L_{1}(\mathcal{M}_{1})}\|y\|_{L_{1}(\mathcal{M}_{2})},
\end{equation}
\begin{equation}\label{Gassump-2}
\|\mathfrak{B}(x,y)\|_{L_{\infty}(\mathcal{M}_{3})}\leq\|x\|_{L_{1}(\mathcal{M}_{1})}\|y\|_{L_{\infty}(\mathcal{M}_{2})},
\end{equation}
\begin{equation}\label{Gassump-3}
\|\mathfrak{B}(x,y)\|_{L_{\infty}(\mathcal{M}_{3})}\leq\|x\|_{L_{\infty}(\mathcal{M}_{1})}\|y\|_{L_{1}(\mathcal{M}_{2})}.
\end{equation}
\end{cond}

We prove the following abstract version of O'Neil's lemma
\cite[Lemma~1.5]{ONeil}, including the measurable-operator extensions.
\begin{lem}\label{final B-lemma}
Suppose that $\mathfrak B$ satisfies Condition~\ref{main condition general}.
For $x\in L_{1}(\mathcal M_1)$ and $y\in\mathcal M_2$, or
$x\in\mathcal M_1$ and $y\in L_{1}(\mathcal M_2)$, one has
\begin{equation*}
\int\limits_0^t\mu(s,\mathfrak B(x,y))\,ds \leq\int\limits_0^t\mu(s,x)\,ds\cdot\int\limits_0^t\mu(s,y)\,ds+t\int\limits_t^1\mu(s,x)\mu(s,y)\,ds,\ \  0<t<1.
\end{equation*}
\end{lem}
\begin{proof}
We may assume the algebras are diffuse. Indeed, replace $\mathcal M_{j}$ by $\widetilde{\mathcal M}_{j} =\mathcal M_{j}\,\overline\otimes\,L_\infty(0,1)$ if necessary,
let $E_{j}=\operatorname{id}\otimes\int\limits_{0}^{1}$ be the partial trace, and set
\begin{equation*}
 \widetilde{\mathfrak B}(a,b)
=\mathfrak B(E_1a,E_2b)\otimes\mathbf{1}.   
\end{equation*}
The partial traces are contractions on $L_{1}$ and $L_\infty$, so all three hypotheses are preserved. Applying the diffuse case to $x\otimes\mathbf{1},y\otimes\mathbf{1}$ gives the original assertion, because tensoring with $\mathbf{1}$ preserves singular values.

Fix $x\in L_{1}(\mathcal M_1)$ and $y\in\mathcal M_2$; the other case is identical. Choose diffuse abelian subalgebras $\mathcal A_{j}\subset\mathcal M_{j}$, $j=1,2$, and trace-preserving $*$-isomorphisms $i_{j}:L_\infty(0,1)\to\mathcal A_{j}$ with
\begin{equation*}
i_1(\mu(x))=|x|,\ \  i_2(\mu(y))=|y|.
\end{equation*}
Here and below $i_{j}$ also denotes the canonical extension $S(0,1)\to S(\mathcal A_{j},\tau_{j})$. This extension preserves singular values and restricts to an isometry on $L_{1}$; the first equality is therefore meaningful for unbounded $x$. Such abelian representations are obtained by splitting spectral projections in the diffuse ambient algebras; see \cite[Lemma~4.4.12 and the subsequent extension argument]{DoddsDePagterSukochev}.

Let $x=u_1|x|$, $y=u_2|y|$ and $\mathfrak B(x,y)=u|\mathfrak B(x,y)|$ be the polar decompositions. Choose a diffuse abelian subalgebra $\mathcal A_3\subset\mathcal M_3$ containing $|\mathfrak B(x,y)|$, and a trace-preserving $*$-isomorphism $i_3:\mathcal A_3\to L_\infty(0,1)$. These abelian algebras may be chosen with separable predual and hence with the indicated standard probability-space representation. Write $i_3$ also for its measurable-operator and $L_{1}$ extensions.

Define the map $\mathfrak{B}_0:L_{\infty}(0,1)\times L_{\infty}(0,1)\to L_{\infty}(0,1)$ by setting
\begin{equation*}
\mathfrak{B}_0(f,g)=i_{3}\Big(\mathbb{E}_{\mathcal{A}_{3}}(u^{\ast}\cdot \mathfrak{B}(u_1\cdot i_1(f), u_2\cdot i_2(g)))\Big), \ \ f,g\in L_{\infty}(0,1),
\end{equation*}
where $\mathbb{E}_{\mathcal{A}_{3}}:\mathcal{M}_{3}\to \mathcal{A}_{3}$ is the trace-preserving conditional expectation \cite[Chapter 7, p. 456]{DoddsDePagterSukochev}. Since $\mathfrak{B}$ is bilinear and $\mathbb{E}_{\mathcal{A}_{3}}$ is linear, it follows that $\mathfrak{B}_{0}$ is bilinear. We have
\begin{equation*}
\begin{split}
\|\mathfrak{B}_0(f,g)\|_{L_{\infty}(0,1)}&=\|\mathbb{E}_{\mathcal{A}_{3}}(u^{\ast}\cdot \mathfrak{B}(u_1\cdot i_1(f), u_2\cdot i_2(g))) \|_{L_{\infty}(\mathcal{A}_3,\tau_3)}\\
&\leq \|u^{\ast}\cdot \mathfrak{B}(u_1\cdot i_1(f), u_2\cdot i_2(g))\|_{L_{\infty}(\mathcal{M}_{3},\tau_3)} \\
   &\leq \|\mathfrak{B}(u_1\cdot i_1(f), u_2\cdot i_2(g))\|_{L_{\infty}(\mathcal{M}_{3},\tau_3)}. 
\end{split}
\end{equation*}

Hence, by \eqref{Gassump-2}, we obtain 
\begin{equation*}
\begin{split}
\| \mathfrak{B}_0(f,g)\|_{L_{\infty}(0,1)}
   &\leq \|\mathfrak{B}(u_1\cdot i_1(f), u_2\cdot i_2(g))\|_{L_{\infty}(\mathcal{M}_{3},\tau_3)}\\
   &\stackrel{\eqref{Gassump-2}}{\leq} \|u_{1}\cdot i_{1}(f)\|_{L_{1}(\mathcal{M}_{1},\tau_1)} \|u_{2}\cdot i_{2}(g)\|_{L_{\infty}(\mathcal{M}_{2},\tau_2)}\\
   &\leq\|i_{1}(f)\|_{L_{1}(\mathcal{M}_{1},\tau_1)} \| i_{2}(g)\|_{L_{\infty}(\mathcal{M}_{2},\tau_2)}\\
   &=\|f\|_{L_{1}(0,1)} \| g\|_{L_{\infty}(0,1)}.
\end{split}
\end{equation*}
Similarly, by \eqref{Gassump-3}, we have
\begin{equation*}
\begin{split}
   \| \mathfrak{B}_0(f,g)\|_{L_{\infty}(0,1)}
   \leq \|f\|_{L_{\infty}(0,1)} \| g\|_{L_{1}(0,1)}.
\end{split}
\end{equation*}
And, finally, since $\mathbb{E}_{\mathcal{A}_{3}}$ can be extended to a contraction from $L_{1}(\mathcal{M}_{3},\tau_3)$ to $L_{1}(\mathcal{A}_{3},\tau_3)$ (\cite[Chapter 7, Proposition 7.1.1, p. 457]{DoddsDePagterSukochev}), it follows from \eqref{Gassump-1} that
\begin{equation*}
\begin{split}
\|\mathfrak{B}_0(f,g)\|_{L_{1}(0,1)}&=\| \mathbb{E}_{\mathcal{A}_{3}}(u^{\ast}\cdot \mathfrak{B}(u_1\cdot i_1(f), u_2\cdot i_2(g))) \|_{L_{1}(\mathcal{A}_3,\tau_3)}\\
&\leq\|u^{\ast}\cdot \mathfrak{B}(u_1\cdot i_1(f), u_2\cdot i_2(g))\|_{L_{1}(\mathcal{M}_{3},\tau_3)} \\
&\leq\|\mathfrak{B}(u_1\cdot i_1(f), u_2\cdot i_2(g))\|_{L_{1}(\mathcal{M}_{3},\tau_3)}\\
& \stackrel{\eqref{Gassump-1}}{\leq} \|u_{1}\cdot i_{1}(f)\|_{L_{1}(\mathcal{M}_{1},\tau_1)}\| u_{2}\cdot i_{2}(g)\|_{L_{1}(\mathcal{M}_{2},\tau_2)}\\
&\leq\|f\|_{L_{1}(0,1)}\|g\|_{L_{1}(0,1)}.
\end{split}
\end{equation*}
Therefore, $\mathfrak{B}_0$ fits into the classical O'Neil setting \cite[Lemma 1.5]{ONeil}, and for any $f,g\in L_{\infty}(0,1)$, we have
\begin{equation}\label{oneil inequality}
  \int\limits_0^t\mu(s,\mathfrak{B}_0(f,g))ds\leq \int\limits_0^t\mu(s,f)ds\cdot\int\limits_0^t\mu(s,g)ds+t\int\limits_t^1\mu(s,f)\mu(s,g)ds, \ \  t\in(0,1).  
\end{equation}
The three bounds extend $\mathfrak B_0$ consistently to $L_{1}(0,1)\times L_{1}(0,1)$ and to both mixed endpoint domains. Inequality~\eqref{oneil inequality} also holds when one input is integrable and the other bounded: approximate the integrable input by its bounded polar truncations, use $L_{1}$ convergence on the left and monotone convergence of the rearrangement expressions on the right. The passage on the left follows from \eqref{eq:K-prelim} and $|K(t,a)-K(t,b)|\leq\|a-b\|_1$. The definition of $\mathfrak B_0$ agrees with this extension because $i_1,i_2,i_3$ and the conditional expectation are continuous on $L_{1}$.
Recall that
\begin{equation*}
  u_1\cdot i_1(\mu(x))=x,\ \ u_2\cdot i_2(\mu(y))=y,\ \ u^{\ast}\cdot\mathfrak{B}(x,y)=|\mathfrak{B}(x,y)|.  
\end{equation*}
Thus,
\begin{equation*}
  \begin{split}
\mathfrak{B}_0(\mu(x),\mu(y))&=i_{3}\Big(\mathbb{E}_{\mathcal{A}_{3}}(u^{\ast}\cdot \mathfrak{B}(u_1\cdot i_1(\mu(x)), u_2\cdot i_2(\mu(y))))\Big)\\
&=i_{3}\Big(\mathbb{E}_{\mathcal{A}_{3}}(u^{\ast}\cdot \mathfrak{B}(x,y))\Big)=i_{3}\Big(\mathbb{E}_{\mathcal{A}_{3}}(|\mathfrak{B}(x,y)|)\Big).
  \end{split}  
\end{equation*}
Since $|\mathfrak{B}(x,y)|$ is affiliated to $\mathcal{A}_3,$ it follows that
\begin{equation*}
  \mathfrak{B}_0(\mu(x),\mu(y))=i_{3}\Big(|\mathfrak{B}(x,y)|\Big).  
\end{equation*}
In particular,
\begin{equation*}
\mu\big(\mathfrak{B}_0(\mu(x),\mu(y))\big)=\mu\big(\mathfrak{B}(x,y)\big).  
\end{equation*}
Using \eqref{oneil inequality} with $f=\mu(x)$ and $g=\mu(y),$ we complete the proof.
\end{proof}

\begin{corollary}[Extension to two integrable inputs]
\label{lem:bilinear-O-Neil}
Under Condition~\ref{main condition general}, the triple $(x,y,\mathfrak B(x,y))$ satisfies the O'Neil relation \eqref{eq:ONeil-relation} for every $x\in L_{1}(\mathcal M_1)$ and $y\in L_{1}(\mathcal M_2)$.
\end{corollary}
\begin{proof}
Write $y=v|y|$ and put $y_N=v\min(|y|,N)$. Then $y_N\in\mathcal M_2$, $\|y_N-y\|_1\to0$, and $\mu(t,y_N)=\min\left\{\mu(t,y),N\right\}$. By \eqref{Gassump-1}, $\mathfrak B(x,y_N)\to\mathfrak B(x,y)$ in $L_{1}$. Apply Lemma~\ref{final B-lemma} to $(x,y_N)$ and pass to the limit using \eqref{eq:K-prelim} on the left and monotone convergence on the right.
\end{proof}

\begin{thm}[Sharp two-operator O'Neil inequality]
\label{thm:operator-oneil}
For $x\in L_{1}(\mathbb{T}_{\theta_1}^{d})$ and $y\in L_{1}(\mathbb{T}_{\theta_2}^{d})$, the triple $(x,y,x\star_{\theta_1,\theta_2}y)$ satisfies \eqref{eq:ONeil-relation}. Equivalently,
\begin{equation}\label{eq:ONeil-Cesaro}
\mathcal{C}\mu(x\star_{\theta_1,\theta_2}y)(t) \leq t\,\mathcal{C}\mu(x)(t)\cdot\mathcal{C}\mu(y)(t)  + \int\limits_t^1\mu(s,x)\cdot\mu(s,y)\,ds.
\end{equation}
The multiplicative constant one cannot be improved.
\end{thm}
\begin{proof}
Apply Corollary~\ref{lem:bilinear-O-Neil} to $\mathfrak B(x,y)=x\star_{\theta_1,\theta_2}y$. The three endpoint hypotheses are precisely \eqref{eq:L1}--\eqref{eq:Linfty-L1-sharp}. Here $\mathcal M_{j}=L_\infty(\mathbb{T}_{\theta_{j}}^{d})$, with $\theta_3=\theta_1+\theta_2$. Dividing the integral inequality by $t$ gives \eqref{eq:ONeil-Cesaro}. For $x=y=\mathbf{1}$ the left side of \eqref{eq:ONeil-relation} equals $t$ and the right side equals $t^2+t(1-t)=t$.
\end{proof}

\begin{corollary}\label{oneil trivial corollary} We have
$$\mu(x\star_{\theta_1,\theta_2}y)
\leq 3\sigma_2\mathfrak{T}(\mathcal{C}\mu(x)\cdot\mathcal{C}\mu(y)).$$
\end{corollary}

\begin{proof}
Put $h=(\mathcal{C}\mu(x))(\mathcal{C}\mu(y))$. Since $h$ is nonnegative and nonincreasing,
\begin{equation*}
 th(t)\leq2\int\limits_{\frac{t}{2}}^{t}h(s)\,ds  \leq2(\mathfrak Th)\left( \frac{t}{2} \right).
\end{equation*}
Moreover, $\mu(x)\mu(y)\leq h$. Thus Theorem~\ref{thm:operator-oneil} gives
\begin{align*}
 \mu(t,x\star_{\theta_1,\theta_2}y)
 &\leq th(t)+\int\limits_t^1\mu(s,x)\mu(s,y)\,ds \leq2(\mathfrak Th)\left( \frac{t}{2} \right)+(\mathfrak Th)(t)\\
 &\leq3(\mathfrak Th)\left( \frac{t}{2} \right)
 =3(\sigma_2\mathfrak Th)(t),
\end{align*}
as required.
\end{proof}
\begin{corollary}[Sharp semi-convolution inequality]
\label{cor:sharp-semi-ONeil}
For $f\in L_{1}(\mathbb{T}^{d})$ and $x\in L_{1}(\mathbb{T}_{\theta}^{d})$, $(f,x,f\ast x)$ satisfies \eqref{eq:ONeil-relation}.
\end{corollary}
\begin{proof}
Use $f\ast x=f\star_{0,\theta}x$ in Theorem~\ref{thm:operator-oneil}.
\end{proof}
\begin{remark}
The sharp semi-convolution inequality gives the Hardy-operator estimate with constant $2$ in Theorem~\ref{thm:fractional-hardiness} below, using the weak-$L_2$ rearrangement quasi-norm.
\end{remark}

\section{Symmetric Sobolev and Besov spaces}
\label{sec:potential-structure}
We now study Sobolev and Besov spaces based on symmetric operator spaces. This extends the $L_p$-based theory of Xiong, Xu and Yin~\cite{XiongXuYin}. We establish completeness, approximation by Fourier polynomials, and basic properties of both scales. For Besov spaces, we also prove multiplication estimates using an ordered version of Bony's decomposition \cite{Bony,BahouriCheminDanchin}. These estimates will be used for nonlinear heat equations.

We present the two scales in parallel: definitions, completeness and approximation, followed by the estimates used in the applications. The Sobolev scale is defined for symmetric Banach bases; the Besov estimates below use fully symmetry. These hypotheses are stated separately.

\subsection{Sobolev spaces: definition and structure}
For a symmetric Banach function space $E$, $s\in\mathbb{R}$, and
$k\in\mathbb N_0$, define
\begin{align}
H_{\mathcal{E}}^{s}(\mathbb{T}_{\theta}^{d})
 &=\left\{x\in\mathcal S':J_\theta^{s}x\in\mathcal{E}\right\},
&\|x\|_{H_{\mathcal{E}}^{s}}&=\|J_\theta^{s}x\|_{\mathcal{E}},
\label{eq:prelim-potential-space}\\
W_{\mathcal{E}}^k(\mathbb{T}_{\theta}^{d})
 &=\left\{x\in\mathcal S':(\partial^{\theta})^\beta x\in\mathcal{E},
 |\beta|\leq k\right\},
&\|x\|_{W_{\mathcal{E}}^k}&=\sum_{|\beta|\leq k}
\|(\partial^{\theta})^\beta x\|_{\mathcal{E}}.
\label{eq:prelim-derivative-space}
\end{align}
We write $H_{p,q}^{s},W_{p,q}^k$ for the Lorentz choices of $E$, and
$H_p^{s},W_p^k$ for $E=L_p$. The potential and derivative definitions
are equivalent for integer order when $1<p<\infty$; no endpoint
identification at $p=1$ or $p=\infty$ is assumed.

\begin{prop}[Completeness and lifting]
\label{prop:prelim-completeness}\label{lem:J^{s}-iso}
For a symmetric Banach function space $E$ and $s,t\in\mathbb{R}$,
$H_{\mathcal{E}}^{s}$ is Banach and
\begin{equation*}
 J_\theta^t:H_{\mathcal{E}}^{s} \to  H_{\mathcal{E}}^{s-t}
\end{equation*}
is a surjective isometry. In particular $J_\theta^{s}$ identifies
$H_{\mathcal{E}}^{s}$ isometrically with $\mathcal{E}$.
\end{prop}
\begin{proof}
The multiplier identity gives
\begin{equation*}
 \|J_\theta^t x\|_{H_{\mathcal{E}}^{s-t}}
 =\|J_\theta^{s-t}J_\theta^tx\|_{\mathcal{E}}
 =\|x\|_{H_{\mathcal{E}}^{s}},\ \  (J_\theta^t)^{-1}=J_\theta^{-t}.
\end{equation*}
If $(x_n)$ is Cauchy in $H_{\mathcal{E}}^{s}$, let
$J_\theta^{s}x_n\to y$ in $\mathcal{E}$ and set $x=J_\theta^{-s}y$.
Then $x\in H_{\mathcal{E}}^{s}$ and
$\|x_n-x\|_{H_{\mathcal{E}}^{s}}=\|J_\theta^{s}x_n-y\|_{\mathcal{E}}\to0$.
\end{proof}

\begin{thm}[Polynomial approximation]\label{thm:HEs-properties}
If $E$ is a separable symmetric Banach function space, Fourier
polynomials are dense in $\mathcal{E}$ and in $H_{\mathcal{E}}^{s}$,
$s\in\mathbb{R}$. Consequently $H_{\mathcal{E}}^{s}$ is separable.
\end{thm}
\begin{proof}
For bounded $x$, its product Fej\'er means $F_Nx$ are Fourier
polynomials with
\begin{equation*}
 \|F_Nx\|_\infty\leq\|x\|_\infty,\ \ 
 \|F_Nx-x\|_1\to0;
\end{equation*}
see \cite[Proposition~3.1]{ChenXuYin}.
For $x\ne0$, put $y_N=(F_Nx-x)/(2\|x\|_\infty)$.
By \cite[Lemma~6]{SukochevTulenovZanin}, there is a function
$\psi_{\mathcal{E}}$ with $\psi_{\mathcal{E}}(0+)=0$ such that
\begin{equation*}
 \|y_N\|_\infty\leq1,\ \  \|y_N\|_1\to0
 \ \ \Rightarrow \ \ 
 \|y_N\|_{\mathcal{E}}\leq\psi_{\mathcal{E}}(\|y_N\|_1)\to0.
\end{equation*}
Since bounded operators are dense in $\mathcal{E}$
\cite[Lemma~7]{SukochevTulenovZanin}, this proves polynomial density.
For $x\in H_{\mathcal{E}}^{s}$, choose polynomials $p_N\to J_\theta^{s}x$
in $\mathcal{E}$. Then $J_\theta^{-s}p_N$ are polynomials and
\begin{equation*}
 \|J_\theta^{-s}p_N-x\|_{H_{\mathcal{E}}^{s}}
 =\|p_N-J_\theta^{s}x\|_{\mathcal{E}}\to0.
\end{equation*}
Polynomials with coefficients in $\mathbb Q+i\mathbb Q$ form a
countable dense set, proving separability.
\end{proof}

\begin{prop}[Derivative spaces]
\label{prop:additional-lifting}
The space
$W_{\mathcal{E}}^k$ is Banach. If $E$ is separable, Fourier polynomials
are dense in $W_{\mathcal{E}}^k$.
\end{prop}
\begin{proof}
For completeness of $W_{\mathcal{E}}^k$, let $(x_{j})$ be Cauchy.
For each $|\beta|\leq k$, there is $y_\beta\in\mathcal{E}$ such that
$(\partial^{\theta})^\beta x_{j}\to y_\beta$ in $\mathcal{E}$.
Taking $\beta=0$ gives $x_{j}\to y_0$ in $\mathcal{E}$ and in $\mathcal S'$.
For every smooth $a$, distributional differentiation therefore gives
\begin{align*}
(y_\beta,a) =\lim_{j\to\infty}((\partial^{\theta})^\beta x_{j},a) = (-1)^{|\beta|}\lim_{j\to\infty}(x_{j},(\partial^{\theta})^\beta a) = ((\partial^{\theta})^\beta y_0,a).
\end{align*}
Consequently $y_\beta=(\partial^{\theta})^\beta y_0$, $y_0\in W_{\mathcal{E}}^k$, and
\begin{equation*}
\|x_{j}-y_0\|_{W_{\mathcal{E}}^k} = \sum_{|\beta|\leq k} \|(\partial^{\theta})^\beta x_{j}-y_\beta\|_{\mathcal{E}}\to 0.
\end{equation*}

Finally, let $E$ be separable. We first justify averaging the torus
action in the $\mathcal{E}$ norm. For $z\in L_\infty$, the action
$t\mapsto\alpha_t^{\theta}(z)$ is continuous in $L_{1}$, and
\begin{equation*}
\|\alpha_t^{\theta}(z)-z\|_\infty\leq2\|z\|_\infty.
\end{equation*}
Applying the function $\psi_{\mathcal{E}}$ used in the proof of
Theorem~\ref{thm:HEs-properties}, after dividing by
$2\|z\|_\infty$ when $z\ne0$, proves continuity in $\mathcal{E}$.
Each $\alpha_t^{\theta}$ preserves the singular value function and is
therefore an isometry on $\mathcal{E}$. Given a general $z\in\mathcal{E}$,
choose $z_m\in L_\infty$ converging to $z$ in $\mathcal{E}$. Then
\begin{equation*}
\|\alpha_t^{\theta}(z)-z\|_{\mathcal{E}}
\leq2\|z-z_m\|_{\mathcal{E}}
   +\|\alpha_t^{\theta}(z_m)-z_m\|_{\mathcal{E}}.
\end{equation*}
First choosing $m$ and then letting $t\to0$ proves continuity of
the action on all of $\mathcal{E}$.

Let $\mathsf F_N$ be the scalar product Fej\'er kernel (see, \cite[Proposition 2.7.(ii)]{XiongXuYin}). It is nonnegative and has normalized integral one. The Bochner integral
\begin{equation*}
F_Nz=\int\limits_{\mathbb{T}^{d}}\mathsf F_N(t)\alpha_t^{\theta}(z)\,dm(t)
\end{equation*}
exists in $\mathcal{E}$ and satisfies
\begin{equation*}
\|F_Nz\|_{\mathcal{E}}
\leq\int\limits_{\mathbb{T}^{d}}\mathsf F_N(t)
       \|\alpha_t^{\theta}(z)\|_{\mathcal{E}}\,dm(t)
=\|z\|_{\mathcal{E}}.
\end{equation*}
The convergence $F_Nz_m\to z_m$ for bounded $z_m$, already proved
in Theorem~\ref{thm:HEs-properties}, now gives
\begin{equation*}
\|F_Nz-z\|_{\mathcal{E}}
\leq2\|z-z_m\|_{\mathcal{E}}+\|F_Nz_m-z_m\|_{\mathcal{E}}
\to 0
\end{equation*}
by taking $N\to\infty$ and then $m\to\infty$.
Since $F_N$ commutes
with the derivations,
\begin{equation*}
\|F_Nx-x\|_{W_{\mathcal{E}}^k}
=\sum_{|\beta|\leq k}
\|F_N((\partial^{\theta})^\beta x)-(\partial^{\theta})^\beta x\|_{\mathcal{E}}
\to 0.
\end{equation*}
Each $F_Nx$ is a Fourier polynomial. This proves the density assertion.
\end{proof}

\label{sec:symmetric-besov-spaces}

\subsection{Besov spaces: definition and dyadic resolution}
\label{subsec:besov-preliminaries}

Let $E=E(0,1)$ be a fully symmetric Banach function space, and let $\mathcal{E}=\mathcal{E}(\mathbb{T}_\theta^{d})$ be its associated operator space, with
\begin{equation*}
 \|x\|_{\mathcal{E}}=\|\mu(x)\|_E.
\end{equation*}
The canonical trace and Haar measure are normalized. We write $E^\times$ for the K\"othe associate space of $E$ and put
\begin{equation*}
 \phi_E(t)=\|\mathbf{1}_{(0,t)}\|_E,
 \ \ 0<t\leq1.
\end{equation*}
We use the standard facts
\begin{equation}\label{eq:symmetric-besov-base-properties}
 \begin{aligned}
 L_\infty(\mathbb{T}_\theta^{d})\hookrightarrow\mathcal{E}\hookrightarrow L_1(\mathbb{T}_\theta^{d}), \ & \ \|azb\|_{\mathcal{E}}\leq        \|a\|_\infty\|z\|_{\mathcal{E}}\|b\|_\infty,\\
 \phi_E(t)\phi_{E^\times}(t)&=t, \ \ 0<t\leq1,
 \end{aligned}
\end{equation}
where $a,b\in L_\infty(\mathbb{T}_\theta^{d})$ and $z\in\mathcal{E}$. See \cite{BennettSharpley,DoddsDePagterSukochev} for the function and operator space background. No order-continuity assumption on $E$ is imposed.

Choose a real even radial function $\chi\in C_c^\infty(\mathbb{R}^{d})$ such that $0\leq\chi\leq1$, $\chi=1$ on $\left\{|\xi|\leq\frac{3}{4}\right\}$, and $\operatorname{supp}\chi\subset\left\{|\xi|<\frac{4}{3}\right\}$. Set
\begin{equation*}
 \varphi_0=\chi, \ \ \varphi(\xi)=\chi(\xi)-\chi(2\xi), \ \ \varphi_{j}(\xi)=\varphi(2^{-j}\xi),
 \ \ j\geq1.
\end{equation*}
The telescoping identity gives
\begin{equation}\label{eq:prelim-dyadic-partition}
 \sum_{j=0}^J\varphi_{j}(\xi)=\chi(2^{-J}\xi), \ \ \sum_{j\geq0}\varphi_{j}(\xi)=1.
\end{equation}
For $j\geq1$,
\begin{equation}\label{eq:symmetric-besov-annuli}
 \operatorname{supp}\varphi_{j}
 \subset\left\{\xi\in\mathbb{R}^{d}:
       \frac38\,2^j\leq|\xi|\leq\frac43\,2^j\right\}.
\end{equation}

For $x\in\mathcal S'(\mathbb{T}_\theta^{d})$, define
\begin{equation}\label{eq:prelim-dyadic-blocks}
 \Delta_{j}^{\theta} x=T_{\varphi_{j}}^{\theta} x
 =\sum_{n\in\mathbb Z^{d}}\varphi_{j}(n)\widehat x(n)e_n^{\theta},
 \ \ j\geq0.
\end{equation}
Each block is a Fourier polynomial. We also put
\begin{equation}\label{eq:symmetric-besov-low-cutoff}
 P_{j}^{\theta}=\sum_{j=0}^J\Delta_{j}^{\theta}
           =T_{\chi(2^{-J}\cdot)}^{\theta},
 \ \ J\geq0.
\end{equation}
The polynomial growth of the Fourier coefficients of a distribution and \eqref{eq:prelim-dyadic-partition} imply $P_{j}^{\theta} x\to x$ in $\mathcal S'(\mathbb{T}_\theta^{d})$. We omit $\theta$ from the dyadic operators when the torus is clear.

For $s\in\mathbb{R}$ and $1\leq r\leq\infty$, the symmetric Besov space $B_{\mathcal{E},r}^{s}(\mathbb{T}_\theta^{d})$ consists of all distributions for which
\begin{equation}\label{eq:besov-app-definition}
 \|x\|_{B_{\mathcal{E},r}^{s}}
 =\left\|\bigl(2^{js}\|\Delta_{j}^{\theta} x\|_{\mathcal{E}}\bigr)_{j\geq0}
   \right\|_{\ell_r}<\infty.
\end{equation}
When $E=L_p$, $1\leq p\leq\infty$, we write $B_{p,r}^{s}$. In particular,
\begin{equation}\label{eq:Besov-infinity-norm}
 \|x\|_{B_{\infty,\infty}^{s}(\mathbb{T}_\theta^{d})}
 =\sup_{j\geq0}2^{js}\|\Delta_{j}^{\theta} x\|_\infty.
\end{equation}
The norm inside this supremum is the ordinary operator norm.

Let $\widetilde\varphi_{j}$ be the scalar trigonometric polynomial whose Fourier coefficients are $\varphi_{j}(n)$. Then $\Delta_{j}^{\theta} x=\widetilde\varphi_{j}\ast x$; in particular, $\widehat{\widetilde\varphi_0}(n)=\chi(n)$. The Fourier coefficient rule for semi-convolution gives
\begin{equation}\label{eq:LP-convolution-compatibility}
 \Delta_{j}^{\theta}(f\ast x)  =(\Delta_{j}^{\mathrm{cl}}f)\ast x  =f\ast (\Delta_{j}^{\theta} x),
\end{equation}
whenever the semi-convolutions are defined.

\subsection{Besov spaces: kernel estimates and structure}
\label{subsec:symmetric-besov-basic}

We first record the kernel estimate that replaces the usual $L_p$ Bernstein estimate in the present setting.

\begin{lem}[Localized kernels]
\label{lem:symmetric-besov-localized-kernels}
Let $m\in C_c^\infty(\mathbb{R}^{d})$, and let $k_{m,j}$ be the scalar periodic kernel with Fourier coefficients $m(2^{-j}n)$, $j\geq0$. Then
\begin{equation}\label{eq:symmetric-besov-kernel-norms}
 \|k_{m,j}\|_{L_1(\mathbb{T}^{d})}\leq C_{d,m},
 \ \ 
 \|k_{m,j}\|_{E^\times(\mathbb{T}^{d})}
 \leq\frac{C_{d,m}}{\phi_E(2^{-jd})}.
\end{equation}
The constants are independent of $j$ and $\theta$.
\end{lem}
\begin{proof}
Poisson summation and the rapid decay of the inverse Fourier transform of $m$ give, for every $N>d$,
\begin{equation*}
 |k_{m,j}(u)|
 \leq C_{d,m,N}2^{jd}
       \bigl(1+2^j\operatorname{dist}(u,0)\bigr)^{-N},
 \ \ u\in\mathbb{T}^{d}.
\end{equation*}
Here the distance is taken on the torus. This is the standard periodization argument; see \cite{G2014} and \cite[Chapter~2]{BahouriCheminDanchin}. Integration gives the first estimate.

Put $v=2^{-jd}$. Taking $N\geq2d$ in the preceding bound and using the volume estimate for balls on $\mathbb{T}^{d}$, we obtain
\begin{equation*}
 \mu(t,k_{m,j})
 \leq Cv^{-1}\left(1+\frac tv\right)^{-2},
 \ \ 0<t<1.
\end{equation*}
Consequently,
\begin{equation*}
 \mu(t,k_{m,j})
 \leq Cv^{-1}\sum_{k\geq0}2^{-2k}
       \mathbf{1}_{(0,\min\left\{2^kv,1\right\})}(t).
\end{equation*}
The fundamental function of $E^\times$ satisfies
\begin{equation*}
 \phi_{E^\times}(\min\left\{2^kv,1\right\})
 \leq2^k\phi_{E^\times}(v).
\end{equation*}
The triangle inequality and
\eqref{eq:symmetric-besov-base-properties} therefore yield
\begin{align*}
 \|k_{m,j}\|_{E^\times}
 &\leq Cv^{-1}\sum_{k\geq0}2^{-2k}
              \phi_{E^\times}(\min\left\{2^kv,1\right\})\\
 &\leq Cv^{-1}\phi_{E^\times}(v)\sum_{k\geq0}2^{-k}
 \leq\frac{C}{\phi_E(v)}.
\end{align*}
\end{proof}

\begin{lem}[Uniform bounds for dyadic cutoffs]
\label{lem:besov-dyadic-cutoffs}
The operators $\Delta_{j}^{\theta}$ and $P_{j}^{\theta}$ are uniformly bounded on $\mathcal{E}$ and on $L_\infty(\mathbb{T}_\theta^{d})$. More precisely,
\begin{equation*}
 \|\Delta_{j}^{\theta} x\|_{\mathcal{E}}
 +\|P_{j}^{\theta} x\|_{\mathcal{E}}
 \leq C\|x\|_{\mathcal{E}},
 \ \ j,J\geq0.
\end{equation*}
The same estimate holds with $\mathcal{E}$ replaced by $L_\infty$. The constant is independent of $j$, $J$, and $\theta$.
\end{lem}
\begin{proof} 
The action formula for semi-convolution gives
\begin{equation*}
 \mu(k\ast x)\prec\prec\|k\|_1\mu(x),
 \ \ k\in L_1(\mathbb{T}^{d}),\ \ x\in\mathcal{E}.
\end{equation*}
Indeed, the action integral is taken in $L_1$, and the triangle inequality for integrated singular values applies because the canonical action preserves $\mu(x)$. Full symmetry therefore gives
\begin{equation}\label{eq:symmetric-besov-kernel-action}
 \|k\ast x\|_{\mathcal{E}}
 \leq\|k\|_1\|x\|_{\mathcal{E}}.
\end{equation}
Applying Lemma~\ref{lem:symmetric-besov-localized-kernels} with $m=\varphi$ and $m=\chi$ proves the estimates for $\Delta_{j}^{\theta}$ and $P_{j}^{\theta}$. Taking $E=L_\infty$ gives the corresponding operator-norm bounds. All constants are independent of $j$, $J$, and $\theta$.
\end{proof}

\begin{lem}[Estimate for localized blocks]
\label{lem:besov-localized-blocks}
For every $x\in\mathcal S'(\mathbb{T}_\theta^{d})$ and $j\geq0$,
\begin{equation}\label{eq:besov-app-block}
 \|\Delta_{j}^{\theta} x\|_\infty
 \leq\frac{C}{\phi_E(2^{-jd})}
          \|\Delta_{j}^{\theta} x\|_{\mathcal{E}},
\end{equation}
where $C$ is independent of $j$ and $\theta$.
\end{lem}
\begin{proof} 
Choose $\eta\in C_c^\infty(\mathbb{R}^{d})$ equal to $1$ on $\left\{|\xi|\leq\frac{4}{3}\right\}$, and put $K_{j}=k_{\eta,j}$. For $z=\Delta_{j}^{\theta} x$, the Fourier support condition gives $z=K_{j}\ast z$. Since $K_{j}$ and $z$ are Fourier polynomials, Theorem~\ref{thm:operator-oneil}, with $t\downarrow0$, yields
\begin{equation*}
 \|z\|_\infty
 \leq\int\limits_{0}^{1}\mu(t,K_{j})\mu(t,z)\,dt
 \leq\|K_{j}\|_{E^\times}\|z\|_{\mathcal{E}}.
\end{equation*}
The last inequality is the associate-space H\"older inequality. Lemma~\ref{lem:symmetric-besov-localized-kernels} now gives
\begin{equation*}
 \|\Delta_{j}^{\theta} x\|_\infty
 \leq\frac{C}{\phi_E(2^{-jd})}
          \|\Delta_{j}^{\theta} x\|_{\mathcal{E}},
\end{equation*}
with $C$ independent of $j$ and $\theta$.
\end{proof}

\begin{prop}[Completeness and approximation]
\label{lem:besov-app-basic}
Let $s\in\mathbb{R}$ and $1\leq r\leq\infty$. Then $B_{\mathcal{E},r}^{s}(\mathbb{T}_\theta^{d})$ is a Banach space, and its norm is independent, up to equivalence, of the chosen smooth inhomogeneous dyadic resolution.
If $r<\infty$, then
\begin{equation*}
 \|P_{j}^{\theta} x-x\|_{B_{\mathcal{E},r}^{s}}\to 0 \text{ as }J\to\infty.
\end{equation*}
In particular, Fourier polynomials are dense in $B_{\mathcal{E},r}^{s}$. The norm-equivalence constants are independent of $\theta$.
\end{prop}
\begin{proof} 
We suppress $\theta$ in the dyadic operators. The norm axioms follow from those of $\mathcal{E}$ and $\ell_r$, together with the reconstruction identity.

Let $\left\{x_m\right\}_{m=1}^{\infty}$ be Cauchy in $B_{\mathcal{E},r}^{s}$. For each $j$, the sequence $\left\{\Delta_{j}x_m\right\}_{m=1}^{\infty}$ converges in $\mathcal{E}$, and its limit has the same finite Fourier support. Moreover, the dyadic partition and the embedding $\mathcal{E}\hookrightarrow L_1$ give
\begin{align*}
 |\widehat x_m(n)| \leq\sum_{j:\,\varphi_{j}(n)\neq0} \|\Delta_{j}x_m\|_1  \leq C_{E,s}(1+|n|)^{-s} \|x_m\|_{B_{\mathcal{E},r}^{s}}.
\end{align*}
Here only a fixed number of indices occur, and $2^j$ is comparable to $1+|n|$ for those indices. Applying the same estimate to $x_m-x_k$ shows that each Fourier coefficient converges. The limiting coefficients have polynomial growth and therefore define a distribution $x$. For each $j$, $\Delta_{j}x_m\to\Delta_{j}x$ in $\mathcal{E}$. Fatou's inequality for the scalar sequence norm gives
\begin{equation*}
 \|x_m-x\|_{B_{\mathcal{E},r}^{s}}  \leq\liminf_{k\to\infty}  \|x_m-x_k\|_{B_{\mathcal{E},r}^{s}}  \to 0.
\end{equation*}
Thus $B_{\mathcal{E},r}^{s}$ is complete.

Let $\left\{\Lambda_{j}\right\}_{j\geq0}$ be another smooth inhomogeneous dyadic resolution. Finite overlap gives
\begin{equation*}
 \Lambda_{j}x  =\sum_{\substack{k\geq0\\|k-j|\leq N_0}}\Lambda_{j}\Delta_kx,
\end{equation*}
where $N_0$ also accounts for the low-frequency blocks. The uniform multiplier bounds imply
\begin{equation*}
 2^{js}\|\Lambda_{j}x\|_{\mathcal{E}}
 \leq C\sum_{\substack{k\geq0\\|k-j|\leq N_0}}
       2^{(j-k)s}2^{ks}\|\Delta_kx\|_{\mathcal{E}}.
\end{equation*}
Taking the $\ell_r$ norm and reversing the roles of the two resolutions proves equivalence of norms. This is the usual finite-overlap argument; compare \cite[Section~3.1]{XiongXuYin}.

Finally, applying the same argument to $x-P_{j}x=\sum_{j>J}\Delta_{j}x$ gives
\begin{equation*}
 \|x-P_{j}x\|_{B_{\mathcal{E},r}^{s}} \leq C\left\|       \bigl(2^{js}\|\Delta_{j}x\|_{\mathcal{E}}\bigr)_{j>J} \right\|_{\ell_r}.
\end{equation*}
The right-hand side tends to zero when $r<\infty$. Since $P_{j}x$ is a Fourier polynomial, this proves density. All multiplier and finite-overlap constants are independent of $\theta$.
\end{proof}

\begin{lem}[Reconstruction for positive smoothness]
\label{lem:besov-positive-reconstruction}
Let $s>0$ and $1\leq r\leq\infty$. For every $x\in B_{\mathcal{E},r}^{s}(\mathbb{T}_\theta^{d})$,
\begin{equation*}
 x=\sum_{j\geq0}\Delta_{j}^{\theta} x  \text{ with absolute convergence in }\mathcal{E}.
\end{equation*}
Consequently, $B_{\mathcal{E},r}^{s}\hookrightarrow\mathcal{E}$, and
\begin{equation}\label{eq:prelim-Besov-full-base}
 \|x\|_{B_{\mathcal{E},r}^{s}}
 \asymp
 \|x\|_{\mathcal{E}}
 +\left\|
   \bigl(2^{js}\|\Delta_{j}^{\theta} x\|_{\mathcal{E}}\bigr)_{j\geq1}
  \right\|_{\ell_r},
\end{equation}
with constants independent of $\theta$.
\end{lem}
\begin{proof} 
Let $r'$ be the conjugate index of $r$. Since $s>0$, H\"older's inequality gives
\begin{equation}\label{eq:symmetric-besov-absolute-reconstruction}
 \sum_{j\geq0}\|\Delta_{j}^{\theta} x\|_{\mathcal{E}}
 \leq
 \|(2^{-js})_{j\geq0}\|_{\ell_{r'}}
 \|x\|_{B_{\mathcal{E},r}^{s}}<\infty.
\end{equation}
Thus the reconstruction series converges absolutely in $\mathcal{E}$. Its sum agrees with $x$ in distributions, so
\begin{equation*}
 \|x\|_{\mathcal{E}}\leq C\|x\|_{B_{\mathcal{E},r}^{s}}.
\end{equation*}
Consequently, the expression on the right-hand side of \eqref{eq:prelim-Besov-full-base} is bounded by $C\|x\|_{B_{\mathcal{E},r}^{s}}$. Conversely, Lemma~\ref{lem:besov-dyadic-cutoffs} gives
\begin{align*}
 \|x\|_{B_{\mathcal{E},r}^{s}}
 &\leq\|\Delta_0^{\theta} x\|_{\mathcal{E}}
 +\left\|
   \bigl(2^{js}\|\Delta_{j}^{\theta} x\|_{\mathcal{E}}\bigr)_{j\geq1}
  \right\|_{\ell_r}\\
 &\leq C\|x\|_{\mathcal{E}}
 +\left\|
   \bigl(2^{js}\|\Delta_{j}^{\theta} x\|_{\mathcal{E}}\bigr)_{j\geq1}
  \right\|_{\ell_r}.
\end{align*}
This proves the norm equivalence, with constants independent
of $\theta$.
\end{proof}
\subsection{Embedding into the operator algebra}
\label{subsec:symmetric-besov-embedding}

\begin{lem}[Embedding criterion]
\label{lem:symmetric-besov-embedding}
Let $s>0$, $1\leq r\leq\infty$, and $\frac{1}{r}+\frac{1}{r'}=1$. Put
\begin{equation}\label{eq:symmetric-besov-embedding-weight}
 \omega_{j}(E,s)=\frac{2^{-js}}{\phi_E(2^{-jd})}, \ \ j\geq0.
\end{equation}
If $\omega(E,s)\in\ell_{r'}$, then
\begin{equation}\label{eq:besov-app-embedding}
 B_{\mathcal{E},r}^{s}(\mathbb{T}_\theta^{d})  \hookrightarrow L_\infty(\mathbb{T}_\theta^{d}),
\end{equation}
and
\begin{equation*}
 \|x\|_\infty \leq C \|\omega(E,s)\|_{\ell_{r'}}
       \|x\|_{B_{\mathcal{E},r}^{s}}.
\end{equation*}
In fact, the reconstruction series converges absolutely in operator norm.
\end{lem}

\begin{proof}
By \eqref{eq:besov-app-block} and H\"older's inequality,
\begin{align*}
 \sum_{j\geq0}\|\Delta_{j}x\|_\infty
 &\leq C\sum_{j\geq0}
       \frac{2^{-js}}{\phi_E(2^{-jd})}
       2^{js}\|\Delta_{j}x\|_{\mathcal{E}}\\
 &\leq C\|\omega(E,s)\|_{\ell_{r'}}
       \|x\|_{B_{\mathcal{E},r}^{s}}.
\end{align*}
The operator-norm limit agrees with $x$ in distributions.
\end{proof}

For $E=L_p$, $1\leq p<\infty$, one has $\phi_E(t)=t^{\frac{1}{p}}$, so $\omega_{j}(E,s)=2^{-j(s-\frac{d}{p})}$. The criterion therefore recovers the familiar sufficient conditions $s>\frac{d}{p}$, or $s=\frac{d}{p}$ and $r=1$. For $E=L_\infty$, it applies whenever $s>0$. Compare the classical Besov embeddings in \cite[Chapter~2]{BahouriCheminDanchin} and the quantum $L_p$ theory in \cite{XiongXuYin}.

\subsection{Ordered Bony decomposition and multiplication}
\label{subsec:symmetric-besov-multiplication}

For Fourier polynomials $a,b$, multiplication of Fourier modes gives
\begin{equation}\label{eq:symmetric-besov-product-support}
 \operatorname{supp}\widehat{ab}
 \subset\operatorname{supp}\widehat a+  \operatorname{supp}\widehat b.
\end{equation}
Indeed, the product of a mode of index $n$ and a mode of index $m$ is a scalar multiple of the mode of index $n+m$. Thus the classical support argument remains valid on a quantum torus, provided that the order of the factors is preserved.

Put
\begin{equation*}
 S_{j-1}=\sum_{0\leq k<j-1}\Delta_k, \ \ j\geq0,
\end{equation*}
with empty sums equal to zero. For $j\geq2$, this is $P_{j-2}$. We use the ordered paraproducts
\begin{align*}
 \mathcal{P}_{\mathrm L}(x,y)
 &=\sum_{j\geq0}(S_{j-1}x)\Delta_{j}y,\\
 \mathcal{P}_{\mathrm R}(x,y)
 &=\sum_{j\geq0}\Delta_{j}x(S_{j-1}y),\\
 \mathcal R(x,y)
 &=\sum_{j\geq0}\sum_{\substack{k\geq0\\|j-k|\leq1}}
                    \Delta_{j}x\Delta_ky.
\end{align*}
In particular, $\mathcal{P}_{\mathrm R}(x,y)$ is not defined by
interchanging $x$ and $y$ in $\mathcal{P}_{\mathrm L}(x,y)$.

\begin{lem}[Ordered Bony decomposition]
\label{lem:symmetric-besov-bony}
Let $s>0$ and $1\leq r\leq\infty$.
For $x,y\in B_{\mathcal{E},r}^{s}\cap L_\infty(\mathbb{T}_\theta^{d})$,
the three series above converge absolutely in $\mathcal{E}$, and
\begin{equation}\label{eq:symmetric-besov-bony}
 xy=\mathcal{P}_{\mathrm L}(x,y)
    +\mathcal{P}_{\mathrm R}(x,y)+\mathcal R(x,y).
\end{equation}
There is a fixed integer $N_0$ such that
\begin{align}
 \Delta_\ell\bigl((S_{j-1}x)\Delta_{j}y\bigr)
 &=\Delta_\ell\bigl(\Delta_{j}x(S_{j-1}y)\bigr)=0
 \ \ \text{if }|\ell-j|>N_0,
 \label{eq:symmetric-besov-low-high-support}\\
 \Delta_\ell(\Delta_{j}x\Delta_ky)&=0
 \ \text{ if } \ |j-k|\leq1 \ \text{ and } \ \ell>j+N_0.
 \label{eq:symmetric-besov-high-high-support}
\end{align}
These support bounds are independent of $\theta$.
\end{lem}

\begin{proof}
For Fourier polynomials, split the pairs of indices $(j,k)$ into $j<k-1$, $k<j-1$, and $|j-k|\leq1$. This gives \eqref{eq:symmetric-besov-bony} without interchanging any factors. It is the ordered form of Bony's decomposition \cite{Bony}; see also \cite[Chapter~2]{BahouriCheminDanchin}.

For $j\geq2$, the Fourier support of $S_{j-1}x$ is contained
in the ball of radius $2^j/3$. By
\eqref{eq:symmetric-besov-annuli} and
\eqref{eq:symmetric-besov-product-support}, both products in
\eqref{eq:symmetric-besov-low-high-support} have support in
\begin{equation*}
 \left\{n\in\mathbb Z^{d}:
       \frac1{24}\,2^j\leq|n|\leq\frac53\,2^j\right\}.
\end{equation*}
For $j=0,1$, these products vanish. If $|j-k|\leq1$, the support of $\Delta_{j}x\Delta_ky$ is contained in the ball of radius $4\cdot2^j$. The stated output restrictions follow.

For general $x,y$ as in the statement, uniform cutoff bounds give
\begin{align*}
 \sum_{j\geq0}\|(S_{j-1}x)\Delta_{j}y\|_{\mathcal{E}}
 &\leq C\|x\|_\infty
             \sum_{j\geq0}\|\Delta_{j}y\|_{\mathcal{E}},\\
 \sum_{j\geq0}\|\Delta_{j}x(S_{j-1}y)\|_{\mathcal{E}}
 &\leq C\|y\|_\infty
             \sum_{j\geq0}\|\Delta_{j}x\|_{\mathcal{E}},\\
 \sum_{j\geq0}\sum_{|k-j|\leq1}
       \|\Delta_{j}x\Delta_ky\|_{\mathcal{E}}
 &\leq C\|y\|_\infty
             \sum_{j\geq0}\|\Delta_{j}x\|_{\mathcal{E}}.
\end{align*}
In the last sum only nonnegative indices are used. All right-hand sides are finite by \eqref{eq:symmetric-besov-absolute-reconstruction}. Moreover, $P_{j}x\to x$ and $P_{j}y\to y$ in $\mathcal{E}$, while their operator norms remain uniformly bounded. Hence
\begin{align*}
 \|(P_{j}x)(P_{j}y)-xy\|_{\mathcal{E}} \leq \|P_{j}x-x\|_{\mathcal{E}}\|P_{j}y\|_\infty       +\|x\|_\infty\|P_{j}y-y\|_{\mathcal{E}}\to 0.
\end{align*}
Splitting the finite double sum for $(P_{j}x)(P_{j}y)$ and passing to the limit proves \eqref{eq:symmetric-besov-bony}.
\end{proof}

\begin{prop}[Multiplication and the algebra property]
\label{prop:besov-app-algebra}
Let $s>0$ and $1\leq r\leq\infty$. For
$x,y\in B_{\mathcal{E},r}^{s}(\mathbb{T}_\theta^{d})\cap L_\infty(\mathbb{T}_\theta^{d})$,
\begin{equation}\label{eq:besov-app-tame}
 \|xy\|_{B_{\mathcal{E},r}^{s}}
 \leq C\bigl(
       \|x\|_\infty\|y\|_{B_{\mathcal{E},r}^{s}}
       +\|y\|_\infty\|x\|_{B_{\mathcal{E},r}^{s}}\bigr).
\end{equation}
If $\omega(E,s)\in\ell_{r'}$, with $\omega$ defined in
\eqref{eq:symmetric-besov-embedding-weight}, then $B_{\mathcal{E},r}^{s}\hookrightarrow L_\infty$ and $B_{\mathcal{E},r}^{s}$ is a unital Banach $*$-algebra with an equivalent submultiplicative norm. The constants are independent of $\theta$.
\end{prop}

\begin{proof}
Set
\begin{equation*}
 b_{j}(x)=2^{js}\|\Delta_{j}x\|_{\mathcal{E}},
 \ \ b_{j}(y)=2^{js}\|\Delta_{j}y\|_{\mathcal{E}}.
\end{equation*}
By \eqref{eq:symmetric-besov-low-high-support}, the ideal property,
and the uniform multiplier bounds,
\begin{align*}
 2^{\ell s}\|\Delta_\ell\mathcal{P}_{\mathrm L}(x,y)\|_{\mathcal{E}}
 &\leq C\|x\|_\infty
       \sum_{|j-\ell|\leq N_0}2^{(\ell-j)s}b_{j}(y),\\
 2^{\ell s}\|\Delta_\ell\mathcal{P}_{\mathrm R}(x,y)\|_{\mathcal{E}}
 &\leq C\|y\|_\infty
       \sum_{|j-\ell|\leq N_0}2^{(\ell-j)s}b_{j}(x).
\end{align*}
For the remainder, \eqref{eq:symmetric-besov-high-high-support}
and $\|\Delta_ky\|_\infty\leq C\|y\|_\infty$ give
\begin{equation*}
 2^{\ell s}\|\Delta_\ell\mathcal R(x,y)\|_{\mathcal{E}}
 \leq C\|y\|_\infty
       \sum_{j\geq\max\left\{0,\ell-N_0\right\}}
                   2^{-(j-\ell)s}b_{j}(x).
\end{equation*}
The first two sequence kernels have finite support, and the
last belongs to $\ell_1$ because $s>0$. Discrete Young's
inequality, valid also for $r=1$ and $r=\infty$, proves
\eqref{eq:besov-app-tame} by
Lemma~\ref{lem:symmetric-besov-bony}.

Under the sequence condition,
Lemma~\ref{lem:symmetric-besov-embedding} gives
\begin{equation*}
 \|xy\|_{B_{\mathcal{E},r}^{s}}
 \leq C_A\|x\|_{B_{\mathcal{E},r}^{s}}
             \|y\|_{B_{\mathcal{E},r}^{s}}.
\end{equation*}
Multiplying the Besov norm by $\max\left\{1,C_A\right\}$ makes it
submultiplicative. Completeness follows from
Proposition~\ref{lem:besov-app-basic}.
Since $\Delta_0\mathbf{1}=\mathbf{1}$ and
$\Delta_{j}\mathbf{1}=0$ for $j\geq1$, the unit belongs to the space.
Finally, the cutoffs are real and even, so
\begin{equation*}
 \Delta_{j}(x^*)=(\Delta_{j}x)^*,
 \ \ \|x^*\|_{B_{\mathcal{E},r}^{s}}
       =\|x\|_{B_{\mathcal{E},r}^{s}}.
\end{equation*}
Thus the involution is isometric.
\end{proof}

\begin{cor}[Multiplication by regular coefficients]
\label{cor:besov-app-coefficient}
Let $\rho>s>0$ and $1\leq r\leq\infty$. If $a\in B_{\infty,\infty}^{\rho}(\mathbb{T}_\theta^{d})$, then, for every $x\in B_{\mathcal{E},r}^{s}$,
\begin{equation}\label{eq:besov-app-coefficient}
 \|ax\|_{B_{\mathcal{E},r}^{s}}
 +\|xa\|_{B_{\mathcal{E},r}^{s}}
 \leq C\|a\|_{B_{\infty,\infty}^{\rho}}
       \|x\|_{B_{\mathcal{E},r}^{s}}.
\end{equation}
Moreover, $a\in B_{\mathcal{E},r}^{s}\cap L_\infty$. The constants are independent of $\theta$.
\end{cor}

\begin{proof}
Write $A=\|a\|_{B_{\infty,\infty}^{\rho}}$. Then
\begin{equation*}
 \|\Delta_{j}a\|_\infty\leq2^{-j\rho}A,
 \ \ \|a\|_\infty\leq\sum_{j\geq0}2^{-j\rho}A\leq C_\rho A.
\end{equation*}
In particular, $P_{j}a\to a$ in operator norm. Since $x\in\mathcal{E}$ and $P_{j}x\to x$ in $\mathcal{E}$, the finite ordered decompositions of $(P_{j}a)(P_{j}x)$ and $(P_{j}x)(P_{j}a)$ pass to the limit in $\mathcal{E}$.

The terms with low-frequency $a$, and the terms with comparable frequencies, are estimated as in Proposition~\ref{prop:besov-app-algebra}, with $\|a\|_\infty\leq C_\rho A$. For the terms with low-frequency $x$, we use
\begin{equation*}
 \|S_{j-1}x\|_{\mathcal{E}}
 \leq C\|x\|_{\mathcal{E}}
 \leq C\|x\|_{B_{\mathcal{E},r}^{s}}.
\end{equation*}
Thus, in both orders of multiplication,
\begin{equation*}
 \begin{aligned}
 2^{js}\|(S_{j-1}x)\Delta_{j}a\|_{\mathcal{E}}
 &\leq C2^{-j(\rho-s)}A\|x\|_{B_{\mathcal{E},r}^{s}},\\
 2^{js}\|\Delta_{j}a(S_{j-1}x)\|_{\mathcal{E}}
 &\leq C2^{-j(\rho-s)}A\|x\|_{B_{\mathcal{E},r}^{s}}.
 \end{aligned}
\end{equation*}
Finite overlap and the summability of
$(2^{-j(\rho-s)})_{j\geq0}$ prove
\eqref{eq:besov-app-coefficient}.
These estimates also justify absolute convergence of the
ordered series in $\mathcal{E}$ when $x$ is not bounded.

Finally, $\|\Delta_{j}a\|_{\mathcal{E}}
\leq\phi_E(1)\|\Delta_{j}a\|_\infty$, so
\begin{equation*}
 \|a\|_{B_{\mathcal{E},r}^{s}}
 \leq\phi_E(1)
       \|(2^{-j(\rho-s)})_{j\geq0}\|_{\ell_r}
       \|a\|_{B_{\infty,\infty}^{\rho}}.
\end{equation*}
\end{proof}

\subsection{Besov--Lorentz spaces}
\label{subsec:symmetric-besov-lorentz}

For $0<p<\infty$, $0<q\leq\infty$, $s\in\mathbb{R}$, and
$1\leq\rho\leq\infty$, define the Besov--Lorentz quasi-norm by
\begin{equation}\label{eq:prelim-Besov}
 \|x\|_{B_{p,q;\rho}^{s}}
 =\left\|\bigl(2^{js}\|\Delta_{j}^{\theta} x\|_{L_{p,q}}\bigr)_{j\geq0}
   \right\|_{\ell_\rho}.
\end{equation}
Here the Lorentz quasi-norm is defined using the singular value
function. If $1<p<\infty$, $1\leq q\leq\infty$, or
$(p,q)=(1,1)$, then $L_{p,q}$ admits an equivalent fully
symmetric Banach norm; see \cite{BennettSharpley} and
\cite{DoddsDePagterSukochev}.
In this Banach range, \eqref{eq:prelim-Besov} agrees, up to
equivalent quasi-norms, with the symmetric Besov construction
above. In particular, these Besov--Lorentz spaces are Banach spaces
after an equivalent renorming. We retain
\eqref{eq:prelim-Besov} for the estimates, and write
$B_{p,\rho}^{s}=B_{p,p;\rho}^{s}$.

For $s>0$ and indices in the Banach range,
Proposition~\ref{lem:besov-app-basic} yields
\begin{equation*}
 \|x\|_{B_{p,q;\rho}^{s}}
 \asymp\|x\|_{L_{p,q}}
 +\left\|\bigl(2^{js}\|\Delta_{j}^{\theta} x\|_{L_{p,q}}\bigr)_{j\geq1}
   \right\|_{\ell_\rho}.
\end{equation*}
The fundamental function of a fully symmetric Banach realization
of $L_{p,q}(0,1)$ is comparable to $t^{\frac{1}{p}}$.
Proposition~\ref{lem:besov-app-basic} and
Lemma~\ref{lem:symmetric-besov-embedding}, together with
Proposition~\ref{prop:besov-app-algebra}, therefore give the
following consequence.

\begin{cor}\label{cor:symmetric-besov-lorentz-algebra}
Let $1<p<\infty$, $1\leq q\leq\infty$, or $(p,q)=(1,1)$,
and let $1\leq\rho\leq\infty$.
Then $B_{p,q;\rho}^{s}$ embeds continuously
into $L_\infty(\mathbb{T}_\theta^{d})$ and is a unital Banach
$*$-algebra after an equivalent renorming whenever
\begin{equation*}
 s>\frac dp,
 \ \ \text{or}\ \  s=\frac dp\text{ and }\rho=1.
\end{equation*}
The embedding and multiplication constants are independent
of $\theta$.
\end{cor}

The constants throughout this section may depend on the dimension,
the chosen cutoffs, the smoothness and summability indices, and,
where applicable, the underlying function space. They do not
depend on the deformation parameter. The algebra criterion and
the coefficient estimate \eqref{eq:besov-app-coefficient} will
be used in the final application section.

\begin{rem}[Approximation in the two scales]
Sobolev polynomial approximation above uses separability of $E$.
For Besov spaces, finite dyadic blocks and the $\ell_r$ tail imply
approximation for $r<\infty$ without this hypothesis. At $r=\infty$,
polynomial density in the full Besov space is not asserted.
\end{rem}

\section{Convolution inequalities in function spaces}
\label{sec:applications}\label{sec:conv-function-spaces}

In this section, we apply the sharp O'Neil inequality to convolution in symmetric operator spaces. The main idea is to reduce the operator estimates to inequalities for decreasing rearrangements, as in the classical theory~\cite{ONeil,BennettSharpley}. We obtain Young--O'Neil inequalities in Lorentz spaces, logarithmic limiting estimates of Br\'ezis--Wainger type, and convolution estimates in Besov--Lorentz spaces; see also~\cite{BW,Krepela,KuhnSchilling}.

Throughout this section, $\theta_{1}$ and $\theta_{2}$ are deformation parameters.  All constants are independent of $\theta_{1}$ and $\theta_{2}$.

\subsection{Symmetric-space transfer and Hardy reduction}

For nonnegative measurable functions $f$ and $g$ on $(0,1)$, define the bilinear O'Neil majorant by
\begin{equation}\label{eq:ONeil-majorant}
\mathscr O(f,g)(t) = t(\C f)(t)(\C g)(t)
+ \int\limits_{t}^{1}f(s)g(s)\,ds, \ \  0<t<1.
\end{equation}

\begin{thm}[O'Neil transfer principle]
\label{thm:ONeil-transfer}
Let $(x,y,z)$ satisfy the O'Neil relation. Then
\begin{equation}\label{eq:pointwise-ONeil}
\mu(t,z)\leq\mathscr O(\mu(x),\mu(y))(t), \ \  0<t<1.
\end{equation}
Let $E$, $F$, and $G$ be symmetric Banach function spaces on $(0,1)$ such that $E,F\subset L_{1}(0,1)$, and suppose that
\begin{equation}\label{eq:scalar-ONeil-bound}
\|\mathscr O(f,g)\|_{G} \leq C_{E,F,G}\|f\|_{E}\|g\|_{F}
\end{equation}
for all nonnegative decreasing $f\in E$ and $g\in F$. If $\mu(x)\in E$ and $\mu(y)\in F$, then $\mu(z)\in G$ and
\begin{equation}\label{eq:abstract-symmetric-transfer}
\|\mu(z)\|_G \leq C_{E,F,G}\|\mu(x)\|_E\|\mu(y)\|_F.
\end{equation}
\end{thm}

\begin{proof}
Dividing \eqref{eq:ONeil-relation} by $t$ gives
\begin{equation*}
\mathcal{C}\mu(z)(t)
\leq\mathscr O(\mu(x),\mu(y))(t).
\end{equation*}
Since $\mu(z)$ is decreasing, $\mu(t,z)\leq\mathcal{C}\mu(z)(t)$, which proves \eqref{eq:pointwise-ONeil}. The ideal property of $G$ and \eqref{eq:scalar-ONeil-bound} then give
\begin{equation*}
\|\mu(z)\|_G
\leq\|\mathscr O(\mu(x),\mu(y))\|_G
\leq C_{E,F,G}\|\mu(x)\|_E\|\mu(y)\|_F.
\end{equation*}
\end{proof}

For a symmetric Banach function space $E$ on $(0,1)$, we continue to write
\begin{equation*}
\mathcal{E}(\mathbb{T}_{\theta}^{d})
=\left\{x\in S(\mathbb{T}_{\theta}^{d}):\mu(x)\in E\right\},
\ \ 
\|x\|_{\mathcal{E}}=\|\mu(x)\|_E,
\end{equation*}
and use analogous notation for $F$ and $G$.  By Theorem~\ref{thm:operator-oneil}, the triple
\begin{equation}\label{eq:convolution-ONeil-triple}
\bigl(x,y,x\star_{\theta_1,\theta_2}y\bigr)
\end{equation}
satisfies the O'Neil relation for every pair of $L_{1}$ inputs.  The same is true of $(f,x,f\star_{0,\theta} x)$ when $\theta_1=0$.  Thus all function-space consequences below follow by applying Theorem~\ref{thm:ONeil-transfer} to the appropriate triple.

\begin{corollary}[Lorentz consequence of the O'Neil relation]
\label{cor:ONeil-Lorentz-transfer}
Let $(x,y,z)$ satisfy the O'Neil relation. Suppose that $1<p_1,p_2,p<\infty$ and
\begin{equation*}
1+\frac1p=\frac1{p_1}+\frac1{p_2},
\end{equation*}
and let $1\leq q_{1},q_{2},q\leq\infty$ satisfy
\begin{equation*}
\frac1q\leq\frac1{q_{1}}+\frac1{q_{2}}.
\end{equation*}
Then
\begin{equation}\label{eq:abstract-Lorentz-ONeil}
\|z\|_{L_{p,q}(\mathcal M_3)}
\leq C
\|x\|_{L_{p_1,q_{1}}(\mathcal M_1)}
\|y\|_{L_{p_2,q_{2}}(\mathcal M_2)},
\end{equation}
where $C$ depends only on the indices and on the chosen equivalent Lorentz norms.
\end{corollary}

\begin{proof}
The Hardy inequalities used in the proof of \cite[Theorem~2.6]{ONeil} give
\begin{equation*}
\|\mathscr O(f,g)\|_{L_{p,q}(0,1)}
\leq C
\|f\|_{L_{p_1,q_{1}}(0,1)}
\|g\|_{L_{p_2,q_{2}}(0,1)}
\end{equation*}
for nonnegative decreasing $f$ and $g$. Apply Theorem~\ref{thm:ONeil-transfer} with $E=L_{p_1,q_{1}}$, $F=L_{p_2,q_{2}}$, and $G=L_{p,q}$.
\end{proof}

\subsection{Young--O'Neil inequalities}\label{subsec:young-oneil}
We use the Lorentz quasi-norms and equivalent Banach norms of Section~\ref{sec:preliminaries}. The sharp Young estimate \eqref{eq:Young} is already proved in Part I \cite{ConvolutionI}; the refinement below uses the O'Neil relation.

\begin{corollary}[Young--O'Neil inequality]
\label{cor:Young-ONeil}
Let $1<p_{1},p_{2},p<\infty$ satisfy
\begin{equation*}
1+\frac{1}{p}
=
\frac{1}{p_{1}}+\frac{1}{p_{2}},
\end{equation*}
and let $1\leq q_{1},q_{2},q\leq\infty$ satisfy
\begin{equation*}
\frac{1}{q}
\leq
\frac{1}{q_{1}}+\frac{1}{q_{2}}.
\end{equation*}
Then
\begin{align*}
\left\|
x\star_{\theta_{1},\theta_{2}}y
\right\|_{L_{p,q}(\mathbb{T}_{\theta_{1}+\theta_{2}}^{d})}
\leq
C\,
\|x\|_{L_{p_{1},q_{1}}(\mathbb{T}_{\theta_{1}}^{d})}
\|y\|_{L_{p_{2},q_{2}}(\mathbb{T}_{\theta_{2}}^{d})},
\end{align*}
where $C$ depends only on the indices and on the chosen equivalent
Lorentz norms.  In particular, $C$ is independent of
$\theta_{1}$ and $\theta_{2}$.
\end{corollary}

\begin{proof}
The triple in \eqref{eq:convolution-ONeil-triple} satisfies the O'Neil
relation. Apply Corollary~\ref{cor:ONeil-Lorentz-transfer}.
\end{proof}

\subsection{Limiting targets and logarithmic estimates}
\label{subsec:working-weak-Linfty}

We use the logarithmic classes introduced in Subsection~\ref{subsec:logarithmic-limiting-classes}. With the endpoint conventions stated there, one has
\begin{equation}\label{eq:Linfty-q-into-L1}
 L_{\infty,q}(0,1)\hookrightarrow L_1(0,1),  \ \ 1\leq q\leq\infty.
\end{equation}
Indeed, let $1<q\leq\infty$, with $q'=1$ when $q=\infty$. Writing $\ell(t)=\log\left( e/t \right)$, H\"older's inequality gives
\begin{align*}
 \|x\|_{L_1(0,1)}
 &=\int\limits_{0}^{1}
   \frac{\mu(t,x)}{\ell(t)}\,t\ell(t)\,\frac{dt}{t}\\
 &\leq
 \|x\|_{L_{\infty,q}(0,1)}
 \left(
 \int\limits_{0}^{1}[t\ell(t)]^{q'}\frac{dt}{t}
 \right)^{\frac{1}{q}'}.
\end{align*}
The last integral is finite, including at $q=\infty$, where it equals $2$. For $q=1$, our convention $L_{\infty,1}=L_\infty$ gives $\|x\|_{L_1(0,1)}\leq\|x\|_{L_\infty(0,1)}$.

\begin{thm}\label{thm:NT31-NC}
Let $\theta_1,\theta_2$ be real skew-symmetric $d\times d$
matrices, and let $1\leq q_1,q_2\leq\infty$ satisfy
\begin{equation*}
 \frac1{q_1}+\frac1{q_2}=1.
\end{equation*}
Suppose that
$x\in L_{\infty,q_1}(\mathbb{T}_{\theta_1}^{d})$ and
$y\in L_1(\mathbb{T}_{\theta_2}^{d})$ satisfy
$\mathcal{C}\mu(y)\in L_{1,q_2}(0,1)$.
Then
$x\star_{\theta_1,\theta_2}y
 \in L_\infty(\mathbb{T}_{\theta_1+\theta_2}^{d})$ and
\begin{equation*}
 \|x\star_{\theta_1,\theta_2}y\|_
 {L_\infty(\mathbb{T}_{\theta_1+\theta_2}^{d})}
 \leq C_{q_1}
 \|x\|_{L_{\infty,q_1}(\mathbb{T}_{\theta_1}^{d})}
 \|\mathcal{C}\mu(y)\|_{L_{1,q_2}(0,1)}.
\end{equation*}
The constant is independent of $\theta_1$ and $\theta_2$.
\end{thm}
\begin{lem}\label{37 lemma}
Let $1\leq q_{1},q_{2}\leq\infty$ satisfy $\frac1{q_{1}}+\frac1{q_{2}}=1.$ Let $x\in L_{\infty,q_{1}}(0,1)$ and let $f\in L_{1}(0,1)$ satisfy $\mathcal{C}\mu(f)\in L_{1,q_{2}}(0,1).$ Then
\begin{equation*}
\int\limits_{0}^{1}\mu(s,f)\mu(s,x)\,ds
\leq
C_{q_{1}}
\|x\|_{L_{\infty,q_{1}}(0,1)}
\|\mathcal{C}\mu(f)\|_{L_{1,q_{2}}(0,1)}.
\end{equation*}
\end{lem}
 
\begin{proof}
First assume that $1<q_{1},q_{2}<\infty$. Since $q_{1}=q_{2}'$, introduce the Lorentz--Gamma space
\begin{equation*}
Y_{q_{2}}(0,1) = \left\{ f: \|\mathcal{C}\mu(f)\|_{L_{1,q_{2}}(0,1)}<\infty \right\}.
\end{equation*}
Its norm can be written as
\begin{equation*}
\|f\|_{Y_{q_{2}}}^{q_{2}}
= \int\limits_{0}^{1} \bigl((\mathcal{C}\mu(f))(t)\bigr)^{q_{2}} t^{q_{2}-1}\,dt.
\end{equation*}
Thus, $Y_{q_{2}}$ is the Lorentz--Gamma space with exponent $q_{2}$ and weight $\Phi(t)=t^{q_{2}-1}$. For this weight,
\begin{equation*}
\int\limits_{0}^{t}\Phi(s)\,ds = \frac{t^{q_{2}}}{q_{2}},
\ \  \int\limits_{t}^{1}s^{-q_{2}}\Phi(s)\,ds = \log\left(\frac{1}{t}\right).
\end{equation*}
Substitution into the associate-weight formula in \cite[Theorem~5.2]{KP2011} gives
\begin{equation*}
\Psi(t) \asymp_{q_{1}}
\frac{1}{t\log^{q_{1}}\left( e/t \right)},
\ \  0<t<1,
\end{equation*}
where, on the interval bounded away from zero, the admissible extension of $\Psi$ is comparable with the displayed weight. Consequently,
\begin{equation}\label{eq:KP-associate-Gamma-norm}
\|x\|_{Y_{q_{2}}'}
\asymp_{q_{1}}
\left(
\int\limits_{0}^{1}
\left[
\frac{(\mathcal{C}\mu(x))(t)}
{\log\left( e/t \right)}
\right]^{q_{1}}
\frac{dt}{t}
\right)^{\frac1{q_{1}}}.
\end{equation}

We claim that the expression on the right-hand side of
\eqref{eq:KP-associate-Gamma-norm} is equivalent to
$\|x\|_{L_{\infty,q_{1}}(0,1)}$. One direction follows immediately from
$\mu(t,x)\leq(\mathcal{C}\mu(x))(t)$. For the converse, set
\begin{equation*}
t=e^{1-u},
\ \ 
X(u)=\frac{\mu(e^{1-u},x)}{u},
\ \  u>1.
\end{equation*}
Then
\begin{equation*}
\frac{(\mathcal{C}\mu(x))(e^{1-u})}{u}
\leq
\int\limits_{0}^{\infty}e^{-v}(1+v)X(u+v)\,dv.
\end{equation*}
Since $e^{-v}(1+v)\in L_{1}(0,\infty)$, Young's inequality yields
\begin{equation*}
\left( \int\limits_{0}^{1}
\left[ \frac{(\mathcal{C}\mu(x))(t)}{\log\left( e/t \right)} \right]^{q_{1}}
\frac{dt}{t} \right)^{\frac1{q_{1}}}
\leq C_{q_{1}}\|x\|_{L_{\infty,q_{1}}(0,1)}.
\end{equation*}
Therefore,
\begin{equation*}
Y_{q_{2}}(0,1)' = L_{\infty,q_{1}}(0,1)
\end{equation*}
with equivalence of norms. K\"othe duality now gives
\begin{equation*}
\int\limits_{0}^{1}\mu(s,f)\cdot\mu(s,x)\,ds \leq
C_{q_{1}} \|x\|_{L_{\infty,q_{1}}(0,1)} \|\mathcal{C}\mu(f)\|_{L_{1,q_{2}}(0,1)}.
\end{equation*}

If $(q_{1},q_{2})=(1,\infty)$, then $L_{\infty,1}=L_{\infty}$ and
\begin{equation*}
\|\mathcal{C}\mu(f)\|_{L_{1,\infty}(0,1)}
= \sup_{0<t<1}t(\mathcal{C}\mu(f))(t)
= \|f\|_{L_{1}(0,1)}.
\end{equation*}
Hence the assertion follows from the usual $L_{\infty}$--$L_{1}$ pairing.

Finally, if $(q_{1},q_{2})=(\infty,1)$, it follows from
\begin{equation*}
(\operatorname{Exp}L)'=L\log L, \ \ 
\|f\|_{L\log L} \asymp \|\mathcal{C}\mu(f)\|_{L_{1}(0,1)},
\end{equation*}
together with $L_{\infty,\infty}=\operatorname{Exp}L$. This completes the proof.
\end{proof}

\begin{lem}\label{38 easy lemma}
For nonnegative nonincreasing functions $a,b$ on $(0,1)$,
\begin{equation}\label{eq:Chebyshev-endpoint-product}
 \frac{1}{t}  \int\limits_0^t a(s)\,ds\cdot \int\limits_0^t b(s)\,ds  \leq \int\limits_0^t a(s)b(s)\,ds,
 \ \ 0<t<1.
\end{equation}
\end{lem}

\begin{proof}
Apply Chebyshev's integral inequality on $(0,t)$ with respect to the normalized measure $ds/t$. For unbounded functions, apply the inequality to their bounded truncations and use monotone convergence.
\end{proof}

\begin{proof}[Proof of Theorem~\ref{thm:NT31-NC}]
By \eqref{eq:Linfty-q-into-L1}, applied to $\mu(x)$, we have $x\in L_1(\mathbb{T}_{\theta_1}^{d})$. Hence the convolution $z=x\star_{\theta_1,\theta_2}y$ is well defined in $L_1(\mathbb{T}_{\theta_1+\theta_2}^{d})$. Theorem~\ref{thm:operator-oneil} and Lemma~\ref{38 easy lemma} give, for $0<t<1$,
\begin{align*}
 (\mathcal{C}\mu(z))(t)
 &\leq
 \frac{1}{t}
 \int\limits_0^t\mu(s,x)\,ds\cdot \int\limits_0^t\mu(s,y)\,ds 
 +\int\limits_t^1\mu(s,x) \cdot \mu(s,y)\,ds\\
 &\leq
 \int\limits_{0}^{1}\mu(s,x)\cdot\mu(s,y)\,ds \leq C_{q_1}
 \|x\|_{L_{\infty,q_1}(\mathbb{T}_{\theta_1}^{d})}
 \|\mathcal{C}\mu(y)\|_{L_{1,q_2}(0,1)},
\end{align*}
where the last step follows from Lemma~\ref{37 lemma}, applied to the  functions $\mu(x)$ and $\mu(y)$. Letting $t\downarrow0$ proves the assertion, since
\begin{equation*}
 \lim_{t\downarrow0}(\mathcal{C}\mu(z))(t)=\|z\|_\infty,
\end{equation*}
with the equality initially understood in $[0,\infty]$.
\end{proof}

\subsubsection{Br\'ezis--Wainger estimates}

Motivated by the limiting Sobolev and convolution inequalities of Br\'ezis and Wainger~\cite{BW}, we establish convolution estimates on noncommutative tori with targets in the logarithmic Lorentz--Zygmund classes introduced in Subsection~\ref{subsec:logarithmic-limiting-classes}.

Recall that
\begin{equation*}
 \ell(t)=\log\left( e/t \right),
 \ \ 
 (\mathcal{C}a)(t)=\frac{1}{t}\int\limits_0^t a(s)\,ds,
 \ \ 0<t<1.
\end{equation*}
For $1<q\leq\infty$, we write
\begin{align*}
 \|h\|_{L_{\infty,q}(0,1)}
 &=\left\|t\mapsto\frac{\mu(t,h)}{\ell(t)}
   \right\|_{L_q((0,1),\frac{dt}{t})},\\
 \|h\|_{L_{1,q}(0,1)}
 &=\left\|t\mapsto t\mu(t,h)
   \right\|_{L_q((0,1),\frac{dt}{t})}.
\end{align*}
At $q=\infty$, the right-hand sides are understood as essential suprema. For operators on noncommutative tori, the same formulas apply with $\mu(t,h)$ denoting the singular value function.

\begin{thm}\label{thm:working-Brezis-Wainger}
Let $\theta_1$ and $\theta_2$ be real antisymmetric $d\times d$ matrices. Let $1<q_1<\infty$ and $1<q_2,q_3\leq\infty$ satisfy
\begin{equation*}
 \frac1{q_1}=\frac1{q_2}+\frac1{q_3},
 \ \ \frac1\infty=0.
\end{equation*}
Then the following assertions hold.
\begin{enumerate}
 \item[(i)] If $x\in L_{\infty,q_2}(\mathbb{T}^{d}_{\theta_2})$ and
 $f\in L_{1}(\mathbb{T}^{d}_{\theta_1})$ satisfy $\mathcal{C}\mu(f)\in L_{1,q_3}(0,1)$, then  $f\star_{\theta_1,\theta_2}x \in L_{\infty,q_1}(\mathbb{T}^{d}_{\theta_1+\theta_2})$ and
 \begin{equation*}
  \|f\star_{\theta_1,\theta_2}x\|_{L_{\infty,q_1}
     (\mathbb{T}^{d}_{\theta_1+\theta_2})}
  \leq C_{q_1,q_2,q_3}
  \|x\|_{L_{\infty,q_2}(\mathbb{T}^{d}_{\theta_2})}
  \|\mathcal{C}\mu(f)\|_{L_{1,q_3}(0,1)}.
 \end{equation*}

 \item[(ii)] Let $1<p<\infty$ and $p'=p/(p-1)$. If
 $x\in L_{p,q_2}(\mathbb{T}^{d}_{\theta_2})$ and
 $f\in L_{p',q_3}(\mathbb{T}^{d}_{\theta_1})$, then
 $f\star_{\theta_1,\theta_2}x
 \in L_{\infty,q_1}(\mathbb{T}^{d}_{\theta_1+\theta_2})$ and
 \begin{equation*}
  \|f\star_{\theta_1,\theta_2}x\|_{L_{\infty,q_1}
     (\mathbb{T}^{d}_{\theta_1+\theta_2})}
  \leq C_{p,q_1,q_2,q_3}
  \|x\|_{L_{p,q_2}(\mathbb{T}^{d}_{\theta_2})}
  \|f\|_{L_{p',q_3}(\mathbb{T}^{d}_{\theta_1})}.
 \end{equation*}
\end{enumerate}
\end{thm}
Taking $\theta_1=0$ gives the corresponding semi-convolution estimates. We first establish the scalar estimates needed in the proof.

\begin{lem}[Logarithmic Hardy estimate]
\label{lem:logarithmic-Hardy-estimate}
Let $1<p\leq\infty$, and let $a$ be a nonnegative measurable function on $(0,1)$. Then
\begin{equation}\label{eq:logarithmic-Hardy-estimate}
 \left\|\frac{\mathcal{C}a}{\ell}\right\|_{L_p((0,1),\frac{dt}{t})}
 \leq 2 \left\|\frac{a}{\ell}\right\|_{L_p((0,1),\frac{dt}{t})}.
\end{equation}
\end{lem}

\begin{proof}
We may assume that the right-hand side is finite. Set
\begin{equation*}
 t=e^{1-u}, \ \ Z(u)=\frac{a(e^{1-u})}{u}, \ \ u>1.
\end{equation*}
Since $\ell(e^{1-u})=u$ and $\frac{dt}{t}=-du$, we have
\begin{equation*}
 \|Z\|_{L_p(1,\infty)} = \left\|\frac{a}{\ell}\right\|_{L_p((0,1),\frac{dt}{t})}.
\end{equation*}
The substitution $s=e^{1-u-v}$ gives
\begin{align*}
 \frac{(\mathcal{C}a)(e^{1-u})}{u} & =\frac1u\int\limits_0^\infty e^{-v}a(e^{1-u-v})\,dv =\int\limits_0^\infty  e^{-v}\left(1+\frac vu\right)Z(u+v)\,dv \\
 & \leq\int\limits_0^\infty e^{-v}(1+v)Z(u+v)\,dv.
\end{align*}
For $p<\infty$ and $v\geq0$,
\begin{equation*}
 \|Z(\,\cdot+v)\|_{L_p(1,\infty)}^p = \int\limits_{1+v}^\infty |Z(r)|^p\,dr  \leq\|Z\|_{L_p(1,\infty)}^p.
\end{equation*}
The same contraction estimate holds for $p=\infty$. Thus Minkowski's integral inequality, or the essential supremum estimate when $p=\infty$, yields
\begin{align*}
 \left\|u\mapsto\frac{(\mathcal{C}a)(e^{1-u})}{u}\right\|_{L_p(1,\infty)}
 &\leq\int\limits_0^\infty e^{-v}(1+v)
       \left\|u\mapsto Z(u+v)\right\|_{L_p(1,\infty)}\,dv\\
 &\leq2\|Z\|_{L_p(1,\infty)}.
\end{align*}
Returning to the variable $t$ proves the assertion.
\end{proof}

For $1<q<\infty$ and a measurable function $h$ on $(0,1)$, put
\begin{equation*}
 \|h\|_{N_q} = \left(\int\limits_{0}^{1}
   \left[\int\limits_0^t\frac{\mu(s,h)}{\ell(s)}\,ds\right]^q
   \frac{dt}{t}\right)^{\frac{1}{q}}.
\end{equation*}
The following estimate controls the endpoint value of the inner integral. Here and below, $c_{\mathrm{abs}}$ denotes an absolute constant, whose value may change from line to line.

\begin{lem}\label{nq vs llogl}
Let $1<q<\infty$. For every measurable function $h$ on $(0,1)$,
\begin{equation*}
 \int\limits_{0}^{1}\frac{\mu(s,h)}{\ell(s)}\,ds
 \leq c_{\mathrm{abs}}\|h\|_{N_q}.
\end{equation*}
\end{lem}

\begin{proof}
We may assume that $\|h\|_{N_q}<\infty$. By positivity,
\begin{align*}
 \|h\|_{N_q}^q \geq \int\limits_{e^{-1}}^1
   \left[\int\limits_0^{e^{-1}}\frac{\mu(s,h)}{\ell(s)}\,ds\right]^q
   \frac{dt}{t} =\left[\int\limits_0^{e^{-1}}\frac{\mu(s,h)}{\ell(s)}\,ds\right]^q.
\end{align*}
Since $\mu(\,\cdot,h)$ is nonincreasing,
\begin{align*}
 \int\limits_0^{e^{-1}}\frac{\mu(s,h)}{\ell(s)}\,ds
\geq\mu(e^{-1},h)\int\limits_0^{e^{-1}}\frac{ds}{\ell(s)},
 \end{align*}
 and
 \begin{align*}
 \int\limits_{e^{-1}}^1\frac{\mu(s,h)}{\ell(s)}\,ds
\leq\mu(e^{-1},h)\int\limits_{e^{-1}}^1\frac{ds}{\ell(s)}.
\end{align*}
Combining these inequalities gives
\begin{equation*}
 \int\limits_{0}^{1}\frac{\mu(s,h)}{\ell(s)}\,ds
 \leq
 \left(1+
 \frac{\int\limits_{e^{-1}}^1\ell(s)^{-1}\,ds}
      {\int\limits_0^{e^{-1}}\ell(s)^{-1}\,ds}\right)
 \|h\|_{N_q},
\end{equation*}
as required.
\end{proof}

\begin{lem}\label{lem:tail-log-equivalence}
Let $1<q<\infty$ and $q'=\frac{q}{q-1}$. For a measurable function $h$ on $(0,1)$, define
\begin{equation*}
 (\mathfrak Th)(t)=\int\limits_t^1\mu(s,h)\,ds,
 \ \ 0<t<1.
\end{equation*}
Then
\begin{equation*}
 q^{-1}\|h\|_{N_q}
 \leq\|\mathfrak Th\|_{L_{\infty,q}(0,1)}
 \leq c_{\mathrm{abs}}q'\|h\|_{N_q}.
\end{equation*}
\end{lem}
\begin{proof}
Let $\mathcal{C}_\infty$ denote the Hardy averaging operator on $(0,\infty)$, extending $\mathcal{C}$ on $(0,1)$, and let $\mathcal{C}_\infty^*$ denote its associate:
\begin{equation*}
 (\mathcal{C}_\infty a)(u)=\frac1u\int\limits_0^u a(v)\,dv,
 \ \ 
 (\mathcal{C}_\infty^*a)(u)=\int\limits_u^\infty\frac{a(v)}v\,dv,
 \ \ u>0.
\end{equation*}
For nonnegative $a$, Tonelli's theorem gives
\begin{equation*}
 \mathcal{C}_\infty\mathcal{C}_\infty^*a
 =\mathcal{C}_\infty a+\mathcal{C}_\infty^*a
 =\mathcal{C}_\infty^*\mathcal{C}_\infty a.
\end{equation*}
Combining the preceding identities with Hardy's inequalities \cite[Chapter~III, Lemma~3.9, p.~124]{BennettSharpley}, we obtain
\begin{align}
 \|\mathcal{C}_\infty a\|_{L_q(0,\infty)}
 &\leq q'\|\mathcal{C}_\infty^*a\|_{L_q(0,\infty)},
 \label{eq:BW-Hardy-comparison-forward}\\
 \|\mathcal{C}_\infty^*a\|_{L_q(0,\infty)}
 &\leq q\|\mathcal{C}_\infty a\|_{L_q(0,\infty)}.
 \label{eq:BW-Hardy-comparison-backward}
\end{align}

Set $a(u)=e^{1-u}\mu(e^{1-u},h)$ for $u>1$, and extend $a$ by zero to $(0,1]$. Substituting $s=e^{1-v}$ yields
\begin{equation*}
 (\mathcal{C}_\infty a)(u) = \frac1u\int\limits_{e^{1-u}}^1\mu(s,h)\,ds, \ \ u>1,
\end{equation*}
whereas $\mathcal{C}_\infty a=0$ on $(0,1]$. Since $\mathfrak Th$ is nonnegative and nonincreasing, the change of variables $t=e^{1-u}$ gives
\begin{equation*}
 \|\mathcal{C}_\infty a\|_{L_q(0,\infty)} = \|\mathfrak Th\|_{L_{\infty,q}(0,1)}.
\end{equation*}

Put
\begin{equation*}
 J(t)=\int\limits_0^t\frac{\mu(s,h)}{\ell(s)}\,ds, \ \ 0<t\leq1.
\end{equation*}
The same substitution gives
\begin{equation*}
 (\mathcal{C}_\infty^*a)(u)=
 \begin{cases}
  J(1),&0<u\leq1,\\ 
  J(e^{1-u}),&u>1.
 \end{cases}
\end{equation*}
Consequently,
\begin{equation*}
 \|\mathcal{C}_\infty^*a\|_{L_q(0,\infty)}  = \bigl(J(1)^q+\|h\|_{N_q}^q\bigr)^{\frac{1}{q}}.
\end{equation*}
By Lemma~\ref{nq vs llogl},
\begin{equation*}
 \|h\|_{N_q}
 \leq\|\mathcal{C}_\infty^*a\|_{L_q(0,\infty)}
 \leq c_{\mathrm{abs}}\|h\|_{N_q}.
\end{equation*}
Combining these estimates with
\eqref{eq:BW-Hardy-comparison-forward}
and \eqref{eq:BW-Hardy-comparison-backward}
proves both inequalities.
\end{proof}

\begin{lem}\label{nq holder inequality}
Let $1<q_1<\infty$, and let $q_2,q_3\in[q_1,\infty]$ satisfy
\begin{equation*}
 \frac1{q_1}=\frac1{q_2}+\frac1{q_3}.
\end{equation*}
For measurable functions $f$ and $g$ on $(0,1)$,
\begin{equation*}
 \|fg\|_{N_{q_1}}
 \leq c_{\mathrm{abs}}q_1
 \|f\|_{L_{\infty,q_2}(0,1)}
 \|\mathcal{C}\mu(g)\|_{L_{1,q_3}(0,1)}.
\end{equation*}
\end{lem}

\begin{proof}
We may assume that the norms on the right-hand side are finite. The rearrangement inequality
\begin{equation*}
 \mu(t,fg)\leq\mu\left(\frac{t}{2},f\right)\cdot \mu\left(\frac{t}{2},g\right),
 \ \ 0<t<1,
\end{equation*}
and $\ell(\frac{t}{2})\leq(1+\log2)\ell(t)$ imply
\begin{equation*}
 \|fg\|_{N_{q_1}}
 \leq2(1+\log2)\|\mu(f)\cdot\mu(g)\|_{N_{q_1}}.
\end{equation*}
It therefore suffices to prove the estimate for nonnegative nonincreasing $f$ and $g$.

For $u>1$, define
\begin{equation*}
 F(u)=f(e^{1-u}), \ \ b(u)=e^{1-u}g(e^{1-u}),
\end{equation*}
and extend both functions by zero to $(0,1)$. Then $F$ is nondecreasing on $(0,\infty)$ and
\begin{equation*}
 \left\|u\mapsto\frac{F(u)}{u}\right\|_{L_{q_2}(0,\infty)} = \|f\|_{L_{\infty,q_2}(0,1)}.
\end{equation*}
Also set
\begin{equation*}
 B(u)=\int\limits_u^\infty b(v)\,dv,
 \ \ u>0.
\end{equation*}
For $u>1$, the substitution $t=e^{1-u}$ gives
$B(u)=t(\mathcal{C}g)(t)$, while $B(u)=\|g\|_1$ for $0<u\leq1$.
If $q_3<\infty$, it follows that
\begin{equation*}
 \|B\|_{L_{q_3}(0,\infty)}^{q_3}
 =\|g\|_1^{q_3}+\|\mathcal{C}g\|_{L_{1,q_3}(0,1)}^{q_3}.
\end{equation*}
Since $g$ is nonincreasing, $(\mathcal{C}g)(t)\geq\|g\|_1$.
Hence
\begin{equation*}
 \|B\|_{L_{q_3}(0,\infty)}
 \leq(q_3+1)^{\frac{1}{q}_3}\|\mathcal{C}g\|_{L_{1,q_3}(0,1)}
 \leq2\|\mathcal{C}g\|_{L_{1,q_3}(0,1)}.
\end{equation*}
For $q_3=\infty$, we instead have
\begin{equation*}
 \|B\|_{L_\infty(0,\infty)}
 =\|g\|_1=\|\mathcal{C}g\|_{L_{1,\infty}(0,1)}.
\end{equation*}

Because $fg$ is nonincreasing, a change of variables yields
\begin{equation*}
 \|fg\|_{N_{q_1}}
 =\left\|u\mapsto \int\limits_u^\infty\frac{F(v)b(v)}{v}\,dv
   \right\|_{L_{q_1}(1,\infty)}.
\end{equation*}
Monotonicity gives
\begin{equation*}
 F(v)\leq\frac1v\int\limits_v^{2v}F(w)\,dw,
 \ \ v>0.
\end{equation*}
Thus Tonelli's theorem implies
\begin{align}
 \int\limits_u^\infty\frac{F(v)b(v)}v\,dv
 &\leq\int\limits_u^\infty\frac{b(v)}{v^2}
       \int\limits_v^{2v}F(w)\,dw\,dv\nonumber =\int\limits_u^\infty F(w)
       \int\limits_{\max(u,\frac{w}{2})}^w\frac{b(v)}{v^2}\,dv\,dw\nonumber\\
 &\leq4\int\limits_u^\infty\frac{F(w)}{w^2}
       \int\limits_{\frac{w}{2}}^\infty b(v)\,dv\,dw\nonumber =4(\mathcal{C}_{\infty}^*K)(u),
 \label{eq:BW-Nq-Hardy-majorant}
\end{align}
where
\begin{equation*}
 K(w)=\frac{F(w)}w B\left(\frac{w}{2}\right),
 \ \ w>0.
\end{equation*}
The Hardy inequality for $C_{\infty}^*$ and H\"older's inequality give
\begin{align*}
 \|fg\|_{N_{q_1}}
 &\leq4q_1\|K\|_{L_{q_1}(0,\infty)}\\
 &\leq4q_1
       \left\|w\mapsto\frac{F(w)}{w}\right\|_{L_{q_2}(0,\infty)}
       \left\|w\mapsto B\left(\frac{w}{2}\right)\right\|_{L_{q_3}(0,\infty)}\\
 &=4q_1 2^{\frac{1}{q_3}}
       \left\|w\mapsto\frac{F(w)}{w}\right\|_{L_{q_2}(0,\infty)}
       \|B\|_{L_{q_3}(0,\infty)}\\
 &\leq c_{\mathrm{abs}}q_1
       \|f\|_{L_{\infty,q_2}(0,1)}
       \|\mathcal{C}g\|_{L_{1,q_3}(0,1)}.
\end{align*}
Together with the initial rearrangement reduction, this proves the lemma.
\end{proof}

\begin{lem}[Tail Hardy estimate]\label{BW T lemma}
Let $1<q<\infty$. For every measurable function $h$ on $(0,1)$,
\begin{equation}\label{eq:BW-tail-Hardy}
 \|\mathfrak Th\|_{L_{\infty,q}(0,1)}
 \leq q'\|h\|_{L_{1,q}(0,1)}.
\end{equation}
\end{lem}

\begin{proof}
Set $H(u)=e^{1-u}\mu(e^{1-u},h)$ for $u>1$, and extend $H$ by zero to $(0,1)$. The substitution $t=e^{1-u}$ gives
\begin{align*}
 \|H\|_{L_q(0,\infty)}=\|h\|_{L_{1,q}(0,1)},
\end{align*} 
and
\begin{align*} 
 (\mathfrak Th)(e^{1-u})=\int\limits_1^u H(v)\,dv,  \ \ u>1.
\end{align*}
Since $\mathfrak Th$ is nonnegative and nonincreasing,
\begin{equation*}
 \|\mathfrak Th\|_{L_{\infty,q}(0,1)}
 =\|\mathcal{C}_{\infty}(H)\|_{L_q(0,\infty)}
 \leq q'\|H\|_{L_q(0,\infty)}
\end{equation*}
by the classical Hardy inequality. This proves
\eqref{eq:BW-tail-Hardy}.
\end{proof}

\begin{proof}[Proof of Theorem~\ref{thm:working-Brezis-Wainger}]
Set
\begin{equation*}
 z=f\star_{\theta_1,\theta_2}x,
 \ \ h(t)=\mu(t,x)\mu(t,f),
 \ \ 0<t<1.
\end{equation*}
By Theorem~\ref{thm:operator-oneil},
\begin{equation*}
 \mu(t,z) \leq t(\mathcal{C}\mu(f))(t)\cdot(\mathcal{C}\mu(x))(t) + (\mathfrak Th)(t), \ \ 0<t<1.
\end{equation*}
Since $\mathfrak Th$ is nonnegative and nonincreasing, Minkowski's
inequality gives
\begin{equation}\label{eq:BW-convolution-reduction}
 \|z\|_{L_{\infty,q_1}}
 \leq
 \left\|t\mapsto
 \frac{t(\mathcal{C}\mu(f))(t)\cdot(\mathcal{C}\mu(x))(t)}{\ell(t)}
 \right\|_{L_{q_1}((0,1),\frac{dt}{t})}
 +\|\mathfrak Th\|_{L_{\infty,q_1}(0,1)}.
\end{equation}

(i). H\"older's inequality and Lemma~\ref{lem:logarithmic-Hardy-estimate} yield
\begin{align*}
 &\left\|t\mapsto
 \frac{t(\mathcal{C}\mu(f))(t)\cdot(\mathcal{C}\mu(x))(t)}{\ell(t)}
 \right\|_{L_{q_1}((0,1),\frac{dt}{t})}\\
 &\leq
 \left\|t\mapsto t(\mathcal{C}\mu(f))(t)
 \right\|_{L_{q_3}((0,1),\frac{dt}{t})}\cdot
 \left\|t\mapsto\frac{(\mathcal{C}\mu(x))(t)}{\ell(t)}
 \right\|_{L_{q_2}((0,1),\frac{dt}{t})}\\
 &\leq
 2\|\mathcal{C}\mu(f)\|_{L_{1,q_3}(0,1)}\cdot
 \|x\|_{L_{\infty,q_2}}.
\end{align*}
Here we used the fact that $\mathcal{C}\mu(f)$ is nonincreasing. Moreover, Lemmas~\ref{lem:tail-log-equivalence} and~\ref{nq holder inequality} give
\begin{align*}
 \|\mathfrak Th\|_{L_{\infty,q_1}(0,1)}
 &\leq c_{\mathrm{abs}}q_1'\|h\|_{N_{q_1}}\\
 &\leq c_{\mathrm{abs}}q_1q_1'
       \|x\|_{L_{\infty,q_2}}
       \|\mathcal{C}\mu(f)\|_{L_{1,q_3}(0,1)}.
\end{align*}
Substituting these estimates into
\eqref{eq:BW-convolution-reduction} proves (i).

(ii). The classical weighted Hardy inequality gives
\begin{align*}
 \left\|t\mapsto t^{\frac{1}{p}}(\mathcal{C}\mu(x))(t)
 \right\|_{L_{q_2}((0,1),\frac{dt}{t})}
 &\leq p'\|x\|_{L_{p,q_2}},\\
 \left\|t\mapsto t^{\frac{1}{p}'}(\mathcal{C}\mu(f))(t)
 \right\|_{L_{q_3}((0,1),\frac{dt}{t})}
 &\leq p\|f\|_{L_{p',q_3}}.
\end{align*}
Since $\ell(t)\geq1$, H\"older's inequality therefore yields
\begin{align*}
 &\left\|t\mapsto \frac{t(\mathcal{C}\mu(f))(t)\cdot(\mathcal{C}\mu(x))(t)}{\ell(t)} \right\|_{L_{q_1}((0,1),\frac{dt}{t})}\\
 &\ \ \leq  \left\|t\mapsto t^{\frac{1}{p}}(\mathcal{C}\mu(x))(t)  \right\|_{L_{q_2}((0,1),\frac{dt}{t})}\cdot
 \left\|t\mapsto t^{\frac{1}{p}'}(\mathcal{C}\mu(f))(t)  \right\|_{L_{q_3}((0,1),\frac{dt}{t})}\\
 &\ \ \leq pp'\|x\|_{L_{p,q_2}}\|f\|_{L_{p',q_3}}.
\end{align*}
Also, $h$ is nonnegative and nonincreasing. Hence another application of H\"older's inequality gives
\begin{align*}
 \|h\|_{L_{1,q_1}(0,1)}
 &=\left\|t\mapsto t\mu(t,x)\mu(t,f)
   \right\|_{L_{q_1}((0,1),\frac{dt}{t})}\\
 &\leq
 \left\|t\mapsto t^{\frac{1}{p}}\mu(t,x)
 \right\|_{L_{q_2}((0,1),\frac{dt}{t})}
 \left\|t\mapsto t^{\frac{1}{p}'}\mu(t,f)
 \right\|_{L_{q_3}((0,1),\frac{dt}{t})}\\
 &=\|x\|_{L_{p,q_2}}\|f\|_{L_{p',q_3}}.
\end{align*}
Lemma~\ref{BW T lemma} now implies
\begin{equation*}
 \|\mathfrak Th\|_{L_{\infty,q_1}(0,1)}
 \leq q_1'\|x\|_{L_{p,q_2}}\|f\|_{L_{p',q_3}}.
\end{equation*}
Combining these estimates with
\eqref{eq:BW-convolution-reduction} proves (ii).
\end{proof}

\subsection{Besov--Lorentz estimates}

For $\alpha>0$, $1\leq\rho\leq\infty$, and either $1<p<\infty$, $1\leq q\leq\infty$, or $(p,q)=(1,1)$, use the norm
\begin{equation}\label{eq:Besov-Lorentz-norm}
\|z\|_{B^{\alpha}_{p,q;\rho}(\mathbb{T}^{d}_{\theta})} = \|z\|_{L_{p,q}(\mathbb{T}^{d}_{\theta})}
+ \left\|\left(2^{j\alpha} \|\widetilde\varphi_{j}\ast z\|_{L_{p,q}(\mathbb{T}^{d}_{\theta})} \right)_{j\geq0}\right\|_{\ell_\rho}.
\end{equation}
Equivalent Banach norms may be used for $L_{p,q}$ when $p>1$. At the endpoint we set $L_{1,1}=L_{1}$, with its usual norm. In the notation without a middle index, $B^\alpha_{p,\rho}$ means $B^\alpha_{p,p;\rho}$.

\begin{thm}[Besov--Lorentz convolution inequality]
\label{thm:Besov-Lorentz}
Let $1<p_{1},p_{2}<\infty$, $1\leq q_{1},q_{2},s\leq\infty$, and suppose that
\begin{equation*}
\frac{1}{p_{1}}+\frac{1}{p_{2}}>1.
\end{equation*}
Define $r$ by
\begin{equation*}
1+\frac{1}{r} = \frac{1}{p_{1}}+\frac{1}{p_{2}},
\end{equation*}
and suppose that
\begin{equation*}
\frac{1}{s} \leq \frac{1}{q_{1}}+\frac{1}{q_{2}}.
\end{equation*}
Let $\alpha>0$ and $1\leq\rho\leq\infty$. If $x\in L_{p_{1},q_{1}}(\mathbb{T}_{\theta_{1}}^{d})$ and  $y\in B^{\alpha}_{p_{2},q_{2};\rho} (\mathbb{T}_{\theta_{2}}^{d}),$ then
\begin{equation*}
x\star_{\theta_{1},\theta_{2}}y \in B^{\alpha}_{r,s;\rho} (\mathbb{T}_{\theta_{1}+\theta_{2}}^{d})
\end{equation*}
and
\begin{align*}
\left\|
x\star_{\theta_{1},\theta_{2}}y
\right\|_
{B^{\alpha}_{r,s;\rho}
(\mathbb{T}_{\theta_{1}+\theta_{2}}^{d})}
\leq
C
\|x\|_{L_{p_{1},q_{1}}(\mathbb{T}_{\theta_{1}}^{d})}
\|y\|_{B^{\alpha}_{p_{2},q_{2};\rho} (\mathbb{T}_{\theta_{2}}^{d})}.
\end{align*}
Similarly, if $x\in B^{\alpha}_{p_{1},q_{1};\rho} (\mathbb{T}_{\theta_{1}}^{d})$ and $y\in L_{p_{2},q_{2}}(\mathbb{T}_{\theta_{2}}^{d}),$ then
\begin{equation*}
x\star_{\theta_{1},\theta_{2}}y
\in B^{\alpha}_{r,s;\rho} (\mathbb{T}_{\theta_{1}+\theta_{2}}^{d})
\end{equation*}
and
\begin{align*}
\left\|
x\star_{\theta_{1},\theta_{2}}y
\right\|_{B^{\alpha}_{r,s;\rho}
(\mathbb{T}_{\theta_{1}+\theta_{2}}^{d})}
\leq
C \|x\|_{B^{\alpha}_{p_{1},q_{1};\rho}
(\mathbb{T}_{\theta_{1}}^{d})}
\|y\|_{L_{p_{2},q_{2}}(\mathbb{T}_{\theta_{2}}^{d})}.
\end{align*}
\end{thm}

\begin{proof}
The Fourier formula \eqref{eq:multiplier-convolution-commutation}
gives, for each $j\geq0$,
\begin{equation*}
\widetilde\varphi_{j}\ast(x\star_{\theta_1,\theta_2}y)
=x\star_{\theta_1,\theta_2}(\widetilde\varphi_{j}\ast y).
\end{equation*}
For every $j$, Theorem~\ref{thm:operator-oneil} shows that
\begin{equation*}
\bigl(x,\widetilde\varphi_{j}\ast y,
x\star_{\theta_1,\theta_2}(\widetilde\varphi_{j}\ast y)\bigr)
\end{equation*}
satisfies the O'Neil relation. Apply
Corollary~\ref{cor:ONeil-Lorentz-transfer}, multiply by
$2^{j\alpha}$, and take the $\ell_\rho$-norm in $j$. Applying the
same corollary to the triple in
\eqref{eq:convolution-ONeil-triple} controls the $L_{r,s}$ term in
\eqref{eq:Besov-Lorentz-norm}.  The second assertion follows by
transferring the Littlewood--Paley multiplier to $x$ instead.
\end{proof}

\begin{corollary}[Besov endpoint inequality]
\label{cor:Besov-endpoint}
Let $\alpha>0$, $1\leq\rho\leq\infty$, and let $1\leq q\leq s\leq\infty$. Suppose that either $1<p<\infty$, or $p=q=s=1$. If $x\in L_{1}(\mathbb{T}_{\theta_{1}}^{d})$ and $y\in B^{\alpha}_{p,q;\rho}(\mathbb{T}_{\theta_{2}}^{d}),$ then
\begin{equation*}
x\star_{\theta_{1},\theta_{2}}y \in B^{\alpha}_{p,s;\rho}(\mathbb{T}_{\theta_{1}+\theta_{2}}^{d})
\end{equation*}
and
\begin{align*}
\left\|
x\star_{\theta_{1},\theta_{2}}y
\right\|_{B^{\alpha}_{p,s;\rho}
(\mathbb{T}_{\theta_{1}+\theta_{2}}^{d})}
\leq
C_{\alpha,\rho,p,q,s}
\|x\|_{L_{1}(\mathbb{T}_{\theta_{1}}^{d})}
\|y\|_{B^{\alpha}_{p,q;\rho}
(\mathbb{T}_{\theta_{2}}^{d})}.
\end{align*}
\end{corollary}

\begin{proof}
By \eqref{eq:HLP-first}--\eqref{eq:HLP-second}, convolution by
an $L_{1}$ element acts on each $L_{p,q}$, $p>1$, with its norm
bounded by a constant depending on the chosen equivalent Lorentz
norm.  When $(p,q)=(1,1)$ this is \eqref{eq:L1}.
Apply the resulting estimate to $y$ and to every block
$\widetilde\varphi_{j}\ast y$, using
\eqref{eq:multiplier-convolution-commutation}.  The inclusion
$L_{p,q}\hookrightarrow L_{p,s}$ for $q\leq s$, together with
\eqref{eq:Besov-Lorentz-norm}, completes the proof.
\end{proof}

\begin{corollary}[Young inequality in Besov spaces]
\label{cor:Young-Besov}
Let $1\leq p_{1},p_{2},r<\infty$ satisfy
\begin{equation*}
1+\frac{1}{r} =\frac{1}{p_{1}}+\frac{1}{p_{2}}.
\end{equation*}
Let $\alpha>0$ and $1\leq\rho\leq\infty$. If $x\in L_{p_{1}}(\mathbb{T}_{\theta_{1}}^{d})$ and $y\in B^{\alpha}_{p_{2},\rho}(\mathbb{T}_{\theta_{2}}^{d}),$ then
\begin{equation*}
x\star_{\theta_{1},\theta_{2}}y \in B^{\alpha}_{r,\rho}
(\mathbb{T}_{\theta_{1}+\theta_{2}}^{d})
\end{equation*}
and
\begin{align*}
\left\|
x\star_{\theta_{1},\theta_{2}}y
\right\|_
{B^{\alpha}_{r,\rho}
(\mathbb{T}_{\theta_{1}+\theta_{2}}^{d})}
\leq
C
\|x\|_{L_{p_{1}}(\mathbb{T}_{\theta_{1}}^{d})}
\|y\|_
{B^{\alpha}_{p_{2},\rho}
(\mathbb{T}_{\theta_{2}}^{d})}.
\end{align*}
The analogous assertion holds when the Besov regularity is imposed on $x$ instead of $y$.
\end{corollary}

\begin{proof}
Use the sharp Young estimate \eqref{eq:Young} for the
low-frequency $L_r$ term and for each identity
$\widetilde\varphi_{j}\ast(x\star y)
=x\star(\widetilde\varphi_{j}\ast y)$.
Multiply the blockwise inequalities by $2^{j\alpha}$ and take the
$\ell_\rho$-norm.  Moving the multiplier to the first input gives
the symmetric assertion.
\end{proof}

\section{Sobolev applications}\label{sec:sobolev-applications}

We now apply the convolution estimates to fractional Sobolev embeddings on quantum tori. The classical approach uses Bessel and Riesz potentials to obtain Sobolev inequalities; see~\cite{AS,Sobolev}. Related embeddings for quantum tori were studied in~\cite{Spera, XiongXuYin}. Here we combine potential-kernel estimates with the sharp O'Neil inequality to obtain embeddings in Lorentz and symmetric operator spaces. At the critical order, we reduce the problem to the boundedness of a scalar Hardy-type operator, following ideas from the classical theory \cite{Klimov1969,Klimov1970,EKP}. 

\subsection{Bessel and Riesz potentials}

We now derive Lorentz refinements of Sobolev embeddings on $\mathbb{T}_{\theta}^{d}$ from the preceding convolution inequalities. The main ingredients are the convolution representations of Bessel and Riesz potentials and the Lorentz estimates for their periodic kernels.

For $\sigma>0$, let $G_{\sigma}$ denote the periodic \textit{Bessel kernel} determined by
\begin{equation}\label{eq:periodic-Bessel-kernel}
\widehat{G_{\sigma}}(n)
= (1+|n|^{2})^{-\frac{\sigma}{2}},
\ \  n\in\mathbb{Z}^{d}.
\end{equation}
It follows from the Fourier multiplier--convolution correspondence that
\begin{equation}\label{eq:Bessel-potential-convolution}
J_{\theta}^{-\sigma}x = G_{\sigma}\ast x.
\end{equation}

We also define the periodic \textit{Riesz kernel} $R_{\sigma}$ by
\begin{equation}\label{eq:periodic-Riesz-kernel}
\widehat{R_{\sigma}}(0)=0,
\ \ 
\widehat{R_{\sigma}}(n)=|n|^{-\sigma}, \ \  n\neq0.
\end{equation}
Therefore, for every $x$ with $\widehat{x}(0)=\tau_{\theta}(x)=0$,
\begin{equation}\label{eq:Riesz-potential-convolution}
(-\Delta_{\theta})^{-\frac{\sigma}{2}}x
= R_{\sigma}\ast x.
\end{equation}

\begin{lem}[Periodic potential kernels]
\label{lem:periodic-potential-kernels}
Let $0<\sigma<d$. Then
\begin{equation}\label{eq:potential-kernel-weak-Lp}
G_{\sigma},R_{\sigma}
\in L_{\frac{d}{d-\sigma},\infty}(\mathbb{T}^{d}),
\end{equation}
and
\begin{equation*}
\|G_{\sigma}\|_{L_{\frac{d}{d-\sigma},\infty}(\mathbb{T}^{d})} + \|R_{\sigma}\|_{L_{\frac{d}{d-\sigma},\infty}(\mathbb{T}^{d})}
\leq C_{d,\sigma}.
\end{equation*}
More generally, let $\beta\in\mathbb{N}_{0}^{d}$ and suppose that $0<\sigma-|\beta|<d.$ The Fourier multiplier with symbol $(in)^{\beta}(1+|n|^{2})^{-\frac{\sigma}{2}}$ is convolution with a kernel $G_{\sigma,\beta}$ satisfying
\begin{equation}\label{eq:differentiated-potential-kernel}
G_{\sigma,\beta}
\in
L_{\frac{d}{d-\sigma+|\beta|},\infty}(\mathbb{T}^{d}).
\end{equation}
\end{lem}

\begin{proof}
The Bessel kernel admits the heat-kernel representation
\begin{equation*}
G_{\sigma}(s)
= \frac{1}{\Gamma(\frac{\sigma}{2})}
\int\limits_{0}^{\infty} t^{\frac{\sigma}{2}-1}e^{-t}H_{t}(s)\,dt.
\end{equation*}
Similarly,
\begin{equation*}
R_{\sigma}(s) = \frac{1}{\Gamma(\frac{\sigma}{2})} \int\limits_{0}^{\infty} t^{\frac{\sigma}{2}-1}
\bigl(H_{t}(s)-1\bigr)\,dt.
\end{equation*}
The standard periodic heat-kernel estimates imply
\begin{equation*}
|G_{\sigma}(s)|+|R_{\sigma}(s)| \leq C_{d,\sigma} \bigl(1+|s|^{\sigma-d}\bigr),
\end{equation*}
where $|s|$ denotes the distance from $s$ to the origin of $\mathbb{T}^{d}$. Consequently,
\begin{equation*}
\mu(t,G_{\sigma})+\mu(t,R_{\sigma}) \leq C_{d,\sigma}t^{-\frac{d-\sigma}{d}}, \ \  0<t<1,
\end{equation*}
which proves \eqref{eq:potential-kernel-weak-Lp}.

Differentiating the heat-kernel representation gives
\begin{equation*}
|G_{\sigma,\beta}(s)| \leq C_{d,\sigma,\beta} \bigl(1+|s|^{\sigma-|\beta|-d}\bigr).
\end{equation*}
This yields \eqref{eq:differentiated-potential-kernel}. The corresponding classical kernel estimates may also be obtained by periodizing the Bessel and Riesz kernels on $\mathbb{R}^{d}$; see \cite{G2014,SS}.
\end{proof}

\subsection{Subcritical Sobolev embeddings}
\begin{thm}[Lorentz--Sobolev embedding]
\label{thm:Lorentz-Sobolev-embedding}
Let $s_{0},s_{1}\in\mathbb{R}$ and put $\sigma=s_{0}-s_{1}.$ Suppose that $0<\sigma<d,  1<p<\frac{d}{\sigma}, 1\leq q\leq\infty,$ and define $r$ by
\begin{equation}\label{eq:Lorentz-Sobolev-exponents}
\frac{1}{r} = \frac{1}{p}-\frac{\sigma}{d}.
\end{equation}
Then
\begin{equation}\label{eq:Lorentz-Sobolev-embedding}
H_{p,q}^{s_{0}}(\mathbb{T}_{\theta}^{d}) \hookrightarrow H_{r,q}^{s_{1}}(\mathbb{T}_{\theta}^{d}).
\end{equation}
More precisely,
\begin{equation}\label{eq:Lorentz-Sobolev-estimate}
\|x\|_{H_{r,q}^{s_{1}}(\mathbb{T}_{\theta}^{d})} \leq C \|x\|_{H_{p,q}^{s_{0}}(\mathbb{T}_{\theta}^{d})},
\end{equation}
where the constant is independent of $\theta$.
\end{thm}

\begin{proof}
By the lifting property of the Bessel potentials,
\begin{equation*}
J_{\theta}^{s_{1}}x = J_{\theta}^{-\sigma}J_{\theta}^{s_{0}}x.
\end{equation*}
Using \eqref{eq:Bessel-potential-convolution}, we obtain
\begin{equation*}
J_{\theta}^{s_{1}}x = G_{\sigma}\ast J_{\theta}^{s_{0}}x.
\end{equation*}
Set $a=\frac{d}{d-\sigma}.$ By Lemma~\ref{lem:periodic-potential-kernels}, $G_{\sigma}\in L_{a,\infty}(\mathbb{T}^{d}).$ Moreover,
\begin{equation*}
1+\frac{1}{r} = \frac{1}{a}+\frac{1}{p}.
\end{equation*}
Corollary~\ref{cor:Young-ONeil}, applied with Lorentz indices $\infty$ and $q$, therefore gives
\begin{equation*}
\|G_{\sigma}\ast J_{\theta}^{s_{0}}x\|_{L_{r,q}(\mathbb{T}_{\theta}^{d})}
\leq C \|G_{\sigma}\|_{L_{a,\infty}(\mathbb{T}^{d})} \|J_{\theta}^{s_{0}}x\|_{L_{p,q}(\mathbb{T}_{\theta}^{d})}.
\end{equation*}
This is precisely \eqref{eq:Lorentz-Sobolev-estimate}.
\end{proof}

\begin{rem}
Taking $q=p$ in Theorem~\ref{thm:Lorentz-Sobolev-embedding} gives $H_p^{s_0}\hookrightarrow H_{r,p}^{s_1}$ and hence the known $H_p^{s_0}\hookrightarrow H_r^{s_1}$ embedding of \cite[Theorem~6.6(i)]{XiongXuYin}, since $p<r$. Likewise, the mean-zero Riesz-potential version follows from \eqref{eq:Riesz-potential-convolution} and Lemma~\ref{lem:periodic-potential-kernels}. The arbitrary secondary Lorentz index in the theorem gives the corresponding refined target.
\end{rem}

\begin{cor}[Sobolev estimates for derivatives]
\label{cor:Lorentz-Sobolev-derivatives}
Let $s\in\mathbb{R}$ and $\beta\in\mathbb{N}_{0}^{d}$. Suppose that $ 0<s-|\beta|<d, 1<p<\frac{d}{s-|\beta|},  1\leq q\leq\infty,$ and define $r$ by
\begin{equation*}
\frac{1}{r} = \frac{1}{p}-\frac{s-|\beta|}{d}.
\end{equation*}
Then
\begin{equation}\label{eq:Lorentz-Sobolev-derivative-estimate}
\|(\partial^{\theta})^{\beta}x\|_{L_{r,q}(\mathbb{T}_{\theta}^{d})} \leq C \|J_{\theta}^{s}x\|_{L_{p,q}(\mathbb{T}_{\theta}^{d})}.
\end{equation}
\end{cor}

\begin{proof}
Write
\begin{equation*}
(\partial^{\theta})^{\beta}x = (\partial^{\theta})^{\beta} J_{\theta}^{-s}J_{\theta}^{s}x.
\end{equation*}
The operator $(\partial^{\theta})^{\beta}J_{\theta}^{-s}$ is convolution with the kernel $G_{s,\beta}$. By Lemma~\ref{lem:periodic-potential-kernels},
\begin{equation*}
G_{s,\beta} \in L_{\frac{d}{d-s+|\beta|},\infty}(\mathbb{T}^{d}).
\end{equation*}
The assertion now follows from Corollary~\ref{cor:Young-ONeil}.
\end{proof}

\begin{rem}
The strict inequality $s>|\beta|$ in Corollary~\ref{cor:Lorentz-Sobolev-derivatives} is essential for the potential-kernel argument. If $s=|\beta|$, then $(\partial^{\theta})^{\beta}J_{\theta}^{-s}$ is an order-zero Fourier multiplier of Riesz-transform type. Its boundedness is no longer a consequence of Young's inequality and instead requires a Fourier multiplier theorem. On general symmetric spaces, such an argument typically requires nontrivial Boyd indices.
\end{rem}

The weak target at $p=1$ follows directly from \eqref{eq:HLP-second}: for $a=\frac{d}{d-\sigma}>1$, the weak-$L_a$ space is fully symmetric under its equivalent Marcinkiewicz norm and $R_\sigma\in L_{a,\infty}$; hence $\|R_\sigma\ast x\|_{L_{a,\infty}}\leq C_a\|R_\sigma\|_{L_{a,\infty}}\|x\|_1$ for every mean-zero $x\in L_{1}$. This is consistent with \cite[Theorem~3.8 and Corollary~6.3]{XiongXuYin}.

\subsection{Critical logarithmic Sobolev estimate}
\begin{thm}[Critical logarithmic Sobolev estimate]
\label{thm:critical-logarithmic-Sobolev}
Let $0<\sigma<d$ and put $p=\frac{d}{\sigma}.$ Let $1<q<\infty$. Then, for every $x\in H_{p,q}^{\sigma}(\mathbb{T}_{\theta}^{d})$,
\begin{equation}\label{eq:critical-logarithmic-Sobolev}
\left\|
\frac{\mu(t,x)}{1+|\log t|}
\right\|_{L_q\left(0,1;\frac{dt}{t}\right)}
\leq C \|J_{\theta}^{\sigma}x\|_{L_{p,q}(\mathbb{T}_{\theta}^{d})}.
\end{equation}
\end{thm}

\begin{proof}
Set $y=J_{\theta}^{\sigma}x.$ Then $x=G_{\sigma}\ast y.$ Since $p'=\frac{d}{d-\sigma},$ Lemma~\ref{lem:periodic-potential-kernels} gives
\begin{equation*}
G_{\sigma}\in L_{p',\infty}(\mathbb{T}^{d}).
\end{equation*}
Apply Theorem~\ref{thm:working-Brezis-Wainger} with
\begin{equation*}
q_{1}=q_{2}=q, \ \ q_{3}=\infty.
\end{equation*}
This gives \eqref{eq:critical-logarithmic-Sobolev}.
\end{proof}

\begin{remark}[Relation to the embeddings of Xiong, Xu and Yin]
\label{rem:XXY-consequences}
The convolution estimates recover several Sobolev and Besov embeddings established by Xiong, Xu and Yin \cite{XiongXuYin}. Indeed, representing Bessel potentials as semi-convolutions with scalar kernels and applying the O'Neil inequality yields the potential-space embeddings of \cite[Theorem~6.6(i)]{XiongXuYin}.

Similarly, the Young inequality from \cite{ConvolutionI}, applied to Littlewood--Paley blocks, gives the Besov embeddings of \cite[Theorem~6.2]{XiongXuYin}. Real interpolation then yields the Besov--Lorentz embeddings in \cite[Corollary~6.3]{XiongXuYin}, while summation of the dyadic estimates gives the limiting embedding into $L_\infty$ in the same corollary. These arguments combine the convolution inequalities with scalar kernel estimates, and the resulting constants are independent of the deformation parameter.

At the critical relation $p\sigma=d$, Theorem~\ref{thm:critical-logarithmic-Sobolev} provides a single limiting rearrangement estimate that implies the corresponding embeddings into every finite Lebesgue space. This is consistent with the critical Sobolev embedding discussion in \cite[Remark~6.7]{XiongXuYin}.
\end{remark}

\begin{rem}
The preceding results use only scalar periodic kernels and the convolution inequalities on $\mathbb{T}_{\theta}^{d}$. Consequently, their constants are independent of the deformation parameter $\theta$. The noncommutative structure enters through the convolution action and the generalized singular value estimates, while the local singularities of the Bessel and Riesz kernels are the same as in the classical case.
\end{rem}

\subsection{A critical fractional symmetric-space estimate}
\label{subsec:fractional-symmetric}

In this subsection we derive fractional Sobolev estimates using the Hardy operator $\mathcal{T}$ from \eqref{operatorT} and the sharp convolution inequality of Theorem~\ref{thm:operator-oneil}.

\begin{lem}[Monotonicity of the Hardy operator]
\label{lem:fractional-Hardy-bounds}
For $h\geq0$ in $L_{1}(0,1)$, $\mathcal{T}h$ is decreasing and
\begin{equation}\label{eq:fractional-Hardy-primitive}
\sqrt t\int\limits_0^t h(s)\,ds+t\int\limits_t^1\frac{h(s)}{\sqrt s}\,ds
\leq\int\limits_0^t(\mathcal{T}h)(u)\,du,\ \  0<t<1.
\end{equation}
\end{lem}
\begin{proof}
Local absolute continuity gives
\begin{equation*}
    (\mathcal{T}h)'(t)=-\frac1{2t^{\frac{3}{2}}}\int\limits_0^t h(s)ds\leq0
\end{equation*}
The derivative of the left-hand side of \eqref{eq:fractional-Hardy-primitive} is
\begin{equation*}
    \frac1{2\sqrt t}\int\limits_0^t h(s)ds+\int\limits_t^1 \frac{h(s)}{\sqrt{s}}  s\,ds
\end{equation*} and is at most $Th(t)$. Both primitives vanish at zero; for the tail term use 
\begin{equation*}
    t\int\limits_t^1\frac{h(s)}{\sqrt{s}} \,ds\leq\sqrt t\|h\|_1.
\end{equation*}
Integration proves the assertion.
\end{proof}

\begin{thm}[Weak-$L_2$ convolution and the Hardy operator]
\label{thm:fractional-hardiness}
Let $d\geq1$, $f\in L_{2,\infty}(\mathbb{T}_{\theta_1}^{d})$ and $x\in L_{1}(\mathbb{T}_{\theta_2}^{d})$. Then
\begin{equation}\label{eq:fractional-hardiness}
\mu(f\star_{\theta_1,\theta_2}x) \prec\!\prec 2\|f\|_{L_{2,\infty}} \mathcal{T}\mu(x).
\end{equation}
The weak-$L_2$ quasi-norm used here is $\sup_{0<t<1}t^{\frac{1}{2}}\mu(t,f)$.
\end{thm}
\begin{proof}
Put $A=\|f\|_{L_{2,\infty}(\mathbb{T}_{\theta_1}^{d})}$. Since $L_{2,\infty}\subset L_{1}$ for a normalized finite trace, Theorem~\ref{thm:operator-oneil} gives
\begin{align*}
\int\limits_0^t\mu(s,f\star_{\theta_1,\theta_2}x)\,ds
\leq \int\limits_0^t\mu(s,f)\,ds\cdot \int\limits_0^t\mu(s,x)\,ds + t\int\limits_t^1\mu(s,f)\mu(s,x)\,ds.
\end{align*}
We estimate the two terms separately. For the first input,
\begin{align*}
\int\limits_0^t\mu(s,f)\,ds
 =\int\limits_0^t\bigl(s^{\frac{1}{2}}\mu(s,f)\bigr)s^{-\frac{1}{2}}\,ds \leq A\int\limits_0^t s^{-\frac{1}{2}}\,ds=2A\sqrt t.
\end{align*}
For the tail term,
\begin{align*}
t\int\limits_t^1\mu(s,f)\mu(s,x)\,ds = t\int\limits_t^1\bigl(s^{\frac{1}{2}}\mu(s,f)\bigr)                    \mu(s,x)s^{-\frac{1}{2}}\,ds\leq At\int\limits_t^1\mu(s,x)s^{-\frac{1}{2}}\,ds.
\end{align*}
Consequently, for every $0<t<1$,
\begin{equation}\label{SubMajor}
\int\limits_0^t\mu(s,f\star_{\theta_1,\theta_2}x)\,ds\leq2A F(t),
\end{equation}
where
\begin{equation*}
F(t)=\sqrt t\int\limits_0^t\mu(s,x)\,ds
    +t\int\limits_t^1\mu(s,x)s^{-\frac{1}{2}}\,ds.
\end{equation*}
The function $F$ is locally absolutely continuous on $(0,1)$, and for almost every $t$,
\begin{align*}
F'(t)
 &=\frac1{2\sqrt t}\int\limits_0^t\mu(s,x)\,ds
   +\sqrt t\,\mu(t,x)+\int\limits_t^1\mu(s,x)s^{-\frac{1}{2}}\,ds
   -\sqrt t\,\mu(t,x)\\
 &=\frac1{2\sqrt t}\int\limits_0^t\mu(s,x)\,ds
   +\int\limits_t^1\mu(s,x)s^{-\frac{1}{2}}\,ds\\
 &\leq\frac1{\sqrt t}\int\limits_0^t\mu(s,x)\,ds
   +\int\limits_t^1\mu(s,x)s^{-\frac{1}{2}}\,ds
  =(\mathcal{T}\mu(x))(t).
\end{align*}
To integrate from zero, note that
\begin{equation*}
0\leq F(t)\leq2\sqrt t\,\|x\|_1\to 0
 \ \ (t\downarrow0).
\end{equation*}
Proposition~\ref{ThmTBdd}(i) also gives $\mathcal{T}\mu(x)\in L_{1}$. Integrating the derivative bound on $(\varepsilon,t)$ and letting
$\varepsilon\downarrow0$ therefore yields
\begin{equation*}
F(t)\leq\int\limits_0^t(\mathcal{T}\mu(x))(s)\,ds.
\end{equation*}
Combining this with \eqref{SubMajor}, we obtain
\begin{equation*}
\int\limits_0^t\mu(s,f\star_{\theta_1,\theta_2}x)\,ds
\leq2A\int\limits_0^t(\mathcal{T}\mu(x))(s)\,ds.
\end{equation*}
The function $\mathcal{T}\mu(x)$ is nonnegative and decreasing by Lemma~\ref{lem:fractional-Hardy-bounds}. The last inequality is therefore precisely the asserted submajorization.
\end{proof}

We next give an explicit periodization argument for the critical Bessel kernel. The periods below agree with the Fourier normalization fixed in the preliminaries.

\begin{lem}[The critical periodic Bessel kernel]\label{weak_L2-function}
Let $d\geq1$. The scalar distribution $f$ with Fourier coefficients
\begin{equation*}
\widehat f(n)=(1+|n|^2)^{-\frac{d}{4}},\ \  n\in\mathbb{Z}^{d},
\end{equation*}
is represented by a function in $L_{2,\infty}(\mathbb{T}^{d})$.
It is the kernel $G_{\frac{d}{2}}$ from \eqref{eq:periodic-Bessel-kernel}.
\end{lem}
\begin{proof}
We identify $\mathbb{T}^{d}$ with the fundamental cube $Q=[-\pi,\pi)^{d}$, equipped with $dm(t)=(2\pi)^{-d}\,dt$. Use the Euclidean Fourier transform
\begin{equation*}
\widehat b(\xi)=\int\limits_{\mathbb{R}^{d}}b(t)e^{-it\cdot\xi}\,dt.
\end{equation*}
Let $b$ be the Euclidean Bessel kernel with transform $(1+|\xi|^2)^{-\frac{d}{4}}$. Its heat-integral representation is
\begin{equation}\label{eq:critical-euclidean-Bessel}
b(t)=\frac1{\Gamma(\frac{d}{4})}\int\limits_0^\infty
 a^{\frac{d}{4}-1}e^{-a}(4\pi a)^{-\frac{d}{2}}
 e^{-\frac{|t|^{2}}{4a}}\,da.
\end{equation}
The integrand is nonnegative. Integrating first in $t$ shows $\int\limits_{\mathbb{R}^{d}}b(t)\,dt=1$, so $b\in L_{1}(\mathbb{R}^{d})$.
Taking the Fourier transform under the integral gives
\begin{align*}
\widehat b(\xi)
 =\frac1{\Gamma(\frac{d}{4})}\int\limits_0^\infty
       a^{\frac{d}{4}-1}e^{-a}e^{-a|\xi|^2}\,da =(1+|\xi|^2)^{-\frac{d}{4}}.
\end{align*}
The standard Macdonald-function formula for this kernel is
\begin{equation}\label{eq:critical-Macdonald}
b(t)=c_d|t|^{-\frac{d}{4}}K_{\frac{d}{4}}(|t|),\ \  t\neq0,
\end{equation}
where $c_d>0$. This follows by the integral representation of
$K_\nu$ in \eqref{eq:critical-euclidean-Bessel}; see
\cite{AS} and \cite[Example~2.2.7]{LMSZ}.

Periodize $b$ by setting
\begin{equation}\label{eq:periodization}
f(t)=(2\pi)^{d}\sum_{n\in\mathbb{Z}^{d}}b(t+2\pi n),\ \  t\in Q.
\end{equation}
The factor $(2\pi)^{d}$ accounts for normalized Haar measure.
Tonelli's theorem gives
\begin{equation*}
\int\limits_Q f(t)\,dm(t)
=\sum_{n\in\mathbb{Z}^{d}}\int\limits_Q b(t+2\pi n)\,dt
=\int\limits_{\mathbb{R}^{d}}b(t)\,dt=1.
\end{equation*}
Thus the series converges almost everywhere and in $L_{1}(Q,dm)$.
For $k\in\mathbb{Z}^{d}$, absolute integrability justifies
\begin{align*}
\widehat f(k) = \sum_{n\in\mathbb{Z}^{d}}\int\limits_Q b(t+2\pi n)e^{-ik\cdot t}\,dt =\int\limits_{\mathbb{R}^{d}}b(s)e^{-ik\cdot s}\,ds = (1+|k|^2)^{-\frac{d}{4}}.
\end{align*}
Consequently the periodized function is the required distribution.
This is the usual periodization construction; see
\cite[Chapter~3]{G2014} and \cite[Section~2.3]{SS}.

Using \eqref{eq:critical-Macdonald}, split the central term from the
remaining translates:
\begin{equation}\label{eq:f-explicit-short}
\begin{split}
f(t)&=(2\pi)^{d}c_d|t|^{-\frac{d}{4}}K_{\frac{d}{4}}(|t|)\\
&  +(2\pi)^{d}c_d\sum_{n\neq0}
|t+2\pi n|^{-\frac{d}{4}}K_{\frac{d}{4}}(|t+2\pi n|)\\
&=:g_0(t)+g_1(t).
\end{split}
\end{equation}
We first estimate $g_1$. For $n\neq0$, choose $j$ with $|n_{j}|\geq1$.
Since $|t_{j}|\leq\pi$ on $Q$,
\begin{equation*}
|t+2\pi n|\geq|t_{j}+2\pi n_{j}|\geq\pi.
\end{equation*}
The large-argument estimate for the Macdonald function gives
$K_{\frac{d}{4}}(r)\leq C_de^{-r}$ for $r\geq\pi$; see \cite{AS}.
Moreover, $|t|\leq\pi\sqrt d$, and therefore
\begin{align*}
|t+2\pi n|^{-\frac{d}{4}}K_{\frac{d}{4}}(|t+2\pi n|)
&\leq C_d\pi^{-\frac{d}{4}}e^{-|t+2\pi n|}\\
&\leq C_d\pi^{-\frac{d}{4}}e^{\pi\sqrt d}e^{-2\pi|n|}.
\end{align*}
The series $\sum_{n\in\mathbb{Z}^{d}}e^{-2\pi|n|}$ converges, so the
series defining $g_1$ converges uniformly and $g_1\in L_\infty(Q)$.

It remains to examine $g_0$ at the origin. For every $\nu>0$,
\begin{equation*}
K_\nu(r)\sim2^{\nu-1}\Gamma(\nu)r^{-\nu}
 \ \ (r\downarrow0).
\end{equation*}
Taking $\nu=\frac{d}{4}$ in \eqref{eq:f-explicit-short} gives
\begin{equation*}
g_0(t)\sim a_d|t|^{-\frac{d}{2}} \ \ (|t|\downarrow0),
\ \  a_d=(2\pi)^{d}c_d2^{\frac{d}{4}-1}\Gamma\left(\frac{d}{4}\right)>0.
\end{equation*}
The function $r\mapsto r^{\frac{d}{4}}K_{\frac{d}{4}}(r)$ is bounded on $(0,\pi\sqrt d]$, by this asymptotic at zero and continuity away from zero. Hence there is $A_d>0$ such that
\begin{equation*}
0\leq g_0(t)\leq A_d|t|^{-\frac{d}{2}},\ \  t\in Q\setminus\left\{0\right\}.
\end{equation*}
Let $v_d=\pi^{\frac{d}{2}}/\Gamma(\frac{d}{2}+1)$ be the volume of the Euclidean unit ball. For every $\lambda>0$, the preceding bound implies
\begin{align*}
m\left\{t\in Q:g_0(t)>\lambda\right\}
&\leq m\left\{t\in Q:A_d|t|^{-\frac{d}{2}}>\lambda\right\}\\
&\leq(2\pi)^{-d}
 \bigl|B(0,(A_d/\lambda)^{\frac{2}{d}})\bigr|\\
&=(2\pi)^{-d}v_d A_d^2\lambda^{-2}.
\end{align*}
These are inequalities: the pointwise singular model is an upper
bound, and its ball may extend beyond $Q$. By \eqref{LorentzEqvNorm},
\begin{align*}
\|g_0\|_{L_{2,\infty}(Q,dm)} = \sup_{\lambda>0}\lambda m\left\{|g_0|>\lambda\right\}^{\frac{1}{2}} \leq(2\pi)^{-\frac{d}{2}}v_d^{\frac{1}{2}}A_d<\infty.
\end{align*}
Finally, write $M=\|g_1\|_\infty$. If $\lambda>2M$, then
$\left\{f>\lambda\right\}\subset\left\{g_0>\lambda/2\right\}$, so
\begin{equation*}
\lambda m\left\{f>\lambda\right\}^{\frac{1}{2}}
\leq2\|g_0\|_{2,\infty}.
\end{equation*}
If $0<\lambda\leq2M$, the same quantity is at most $2M$, since
$m(Q)=1$. Taking the supremum proves $f\in L_{2,\infty}(\mathbb{T}^{d})$.
\end{proof}

For the rest of this section, put
$C_d=2\|G_{\frac{d}{2}}\|_{L_{2,\infty}(\mathbb{T}^{d})}$.

\begin{thm}[Distributional fractional Sobolev estimate]
\label{thm:critical-fractional-distribution}
For $x\in L_{1}(\mathbb{T}_{\theta}^{d})$,
\begin{equation}\label{eq:critical-fractional-distribution}
\mu(J_\theta^{-\frac{d}{2}}x)\prec\!\prec
C_d\mathcal{T}\mu(x).
\end{equation}
The constant is independent of $\theta$.
\end{thm}
\begin{proof}
Let $f=G_{\frac{d}{2}}$ be the scalar kernel of Lemma~\ref{weak_L2-function}. Its Fourier coefficients satisfy
\begin{equation*}
\widehat f(n)=(1+|n|^2)^{-\frac{d}{4}},\ \  n\in\mathbb{Z}^{d}.
\end{equation*}
For $x\in L_{1}(\mathbb{T}_{\theta}^{d})$, the semi-convolution $f\ast x$
is an $L_{1}$ operator by Lemma~\ref{FirstTechLemma}, and
\begin{align*}
\widehat{f\ast x}(n)
=\widehat f(n)\widehat x(n) =(1+|n|^2)^{-\frac{d}{4}}\widehat x(n) = \widehat{J_\theta^{-\frac{d}{2}}x}(n).
\end{align*}
Uniqueness of Fourier coefficients in $\mathcal S'$ gives the representation
\begin{equation}\label{eq:critical-potential-convolution}
J_\theta^{-\frac{d}{2}}x=(1-\Delta_\theta)^{-\frac{d}{4}}x=f\ast x.
\end{equation}
Lemma~\ref{weak_L2-function} gives $f\in L_{2,\infty}(\mathbb{T}^{d})$. Applying Theorem~\ref{thm:fractional-hardiness} with $\theta_1=0$, $\theta_2=\theta$, we obtain
\begin{align*}
\mu(J_\theta^{-\frac{d}{2}}x) =\mu(f\ast x) \prec\!\prec 2\|f\|_{L_{2,\infty}(\mathbb{T}^{d})}\mathcal{T}\mu(x) = C_d \mathcal{T}\mu(x).
\end{align*}
The kernel is scalar and depends only on $d$, so the constant is independent of $\theta$.
\end{proof}

\begin{thm}[Fractional embedding in symmetric spaces]
\label{thm:fractional-symmetric-embedding}
Let $E,F$ be symmetric Banach function spaces on $(0,1)$,
with $F$ fully symmetric, and suppose
$T:E\to F$ is bounded. Then
\begin{equation}\label{eq:fractional-symmetric-embedding}
\|J_\theta^{-\frac{d}{2}}x\|_{\mathcal F}
\leq C_d\|\mathcal{T}\|_{E\to F}\|x\|_{\mathcal{E}},
\ \ 
H_{\mathcal{E}}^{\frac{d}{2}}(\mathbb{T}_{\theta}^{d})
\hookrightarrow\mathcal F(\mathbb{T}_{\theta}^{d}).
\end{equation}
\end{thm}
\begin{proof}
Let $x\in\mathcal{E}(\mathbb{T}_{\theta}^{d})$. Since $E\hookrightarrow L_{1}$ on $(0,1)$, $x$ belongs to $L_{1}(\mathbb{T}_{\theta}^{d})$ and
Theorem~\ref{thm:critical-fractional-distribution} applies:
\begin{equation*}
\mu(J_\theta^{-\frac{d}{2}}x)\prec\!\prec C_d\mathcal{T}\mu(x).
\end{equation*}
By assumption $\mathcal{T}\mu(x)\in F(0,1)$. Full symmetry of $F$ and
\eqref{eq:operator-space-norm} therefore imply
\begin{align*}
\|J_\theta^{-\frac{d}{2}}x\|_{\mathcal F} = \|\mu(J_\theta^{-\frac{d}{2}}x)\|_{F(0,1)} \leq  C_d \|T\mu(x)\|_{F(0,1)}.
\end{align*}
Using boundedness of $T:E\to F$, we obtain
\begin{equation}\label{xinTh4.10}
\begin{split}
\|J_\theta^{-\frac{d}{2}}x\|_{\mathcal F} \leq C_d\|\mathcal{T}\|_{E\to F}\|\mu(x)\|_{E(0,1)} = C_d\|\mathcal{T}\|_{E\to F}\|x\|_{\mathcal{E}}.
\end{split}
\end{equation}
This proves the operator estimate with its stated domain.

Now let $y\in H_{\mathcal{E}}^{\frac{d}{2}}(\mathbb{T}_{\theta}^{d})$. By definition $x=J_\theta^{\frac{d}{2}}y$ belongs to $\mathcal{E}$, and Lemma~\ref{lem:J^{s}-iso} gives the distributional identity $y=J_\theta^{-\frac{d}{2}}J_\theta^{\frac{d}{2}}y$. In view of \eqref{xinTh4.10}, this distribution is represented by an element of $\mathcal{F}$.  Moreover,
\begin{align*}
\|y\|_{\mathcal F} =\|J_\theta^{-\frac{d}{2}}J_\theta^{\frac{d}{2}}y\|_{\mathcal F} \leq C_d\|\mathcal{T}\|_{E\to F}\|J_\theta^{\frac{d}{2}}y\|_{\mathcal{E}} = C_d\|\mathcal{T}\|_{E\to F}\|y\|_{H_{\mathcal{E}}^{\frac{d}{2}}}.
\end{align*}
Thus the inclusion is continuous:
\begin{equation*}
H_{\mathcal{E}}^{\frac{d}{2}}(\mathbb{T}_{\theta}^{d})
\hookrightarrow\mathcal F(\mathbb{T}_{\theta}^{d}).
\end{equation*}
\end{proof}

We now use the concave gauge $\psi$ from \eqref{PhiPsiFunc}
and the $2$-convexified Lorentz space of the preliminaries.

\begin{corollary}[Fractional Sobolev consequences]
\label{ThmSEmb}\label{cor:critical-L2-Lambda}
Let $d\geq1$, $1\leq p\leq\infty$, and let $\psi$ be defined by
\eqref{PhiPsiFunc}. Then:
\begin{itemize}
\item[(i)] For $x\in L_p(\mathbb{T}_{\theta}^{d})$,
\begin{equation*}
\|J_\theta^{-\frac{d}{2}}x\|_{L_p(\mathbb{T}_{\theta}^{d})}
\leq 2C_d\|x\|_{L_p(\mathbb{T}_{\theta}^{d})}.
\end{equation*}
\item[(ii)] For $x\in L_2(\mathbb{T}_{\theta}^{d})$,
\begin{equation}\label{eq:critical-L2-Lambda}
\|J_\theta^{-\frac{d}{2}}x\|_{\Lambda_\psi^{(2)}(\mathbb{T}_{\theta}^{d})}
\leq C_d'\|x\|_2.
\end{equation}
Consequently $H_2^{\frac{d}{2}}(\mathbb{T}_{\theta}^{d})
\hookrightarrow\Lambda_\psi^{(2)}(\mathbb{T}_{\theta}^{d})$.
\end{itemize}
Both constants are independent of $\theta$.
\end{corollary}
\begin{proof}
(i) Proposition~\ref{ThmTBdd}(i) gives
$\|\mathcal{T}\|_{L_p(0,1)\to L_p(0,1)}\leq2$ for every $1\leq p\leq\infty$.
The Lebesgue spaces in this range are fully symmetric; see \cite[Chapter~6, Theorem~6.2.1]{DoddsDePagterSukochev}. Taking $E(0,1)=F(0,1)=L_p(0,1)$ in Theorem~\ref{thm:fractional-symmetric-embedding} gives
\begin{align*}
\|J_\theta^{-\frac{d}{2}}x\|_{p} \leq C_d\|\mathcal{T}\mu(x)\|_{L_p(0,1)} \leq2C_d\|\mu(x)\|_{L_p(0,1)}=2C_d\|x\|_p.
\end{align*}
The power of $1-\Delta_\theta$ here is negative, namely $-\frac{d}{4}$.

(ii) Proposition~\ref{ThmTBdd}(ii) states that
$\mathcal{T}:L_2(0,1)\to\Lambda_\psi^{(2)}(0,1)$ is bounded. The target is fully symmetric by \cite[Chapter~6, Theorem~6.3.8]{DoddsDePagterSukochev}.
Apply Theorem~\ref{thm:fractional-symmetric-embedding} with
$E=L_2$ and $F=\Lambda_\psi^{(2)}$. This yields
\begin{align*}
\|J_\theta^{-\frac{d}{2}}x\|_{\Lambda_\psi^{(2)}} \leq C_d\|\mathcal{T}\mu(x)\|_{\Lambda_\psi^{(2)}(0,1)} \leq C_d\|\mathcal{T}\|_{L_2\to\Lambda_\psi^{(2)}}\|\mu(x)\|_2 = C_d'\|x\|_2.
\end{align*}
For $y\in H_2^{\frac{d}{2}}$, apply this estimate to $x=J_\theta^{\frac{d}{2}}y$:
\begin{equation*}
\|y\|_{\Lambda_\psi^{(2)}}
=\|J_\theta^{-\frac{d}{2}}J_\theta^{\frac{d}{2}}y\|_{\Lambda_\psi^{(2)}}
\leq C_d'\|y\|_{H_2^{\frac{d}{2}}}.
\end{equation*}
This proves the final embedding.
\end{proof}

\section{Evolution equations}\label{sec:evolution}
Finally, we apply our estimates to evolution equations on quantum tori. We first study linear fractional diffusion using spectral calculus and Mittag--Leffler functions \cite{Podlubny}. We then consider semilinear heat equations, building on the classical heat-semigroup and Besov-space methods \cite{DaviesHeat,BahouriCheminDanchin}.

The Besov algebra property gives local well-posedness for sufficiently regular initial data. Heat smoothing allows us to treat less regular data below the algebra threshold. Our results also extend the quantum-torus setting considered in~\cite{YY}.

\subsection{Linear fractional diffusion}

We use the two-parameter Mittag--Leffler function \cite{Podlubny}
\begin{equation*}
E_{\alpha,\beta}(z)=\sum_{j=0}^\infty
\frac{z^j}{\Gamma(\alpha j+\beta)},\ \ 
\alpha>0, \ \ \beta\in\mathbb{C},
\end{equation*}
and abbreviate $E_{\alpha,1}$ to $E_\alpha$. The series is locally
uniformly convergent. The negative-real-axis estimate needed below is
\begin{equation}\label{ML-estimate}
|E_\alpha(-r)|\leq\frac{C_\alpha}{1+r},\ \ 
r\geq0, \ \  0<\alpha\leq1.
\end{equation}
For $0<\alpha<1$ this follows from
\cite[Chapter~1, Theorem~1.4]{Podlubny}; for $\alpha=1$, use
$E_1(-r)=e^{-r}$.

For $0<\alpha<1$, write the Riemann--Liouville derivative as
\begin{equation*}
D_t^\alpha v(t)=\frac{d}{dt}
\frac1{\Gamma(1-\alpha)}\int\limits_0^t(t-s)^{-\alpha}v(s)\,ds,
\end{equation*}
whenever this expression is defined, or in the distributional sense.
For an absolutely continuous $u$ the Caputo derivative is
\begin{equation*}
{}^CD_t^\alpha u(t)=D_t^\alpha(u-u(0))(t)
=\frac1{\Gamma(1-\alpha)}\int\limits_0^t(t-s)^{-\alpha}u'(s)\,ds.
\end{equation*}
Accordingly, the notation $\partial_t^\alpha(u-u_0)$ in the diffusion
equation below denotes $D_t^\alpha(u-u_0)$, the Caputo formulation
with initial value $u_0$. For $\alpha=1$ it is the ordinary derivative.
This convention also applies to Banach-valued functions; for mild
solutions the spectral formula supplies the definition.

\subsubsection{Solution formula and decay}\label{subsec:linear-fractional-diffusion} Consider the following heat-type equation for $s>0$
\begin{equation}\label{Main-equation-1}
\partial_{t}^{\alpha}(u(t)-u_{0}) + (1-\Delta_{\theta})^{\frac{s}{2}}u(t) = 0, \ \ t > 0, \ \ 0 < \alpha \leq 1, 
\end{equation}
with the initial data
\begin{equation}\label{Initial-date-1}
u(0) = u_{0}\in L_{2}(\mathbb{T}^{d}_{\theta}).   
\end{equation} 

\begin{definition} Let $s>0$ and $u_{0}\in L_{2} (\mathbb{T}^{d}_{\theta}).$ Then the operator $u$ defined by the formula
\begin{equation*}
    u(t)=E_{\alpha}(-t^{\alpha}(1-\Delta_{\theta})^{\frac{s}{2}})u_0:=\sum\limits_{n\in\mathbb{Z}^{d}}E_{\alpha}(-t^{\alpha}(1+|n|^{2})^{\frac{s}{2}})\widehat{u}_{0}(n)e^{\theta}_{n}, \ t>0,
\end{equation*} 
is called a mild solution to the problem \eqref{Main-equation-1}-\eqref{Initial-date-1}. When $\theta=0$, the same definition is used in \cite{VergaraPresentAuthors,VergaraZacher2010}.

\end{definition}

Related $L_p$--$L_q$ decay estimates were studied in \cite{STT}.
Estimate~\eqref{ML-estimate} and Parseval's identity show that the
series defining $u(t)$ converges in $L_2$. Dominated convergence,
together with $E_\alpha(0)=1$, gives
$u\in C([0,\infty);L_2)$ and $u(0)=u_0$.
\begin{thm}\label{L2timeest}\label{prop:fractional-diffusion-decay}Let $s=\frac{d}{2}$, $0<\alpha\leq1$ and $u_{0} \in L_{2}(\mathbb{T}^{d}_{\theta}).$ Assume that $u$ is a mild solution to the Cauchy problem \eqref{Main-equation-1}-\eqref{Initial-date-1}. Then we have the following $L_{2}$-time decay for the mild solution
\begin{equation*}
    \|u(t)\|_{L_{2}(\mathbb{T}^{d}_{\theta})}\leq C_{d,\alpha}t^{-\alpha}\|u_0\|_{L_2(\mathbb{T}_{\theta}^{d})},\ \  t>0.
\end{equation*}
\end{thm}
\begin{proof}
Set $\lambda_n=(1+|n|^2)^{\frac{d}{4}}\geq1$. By \eqref{ML-estimate}
and Plancherel,
\begin{equation*}
 \|u(t)\|_2^2\leq C_\alpha^2\sum_n
 \frac{|\widehat u_0(n)|^2}{(1+t^\alpha\lambda_n)^2}
 \leq \frac{C_\alpha^2}{(1+t^\alpha)^2}\|u_0\|_2^2.
\end{equation*}
Since $(1+t^\alpha)^{-1}\leq t^{-\alpha}$, the asserted estimate follows.
In fact, the same argument applies to every $s>0$ in
\eqref{Main-equation-1}. For $\alpha=1$, the spectral formula also
gives $\|u(t)\|_2\leq e^{-t}\|u_0\|_2$.
\end{proof}

\begin{rem}[Connection with the symmetric Sobolev estimate]
The critical choice $s=\frac{d}{2}$ gives the stronger distributional bound
\begin{equation*}
 \mu(u(t))\prec\!\prec
 C_d\,\mathcal{T}\mu(J_\theta^{\frac{d}{2}}u(t)),\ \ 
 \|J_\theta^{\frac{d}{2}}u(t)\|_2\leq C_\alpha t^{-\alpha}\|u_0\|_2.
\end{equation*}
Indeed, apply Theorem~\ref{thm:critical-fractional-distribution}
to $J_\theta^{\frac{d}{2}}u(t)$ and use \eqref{ML-estimate}.
Proposition~\ref{ThmTBdd}(ii) therefore yields
\begin{equation*}
 \|u(t)\|_{\Lambda_\psi^{(2)}(\mathbb{T}_\theta^{d})}
 \leq C_{d,\alpha}t^{-\alpha}\|u_0\|_2.
\end{equation*}
This last conclusion uses the fractional embedding theory; the
elementary $L_2$ decay itself follows directly from the spectrum.
\end{rem}

\subsection{Semilinear heat equations on the symmetric Besov scale}
\label{sec:pde-application}

Throughout this subsection, let $E$ be fully symmetric, $s>0$, and $1\leq r<\infty$. Polynomial density in $B_{\mathcal{E},r}^{s}(\mathbb{T}_\theta^d)$ gives strong continuity of the heat semigroup at zero, without requiring separability of $E$. We use the endpoint convolution bounds and the product estimates of Section~\ref{sec:symmetric-besov-spaces}.

\subsubsection{The scalar heat kernel}

Let $H_t$ be the nonnegative classical heat kernel on $\mathbb{T}^d$, normalized by
\begin{equation*}
 \widehat H_t(n)=e^{-t|n|^2},
 \ \  \|H_t\|_1=1.
\end{equation*}
Then $P_{t,\theta}=e^{t\Delta_\theta}$ is given by $P_{t,\theta}x=H_t\ast x$. The classical Gaussian bounds, obtained by periodizing the Euclidean heat kernel, imply that
\begin{equation}\label{eq:prelim-heat-kernel}
 \|H_t\|_{L_p(\mathbb{T}^d)}
 \leq C_{d,p}t^{-\frac d2(1-\frac1p)},
 \ \  0<t\leq1,\quad 1\leq p\leq\infty;
\end{equation}
and
\begin{equation*}
 \|\partial^\alpha H_t\|_{L_p(\mathbb{T}^d)}
 \leq C_{d,p,\alpha}
 t^{-\frac{|\alpha|}{2}-\frac d2(1-\frac1p)},
 \ \  0<t\leq1.
\end{equation*}
See \cite{DaviesHeat} for general theory on heat-kernel estimates.

\subsubsection{Heat estimates}
For fixed $\theta$, write $P_t=P_{t,\theta}$. We combine the scalar heat-kernel estimates with the convolution bounds to obtain estimates on $B_{\mathcal E,r}^{s}(\mathbb{T}_\theta^d)$.

\begin{lem}\label{lem:besov-app-heat}
The heat operators form a strongly continuous contraction semigroup
on $B_{\mathcal{E},r}^{s}(\mathbb{T}_{\theta}^{d})$. For $\beta\geq0$ and $0<t\leq1$,
\begin{equation}\label{eq:besov-app-smoothing}
 \|P_{t}x\|_{B_{\mathcal{E},r}^{s+\beta}} \leq C_\beta t^{-\frac{\beta}{2}}\|x\|_{B_{\mathcal{E},r}^{s}}.
\end{equation}
For some fixed $c>0$, define
\begin{equation*}
 \mathcal H_{E,s,r}(t)=\phi_E(1)^{-1}
 +\left\|\left(\omega_{j}(E,s)e^{-ct2^{2j}}\right)_{j\geq1}
   \right\|_{\ell_{r'}}.
\end{equation*}
Then
\begin{equation}\label{eq:besov-app-profile}
 \|P_{t}x\|_\infty\leq C\mathcal H_{E,s,r}(t)\|x\|_{B_{\mathcal{E},r}^{s}}.
\end{equation}
If $\phi_E(v)\geq c_Ev^{\frac{1}{p}}$ for $0<v\leq1$, $1\leq p<\infty$, then
\begin{equation}\label{eq:besov-app-power-profile}
 \mathcal H_{E,s,r}(t)\leq C
 \begin{cases}
 t^{-\frac{1}{2}(\frac{d}{p}-s)},&s<\frac{d}{p},\\
 (1+\log\left(1/t\right))^{\frac{1}{r'}},&s=\frac{d}{p},\\
 1,&s>\frac{d}{p}.
 \end{cases}
\end{equation}
\end{lem}
\begin{proof}
 Since $P_{t}$ is a contractive convolution operator on $\mathcal{E}$ and commutes with $\Delta_j$,
\begin{equation*}
 \|\Delta_jP_{t}x\|_{\mathcal{E}} = \|P_{t}\Delta_jx\|_{\mathcal{E}} \leq \|\Delta_jx\|_{\mathcal{E}}.
\end{equation*}
Taking the weighted $\ell_r$ norm gives $\|P_{t}x\|_{B_{\mathcal{E},r}^{s}}\leq\|x\|_{B_{\mathcal{E},r}^{s}}$. For $r<\infty$, choose a Fourier polynomial $p$ approximating $x$ in $B_{\mathcal{E},r}^{s}$. Then
\begin{equation*}
 \|P_{t}x-x\|_{B_{\mathcal{E},r}^{s}} \leq 2\|x-p\|_{B_{\mathcal{E},r}^{s}}+\|P_{t}p-p\|_{B_{\mathcal{E},r}^{s}}  \to  2\|x-p\|_{B_{\mathcal{E},r}^{s}}
 \ \ (t\downarrow0).
\end{equation*}
Polynomial density therefore gives strong continuity.

We next prove the localized estimate. Choose $\widetilde\varphi\in C_c^\infty(\mathbb{R}^{d}\setminus\{0\})$ equal to $1$ on the support of the annular cutoff defining $\Delta_j$, after rescaling. Set
\begin{equation*}
 v=t2^{2j},  \ \ 
 m_v(\xi)=e^{-v|\xi|^2}\widetilde\varphi(\xi).
\end{equation*}
A fixed normalization constant in the heat symbol can be absorbed into $v$. Since $|\xi|\geq a>0$ on $\operatorname{supp}\widetilde\varphi$, for every multi-index $\alpha$,
\begin{equation*}
 \|\partial^\alpha m_v\|_\infty
 \leq C_\alpha(1+v)^{|\alpha|}e^{-a^2v}
 \leq C_\alpha'e^{-cv},
 \ \  v\geq0.
\end{equation*}
Choose an integer $M>\frac{d}{2}$. Integration by parts yields
\begin{equation*}
 |\check m_v(y)| \leq C(1+|y|^2)^{-M} \|(1-\Delta_\xi)^Mm_v\|_{L_1(\mathbb{R}^{d})}
 \leq Ce^{-cv}(1+|y|^2)^{-M},
\end{equation*}
and hence
\begin{equation*}
 \|\check m_v\|_{L_1(\mathbb{R}^{d})}\leq Ce^{-cv}.
\end{equation*}
The periodized kernel
\begin{equation*}
 K_{t,j}(u) = \sum_{k\in\mathbb Z^{d}} 2^{jd}\check m_v\bigl(2^j(u+k)\bigr)
\end{equation*}
therefore satisfies
\begin{equation*}
 \|K_{t,j}\|_{L_1(\mathbb{T}^{d})}
 \leq \|\check m_v\|_{L_1(\mathbb{R}^{d})} \leq Ce^{-ct2^{2j}}.
\end{equation*}
Its Fourier multiplier agrees with that of $P_{t}$ on the support of $\Delta_j$. The convolution-module estimate gives
\begin{equation*}
 \|\Delta_jP_{t}x\|_{\mathcal{E}} \leq Ce^{-ct2^{2j}}\|\Delta_jx\|_{\mathcal{E}}, \ \ j\geq1.
\end{equation*}

Consequently, for $\beta\geq0$,
\begin{equation*}
 \begin{aligned}
 \left\|
 \bigl(2^{j(s+\beta)} \|\Delta_jP_{t}x\|_{\mathcal{E}}\bigr)_{j\geq1} \right\|_{\ell_r}& 
 \leq C\sup_{j\geq1} \bigl(2^{j\beta}e^{-ct2^{2j}}\bigr)  \left\| \bigl(2^{js}\|\Delta_jx\|_{\mathcal{E}}\bigr)_{j\geq1} \right\|_{\ell_r}\\
 & \leq C_\beta t^{-\frac{\beta}{2}}\|x\|_{B_{\mathcal{E},r}^{s}}.
 \end{aligned}
\end{equation*}
Together with contractivity of the low block, this proves \eqref{eq:besov-app-smoothing}.
\end{proof}

\subsubsection{Polynomial nonlinearities and local solutions}
Consider an ordered polynomial
\begin{equation}\label{eq:besov-app-polynomial}
 \mathcal{P}(u,u^*)=a_0+
 \sum_{\nu=1}^{N}\sum_{\lambda=1}^{M_\nu}
 a_{\nu,\lambda,0}u^{\varepsilon_{\nu,\lambda,1}}
 a_{\nu,\lambda,1}\cdots
 u^{\varepsilon_{\nu,\lambda,\nu}}a_{\nu,\lambda,\nu},
\end{equation}
where $u^1=u$, $\varepsilon_{\nu,\lambda,j}\in\left\{1,*\right\}$, and the coefficients are fixed operators. We study
\begin{equation}\label{eq:besov-app-equation}
 \partial_tu-\Delta_\theta u=\mathcal{P}(u,u^*),\ \  u(0)=u_0.
\end{equation}
A mild solution satisfies
\begin{equation}\label{eq:besov-app-mild}
 u(t)=P_{t} u_0 + \int\limits_0^tP_{t-\tau}\mathcal{P}(u(\tau),u(\tau)^*)\,d\tau,
\end{equation}
with a Bochner integral in $B_{\mathcal{E},r}^{s}(\mathbb{T}_\theta^{d})$.

\begin{thm}[Local solutions in the algebra range]
\label{thm:besov-app-algebra-pde}
Assume $\omega(E,s)\in\ell_{r'}$ and all coefficients in \eqref{eq:besov-app-polynomial} belong to $B_{\mathcal{E},r}^{s}$. Every $u_0\in B_{\mathcal{E},r}^{s}(\mathbb{T}_\theta^{d})$ has a unique maximal mild solution $u\in C([0,T_{\max});B_{\mathcal{E},r}^{s})$. Existence times are uniform on bounded sets of initial data, and the solution map is locally Lipschitz. If $T_{\max}<\infty$, then
\begin{equation*}
 \lim_{t\uparrow T_{\max}}\|u(t)\|_{B_{\mathcal{E},r}^{s}}=\infty.
\end{equation*}
For each $0<\beta<2$, the solution belongs to $C((0,T_{\max});B_{\mathcal{E},r}^{s+\beta})$.
\end{thm}
\begin{proof}

By the algebra estimate and boundedness of the involution,
\begin{equation*}
 \|xy\|_{B_{\mathcal{E},r}^{s}}\leq C_A\|x\|_{B_{\mathcal{E},r}^{s}}\|y\|_{B_{\mathcal{E},r}^{s}}, \ \
 \|x^*\|_{B_{\mathcal{E},r}^{s}}\leq C_*\|x\|_{B_{\mathcal{E},r}^{s}}.
\end{equation*}
Each nonconstant term of $\mathcal P(u,u^*)$ has the form
\begin{equation*}
 W(u)=a_0u^{\varepsilon_1}a_1 \cdots  u^{\varepsilon_m}a_m,
 \ \ m\geq1,\ \ \varepsilon_k\in\{1,*\},
\end{equation*}
where $u^1=u$ and $a_0,\ldots,a_m$ are fixed coefficients. Successive replacement of the unknown factors gives
\begin{equation*}
 W(u)-W(v) = \sum_{k=1}^m a_0
 \left(\prod_{i=1}^{k-1}u^{\varepsilon_i}a_i\right)
 (u^{\varepsilon_k}-v^{\varepsilon_k})a_k
 \left(\prod_{i=k+1}^{m}v^{\varepsilon_i}a_i\right),
\end{equation*}
where products are taken in increasing order. Consequently, for $\|u\|_{B_{\mathcal{E},r}^{s}},\|v\|_{B_{\mathcal{E},r}^{s}}\leq R$,
\begin{equation*}
 \|W(u)\|_{B_{\mathcal{E},r}^{s}}\leq C_W R^m,
 \ \   \|W(u)-W(v)\|_{B_{\mathcal{E},r}^{s}} \leq mC_WR^{m-1}\|u-v\|_{B_{\mathcal{E},r}^{s}}.
\end{equation*}
Summing over the terms of $\mathcal{P}$, including its constant term, yields
\begin{equation*}
 \|\mathcal{P}(u,u^*)\|_{B_{\mathcal{E},r}^{s}}\leq M_R,
 \ \ 
 \|\mathcal{P}(u,u^*)-\mathcal{P}(v,v^*)\|_{B_{\mathcal{E},r}^{s}}
 \leq L_R\|u-v\|_{B_{\mathcal{E},r}^{s}}.
\end{equation*}

Let $Y_T=C([0,T];X)$ and define
\begin{equation*}
 (\Phi u)(t)
 =P_tu_0+\int\limits_0^tP_{t-\tau}\mathcal{P}(u(\tau),u(\tau)^*)\,d\tau.
\end{equation*}
Strong continuity of $P_t$ on $B_{\mathcal{E},r}^{s}$ and continuity of the nonlinearity imply $\Phi u\in Y_T$. On the closed ball $\mathbb B_R=\{u\in Y_T:\|u\|_{Y_T}\leq R\}$, contractivity gives
\begin{equation*}
 \|\Phi u\|_{Y_T}\leq\|u_0\|_{B_{\mathcal{E},r}^{s}}+TM_R,
 \ \ 
 \|\Phi u-\Phi v\|_{Y_T}
 \leq TL_R\|u-v\|_{Y_T}.
\end{equation*}
Choose $R>2\|u_0\|_{B_{\mathcal{E},r}^{s}}+1$ and $0<T\leq1$ such that
\begin{equation*}
 TM_R\leq \frac{R}{2}, \ \ TL_R\leq\frac{1}{2}.
\end{equation*}
The contraction principle gives a mild solution in $\mathbb B_R$. Taking $R=2K+2$ shows that the same $T=T_K>0$ works for all $\|u_0\|_{B_{\mathcal{E},r}^{s}}\leq K$.

If $u,v$ are two mild solutions bounded by $R$ on $[0,T]$, then
\begin{equation*}
 \|u(t)-v(t)\|_{B_{\mathcal{E},r}^{s}}
 \leq \|u_0-v_0\|_{B_{\mathcal{E},r}^{s}}
      +L_R\int\limits_0^t\|u(\tau)-v(\tau)\|_{B_{\mathcal{E},r}^{s}}\,d\tau.
\end{equation*}
Gronwall's inequality yields
\begin{equation*}
 \sup_{0\leq t\leq T}\|u(t)-v(t)\|_{B_{\mathcal{E},r}^{s}}
 \leq e^{L_RT}\|u_0-v_0\|_{B_{\mathcal{E},r}^{s}}.
\end{equation*}
This proves uniqueness and local Lipschitz dependence on the common local existence interval. Repeated application gives the same conclusion on compact subintervals of a solution's lifespan.

The semigroup property gives the restart identity
\begin{equation*}
 u(t)=P_{t-t_0}u(t_0) + \int\limits_{t_0}^tP_{t-\tau} \mathcal{P}(u(\tau),u(\tau)^*)\,d\tau,  \ \  t\geq t_0.
\end{equation*}
Local existence and uniqueness therefore define a unique
maximal solution on $[0,T_{\max})$. If $T_{\max}<\infty$ and the asserted limit failed, there would exist $K<\infty$ and $t_n\uparrow T_{\max}$ such that $\|u(t_n)\|_{B_{\mathcal{E},r}^{s}}\leq K$. Restarting at $t_n$ gives a solution on $[t_n,t_n+T_K]$, with $T_K$ independent of $n$. For large $n$, this extends $u$ beyond $T_{\max}$, a contradiction.

Finally, fix $0<T<T_{\max}$ and $0<\beta<2$, and put
\begin{equation*} 
F(t)=\mathcal{P}(u(t),u(t)^*),  \ \  R=\sup_{0\leq t\leq T}\|u(t)\|_{B_{\mathcal{E},r}^{s}}.
\end{equation*}
The smoothing estimate implies, for $0<t\leq T$,
\begin{equation*}
 \begin{aligned}
 \|u(t)\|_{B_{\mathcal{E},r}^{s+\beta}(\mathbb{T}_\theta^d)}
 &\leq C_Tt^{-\frac{\beta}{2}}\|u_0\|_{B_{\mathcal{E},r}^{s}}
       +C_TM_R\int\limits_0^t(t-\tau)^{-\frac{\beta}{2}}\,d\tau\\
 &\leq C_Tt^{-\frac{\beta}{2}}\|u_0\|_{B_{\mathcal{E},r}^{s}}
       +\frac{C_TM_R}{1-\frac{\beta}{2}}\,t^{1-\frac{\beta}{2}}.
 \end{aligned}
\end{equation*}
To verify continuity, fix $0<a<T$ and $0<\varepsilon<a$. For $t\in[a,T]$, decompose the Duhamel term as
\begin{equation*}
 \int\limits_0^tP_{t-\tau}F(\tau)\,d\tau = \underbrace{\int\limits_0^{t-\varepsilon}  P_{t-\tau}F(\tau)\,d\tau}_{I_\varepsilon(t)}
 + \underbrace{\int\limits_{t-\varepsilon}^{t} P_{t-\tau}F(\tau)\,d\tau}_{J_\varepsilon(t)}.
\end{equation*}
Operator-norm continuity of $P_h:B_{\mathcal{E},r}^{s}\to B_{\mathcal{E},r}^{s+\beta}$ for $h>0$ implies $I_\varepsilon\in C([a,T];B_{\mathcal{E},r}^{s+\beta})$, whereas
\begin{equation*}
 \sup_{a\leq t\leq T}\|J_\varepsilon(t)\|_{B_{\mathcal{E},r}^{s+\beta}(\mathbb{T}_\theta^d)}
 \leq  \frac{C_TM_R}{1-\frac{\beta}{2}}\, \varepsilon^{1-\frac{\beta}{2}}  \to 0.
\end{equation*}
Thus the Duhamel term is continuous into $B_{\mathcal{E},r}^{s+\beta}$ on $[a,T]$. The same holds for $t\mapsto P_tu_0$. Since $a$ and $T$ are arbitrary, $u\in C((0,T_{\max});B_{\mathcal{E},r}^{s+\beta})$.
\end{proof}
The next result does not require $B_{\mathcal{E},r}^{s}$ to embed into $L_\infty$. It uses the smoothing of the heat semigroup to control multiplication at positive times. Outside the algebra range, the polynomial nonlinearity need not define a locally Lipschitz map from $X$ into itself. We use the fact that $P_t$ maps $B_{\mathcal{E},r}^{s}$ into $L_\infty$ for every $t>0$.

\begin{thm}[Local solutions under a heat regularization estimate]
\label{thm:besov-app-weighted}
Let $N\geq2$, and assume that all coefficients in \eqref{eq:besov-app-polynomial} belong to
$B_{\infty,\infty}^{\rho}$ for some $\rho>s$. Suppose that
\begin{equation}\label{eq:besov-app-weight-assumption}
 \|P_tx\|_\infty\leq C_Ht^{-a}\|x\|_{B_{\mathcal{E},r}^{s}},
 \ \ 
 0<a<\frac1{N-1},
 \ \  0<t\leq1.
\end{equation}
For every $u_0\in X$, there exists a local mild solution in
\begin{equation*}
 \mathcal{Z}_T^a = 
 \left\{ u\in C([0,T];B_{\mathcal{E}.r}^{s})\cap C((0,T];L_\infty): \sup_{0<t\leq T}t^a\|u(t)\|_\infty<\infty  \right\},
\end{equation*}
equipped with the norm
\begin{equation*}
 \|u\|_{\mathcal{Z}_T^a}
 = \sup_{0\leq t\leq T}\|u(t)\|_{B_{\mathcal{E},r}^{s}}
 + \sup_{0<t\leq T}t^a\|u(t)\|_\infty.
\end{equation*}
The solution is unique in this class, depends locally Lipschitz continuously on $u_0$, and satisfies
\begin{equation*}
 \lim_{t\downarrow0}t^a\|u(t)\|_\infty=0.
\end{equation*}
The local existence time is uniform on bounded subsets of $B_{\mathcal{E},r}^{s}$. The solution extends uniquely to a maximal interval $[0,T_{\max})$, and
\begin{equation*}
 T_{\max}<\infty \Rightarrow 
 \lim_{t\uparrow T_{\max}}\|u(t)\|_{B_{\mathcal{E},r}^{s}}=\infty.
\end{equation*}
\end{thm}
\begin{proof}
Put $b=(N-1)a<1$. Since $N\geq2$, we also have $a<1$. The space $\mathcal{Z}_T^a$ is Banach. Indeed, a Cauchy sequence converges in $C([0,T];B_{\mathcal{E},r}^{s})$ and in the weighted $L_\infty$ norm. For each $t>0$, these limits agree in $L_1$, because
\begin{equation*}
 B_{\mathcal{E},r}^{s}\hookrightarrow\mathcal E\hookrightarrow L_1,
 \ \  L_\infty\hookrightarrow L_1.
\end{equation*}
On every closed interval contained in $(0,T]$, convergence is uniform in $L_\infty$, so the limit is continuous there. By contractivity and \eqref{eq:besov-app-weight-assumption},
\begin{equation*}
 \|P_\cdot u_0\|_{\mathcal{Z}_T^a}
 \leq (1+C_H)\|u_0\|_{B_{\mathcal{E},r}^{s}},
 \ \  0<T\leq1.
\end{equation*}

For a term of the form
\begin{equation*}
 W(z_1,\ldots,z_\nu)=a_0z_1a_1\cdots z_\nu a_\nu,
\end{equation*}
repeated application of \eqref{eq:besov-app-tame} and \eqref{eq:besov-app-coefficient} gives
\begin{equation*}
 \|W(z_1,\ldots,z_\nu)\|_{B_{\mathcal{E},r}^{s}}
 \leq C_W\sum_{j=1}^{\nu}
 \|z_j\|_{B_{\mathcal{E},r}^{s}}
 \prod_{k\ne j}\|z_k\|_\infty.
\end{equation*}
The assumptions on the coefficients also imply that the constant term of $\mathcal P$ belongs to $B_{\mathcal{E},r}^{s}$.

Let $\|u\|_{\mathcal{Z}_T^a},\|v\|_{\mathcal{Z}_T^a}\leq R$, where $R\geq1$, and put $d=\|u-v\|_{\mathcal{Z}_T^a}$. For $0<\tau\leq T\leq1$, we have
\begin{equation*}
 \|u(\tau)\|_{B_{\mathcal{E},r}^{s}},
 \|v(\tau)\|_{B_{\mathcal{E},r}^{s}}\leq R,
 \ \ 
 \|u(\tau)\|_\infty,\|v(\tau)\|_\infty\leq R\tau^{-a},
\end{equation*}
and
\begin{equation*}
 \|u(\tau)-v(\tau)\|_{B_{\mathcal{E},r}^{s}}\leq d,
 \ \ 
 \|u(\tau)-v(\tau)\|_\infty\leq d\tau^{-a}.
\end{equation*}
The same bounds hold for adjoints, up to fixed constants.

To estimate the difference of two terms of degree $\nu$, replace the factors one at a time. Each resulting term contains exactly one difference factor. If this factor is measured in $B_{\mathcal{E},r}^{s}$, the product estimate gives
\begin{equation*}
 Cd(R\tau^{-a})^{\nu-1}.
\end{equation*}
If another factor is measured in $B_{\mathcal{E},r}^{s}$, it gives, for $\nu\geq2$,
\begin{equation*}
 CR(d\tau^{-a})(R\tau^{-a})^{\nu-2}.
\end{equation*}
Both expressions are bounded by $CR^{\nu-1}\tau^{-(\nu-1)a}d$. Since $\nu\leq N$, summing over the terms of $\mathcal P$ gives
\begin{align*}
 \|\mathcal P(u(\tau),u(\tau)^*)\|_{B_{\mathcal{E},r}^{s}}
 &\leq C(1+R^N\tau^{-b}),\\
 \|\mathcal P(u(\tau),u(\tau)^*)  -\mathcal{P}(v(\tau),v(\tau)^*)\|_{B_{\mathcal{E},r}^{s}}
 &\leq CR^{N-1}\tau^{-b}d.
\end{align*}

Write $F_u(\tau)=\mathcal P(u(\tau),u(\tau)^*)$ and
\begin{equation*}
 D_u(t)=\int\limits_0^tP_{t-\tau}F_u(\tau)\,d\tau.
\end{equation*}
Contractivity on $B_{\mathcal{E},r}^{s}$ and the heat estimate give
\begin{align*}
 \|D_u(t)\|_{B_{\mathcal{E},r}^{s}}
 &\leq C\left(t+\frac{R^N}{1-b}t^{1-b}\right),\\
 t^a\|D_u(t)\|_\infty
 &\leq C\left(
 \frac{t}{1-a}+R^N \mathrm{B}(1-a,1-b)t^{1-b}\right).
\end{align*}
Here $\mathrm{B}$ is the beta function, and we used
\begin{equation*}
 \int\limits_0^t\tau^{-b}\,d\tau=\frac{t^{1-b}}{1-b},
 \ \ 
 t^a\int\limits_0^t(t-\tau)^{-a}\tau^{-b}\,d\tau = \mathrm{B}(1-a,1-b)t^{1-b}.
\end{equation*}
Applying the same calculation to $F_u-F_v$, we obtain
\begin{align}
 \|D_u\|_{\mathcal{Z}_T^a}
 &\leq C(T+R^NT^{1-b}),\nonumber\\
 \|D_u-D_v\|_{\mathcal{Z}_T^a}
 &\leq CR^{N-1}T^{1-b}\|u-v\|_{\mathcal{Z}_T^a}.
 \label{eq:besov-app-contraction}
\end{align}

We next check that $D_u\in\mathcal{Z}_T^a$. The product estimate shows that $F_u:(0,T]\to B_{\mathcal{E},r}^{s}$ is continuous. For fixed $t>0$, the function $\tau\mapsto P_{t-\tau}F_u(\tau)$ is continuous on $(0,t)$ in both $B_{\mathcal{E},r}^{s}$ and $L_\infty$. The preceding bounds therefore ensure that its Bochner integral exists in both spaces. These integrals agree in $L_1$.

To prove continuity for $t>0$, fix $0<\eta<T$ and $0<\delta<\eta/2$. For $t\in[\eta,T]$, write
\begin{equation*}
 D_u(t) =  \int\limits_0^{t-\delta}P_{t-\tau}F_u(\tau)\,d\tau
 +  \int\limits_{t-\delta}^{t}P_{t-\tau}F_u(\tau)\,d\tau.
\end{equation*}
The first integral depends continuously on $t$ in both spaces, by operator-norm continuity of $P_h$ for $h>0$ and dominated convergence. Since $F_u$ is bounded on $[\eta/2,T]$, the second integral has norm at most $C_\eta\delta$ in $B_{\mathcal{E},r}^{s}$ and $C_\eta\delta^{1-a}$ in $L_\infty$, uniformly for $t\in[\eta,T]$. Letting $\delta\downarrow0$ proves continuity on $[\eta,T]$.

At zero, the preceding estimates give
\begin{equation*}
 \|D_u(t)\|_{B_{\mathcal{E},r}^{s}}
 +t^a\|D_u(t)\|_\infty \leq C(t+R^Nt^{1-b})\to 0
 \ \ (t\downarrow0).
\end{equation*}
Thus $D_u$ is continuous at zero in
$B_{\mathcal{E},r}^{s}$, with $D_u(0)=0$, and belongs to
$\mathcal{Z}_T^a$.

Define $\Phi u=P_\cdot u_0+D_u$. Choose
\begin{equation*}
 R>2(1+C_H)\|u_0\|_{B_{\mathcal{E},r}^{s}}+1
\end{equation*}
and then $0<T\leq1$ such that
\begin{equation*}
 C(T+R^NT^{1-b})\leq\frac{R}{2},
 \ \ 
 CR^{N-1}T^{1-b}\leq\frac12.
\end{equation*}
The map $\Phi$ sends the closed ball of radius $R$ in
$\mathcal{Z}_T^a$ into itself and is a contraction there.
Banach's fixed-point theorem gives a mild solution.

For solutions constructed in the same ball,
\begin{equation*}
 \|u-v\|_{\mathcal{Z}_T^a}
 \leq (1+C_H)\|u_0-v_0\|_{B_{\mathcal{E},r}^{s}}
 +\frac12\|u-v\|_{\mathcal{Z}_T^a}.
\end{equation*}
This proves local Lipschitz dependence.
The same $R$ and $T$ can be chosen for initial data in
any fixed bounded subset of $B_{\mathcal{E},r}^{s}$.

Any two solutions with the same initial data in the stated
class have finite weighted norms. Choose a ball containing
both and shorten the time interval so that the constant
in \eqref{eq:besov-app-contraction} is less than one.
They therefore agree near zero.
For every positive time $t_0$, the mild identity gives
\begin{equation*}
 u(t_0+h)=P_hu(t_0)
 +\int\limits_0^hP_{h-\sigma}
 \mathcal P(u(t_0+\sigma),u(t_0+\sigma)^*)\,d\sigma.
\end{equation*}
Since the solutions are continuous in
$B_{\mathcal{E},r}^{s}$ and $L_\infty$ near $t_0$,
the same argument applies with $t_0$ as the initial time.
This proves uniqueness throughout their common interval
of existence.

We can now extend the solution uniquely to a maximal
interval $[0,T_{\max})$.
If $T_{\max}<\infty$ and the norm did not tend to infinity,
there would be a sequence $t_n\uparrow T_{\max}$ such that
$\|u(t_n)\|_{B_{\mathcal{E},r}^{s}}\leq K$ for some $K$.
The uniform local existence time would allow us to extend
the solution from each $t_n$ by a fixed time $T_K>0$.
For large $n$, this extends it beyond $T_{\max}$, a
contradiction. Hence
\begin{equation*}
 \lim_{t\uparrow T_{\max}}
 \|u(t)\|_{B_{\mathcal{E},r}^{s}}=\infty.
\end{equation*}

Finally, for any Fourier polynomial $v$, the heat estimate
and contractivity on $L_\infty$ give
\begin{equation*}
 t^a\|P_tu_0\|_\infty
 \leq C_H\|u_0-v\|_{B_{\mathcal{E},r}^{s}}
 +t^a\|v\|_\infty.
\end{equation*}
Since $a>0$, the second term tends to zero.
Polynomial density in $B_{\mathcal{E},r}^{s}$ then gives
$t^a\|P_tu_0\|_\infty\to0$.
Together with the estimate for $D_u$, this proves
$t^a\|u(t)\|_\infty\to0$ as $t\downarrow0$.
\end{proof}

\begin{cor}[Positive-time regularity]
\label{cor:besov-app-positive-regularity}
Assume that the coefficients of $\mathcal P$ are smooth.
The solutions of
Theorems~\ref{thm:besov-app-algebra-pde}
and~\ref{thm:besov-app-weighted} satisfy
\begin{equation*}
 u\in C\bigl((0,T_{\max});
       B_{\mathcal{E},r}^{s+\gamma}(\mathbb{T}_\theta^d)\bigr)
 \ \ \text{for every }\gamma>0.
\end{equation*}
In particular, they are smooth in space for positive times,
belong to $C^1((0,T_{\max});X)$, and satisfy
\begin{equation*}
 \partial_tu-\Delta_\theta u=\mathcal P(u,u^*)
 \ \ \text{in }X,\quad t>0.
\end{equation*}
\end{cor}

\begin{proof}
Fix $0<t_0<T<T_{\max}$.
Both solutions are continuous into $X\cap L_\infty$ on
$[t_0,T]$: this follows from the algebra embedding in the
first case and from the weighted solution class in the second.
Thus $F(t)=\mathcal P(u(t),u(t)^*)$ is continuous into $X$ there.

Restarting at $t_0$, the smoothing estimate gives, for
$0<\beta<2$ and $t>t_0$,
\begin{equation*}
 \begin{aligned}
 \|u(t)\|_{B_{\mathcal{E},r}^{s+\beta}}
 & \leq C(t-t_0)^{-\beta/2}\|u(t_0)\|_{B_{\mathcal{E},r}^{s}}\\
 & + C\int\limits_{t_0}^t(t-\tau)^{-\beta/2}\|F(\tau)\|_{B_{\mathcal{E},r}^{s}}\,d\tau.
 \end{aligned}
\end{equation*}
The integral is finite, and splitting its terminal part as above
proves continuity in $B_{\mathcal{E},r}^{s+\beta}$ away from $t_0$.

At every increased positive order $\sigma$, the tame and coefficient estimates give
\begin{equation*}
 \|\mathcal P(u,u^*)\|_{B_{\mathcal{E},r}^{\sigma}}
 \leq C_\sigma
 \left(1+\|u\|_{B_{\mathcal{E},r}^{\sigma}}
                 (1+\|u\|_\infty)^{N-1}\right).
\end{equation*}
Their difference versions also give continuity of the forcing at that order. Smoothness of the coefficients allows this argument at every finite $\sigma$. Iteration on nested positive-time intervals therefore yields every gain $\gamma>0$.

Since $\mathcal{E}\hookrightarrow L_1$, arbitrarily high Besov regularity gives rapid decay of the Fourier coefficients and hence spatial smoothness. In particular, $u$ is continuous into the domain of the heat generator on $B_{\mathcal{E},r}^{\sigma}$, with $\Delta_\theta u$ continuous into $B_{\mathcal{E},r}^{\sigma}$. The mild identity then gives $u\in C^1((0,T_{\max});B_{\mathcal{E},r}^{\sigma})$ and the differential equation.
\end{proof}

\begin{cor}[Local solutions below the algebra threshold]
\label{cor:besov-app-power}
Suppose that
\begin{equation*}
 \phi_E(v)\geq c_Ev^{\frac{1}{p}},
 \ \  0<v\leq1,\ \  1\leq p<\infty.
\end{equation*}
For smooth coefficients and degree $N\geq2$, local well-posedness holds in $B_{\mathcal{E},r}^{s}$, $1\leq r<\infty$, whenever
\begin{equation}\label{eq:besov-app-subcritical}
 s>0,  \ \  s>\frac{d}{p}-\frac{2}{N-1}.
\end{equation}
Uniqueness is understood in the solution class of the applicable theorem above. The solutions have the corresponding maximal continuation and blow-up alternative, and are smooth for positive times.
\end{cor}

\begin{proof}
If $s<\frac{d}{p}$, the heat profile estimate gives
\begin{equation*}
 \|P_tx\|_\infty \leq Ct^{-(\frac{1}{2}\frac{d}{p}-s)}\|x\|_{B_{\mathcal{E},r}^{s}}.
\end{equation*}
Taking $a=\frac{1}{2}(\frac{d}{p}-s)$, condition \eqref{eq:besov-app-subcritical} is precisely $(N-1)a<1$, so Theorem~\ref{thm:besov-app-weighted} applies.

If $s=\frac{d}{p}$, the logarithmic bound in
\eqref{eq:besov-app-power-profile} is bounded by $C_at^{-a}$
for every $a>0$. Choose $a<1/(N-1)$ and apply the same theorem.
For $r=1$, this endpoint is also covered by the algebra theorem.

If $s>\frac{d}{p}$, then
\begin{equation*}
 \omega_j(E,s)\lesssim 2^{-j(s-\frac{d}{p})},
\end{equation*}
so $\omega(E,s)\in\ell_{r'}$ and
Theorem~\ref{thm:besov-app-algebra-pde} applies.
Positive-time regularity follows from
Corollary~\ref{cor:besov-app-positive-regularity}.
\end{proof}

\subsubsection{Range of the result and comparison}

For $E=L_2$, Plancherel's identity and finite overlap of the dyadic supports give
\begin{equation*}
B_{2,2}^{s}(\mathbb{T}_\theta^d)=H_\theta^{s}(\mathbb{T}_\theta^d)
\end{equation*}
with equivalent norms. Thus Theorem~\ref{thm:besov-app-algebra-pde} includes the local polynomial heat theory in $H_\theta^k$, $k\in\mathbb N$, $k>\frac{d}{2}$, considered by Yang and Yang~\cite[Theorem~1.2]{YY}. The Besov formulation also permits noninteger smoothness, fully symmetric bases, and the endpoint algebras $B_{p,1}^{\frac{d}{p}}$.

Corollary~\ref{cor:besov-app-power} additionally reaches below the algebra threshold. For example, when $d=3$, $E=L_2$, and $N=3$, initial data in $H_\theta^1=B_{2,2}^1$ are admissible:
\begin{equation*}
 a=\frac{3/2-1}{2}=\frac{1}{4},
 \ \  (N-1)a=\frac{1}{2}<1.
\end{equation*}
This applies, in particular, to
\begin{equation*}
 \partial_tu-\Delta_\theta u=u^3.
\end{equation*}
Such initial data need not be bounded operators, whereas the integer Sobolev algebra condition requires $k\geq2$. The corollary also applies to fully symmetric Banach realizations of Lorentz spaces with $1<p<\infty$, including weak-$L_p$ bases, whose fundamental functions are comparable to $v^{\frac{1}{p}}$.

The inequalities in \eqref{eq:besov-app-subcritical} are strict. At the critical equality $s=\frac{d}{p}-\frac{2}{N-1}>0$, the present weighted estimate gives the nonintegrable bound $\tau^{-1}$. At $s=0$, the summability argument used in the tame product estimate no longer applies. Neither boundary is treated here.

{\bf Conflict of Interest Statement:} The author declares that there is no conflict of interest.

{\bf Data Availability Statement:} No datasets were generated or analyzed during the current study.

{\bf Acknowledgement}

The authors thank Dr.~Edward McDonald for helpful discussions and valuable comments concerning Section~\ref{sec:evolution}. This work was partially supported by the Australian Research Council. F.S. acknowledges support from ARC grants FL17010005 and DP230100434. K.T. acknowledges partial support from ARC grant FL17010005.

An earlier version of part of this work was circulated as the preprint arXiv:2606.03185. The present manuscript has been substantially rewritten and reorganized, with additional results and applications. A substantial part of its mathematical content was already established in the earlier preprint.

During the revision, ChatGPT was used to assist with language, organization and the identification of possible references. We have not used any AI tools for developing the mathematical results of the paper. The authors are responsible for verifying all arguments and references and for the final content of the manuscript.

\end{document}